\documentclass[12pt,  reqno]{amsart}
\usepackage{amsmath,amssymb,amsfonts, amscd,   
pifont, amsthm,  amsxtra, latexsym, mathrsfs,   eufrak} 
\makeindex
\usepackage[dvips,hyperindex,breaklinks]{hyperref}
\hypersetup{linkbordercolor=0.2 1 0.5}
\newcommand*{\main}[1]{\textbf{\hyperpage{#1}}}
\hypersetup{linkbordercolor=0.2 1 0.5}
\usepackage{geometry, txfonts}
\usepackage{array}
\usepackage{epic,eepic}
\def\itbf{\itshape\bfseries }
\numberwithin{equation}{section} 
\theoremstyle{plain}
\newtheorem{thm}{Theorem}[section]
\newtheorem{coro}[thm]{Corollary}
\newtheorem{lem}[thm]{Lemma}

\newtheorem{prop}[thm]{Proposition}
\newtheorem{claim}[thm]{Claim}
\newtheorem{fact}[thm]{Fact}
\newtheorem{defe}[thm]{Definition}
\newtheorem{exa}[thm]{Example}
\newtheorem{rem}[thm]{Remark}

\newtheorem{thmalph}{Theorem}

\newtheorem{lemalph}[thmalph]{Lemma}
\newcommand{\cal}{\mathscr}
\def\Fka{\mathcal F_{k,a}}
\def\Hkm{\mathcal H_{k}^m}

\def\CC{\mathfrak{C}}

\def\a{\mathfrak{a}}
 
\def\g{\mathfrak{g}}

\def\l{\mathfrak{l}}

\def\p{\mathfrak{p}}
\def\k{\mathfrak{k}}

\def\s{\mathfrak{s}}
\def\u{\mathfrak{u}}
\def\oo{\mathfrak{o}}
\DeclareMathOperator{\Exp}{Exp}

\newcommand\C{{\mathbb{C}}}

\newcommand\N{{\mathbb{N}}}

\newcommand\Q{{\mathbb{Q}}}
\newcommand\R{{\mathbb{R}}}
\newcommand\Z{{\mathbb{Z}}}

\newcommand\HH{{\mathbb{H}}}
\newcommand\EE{{\mathbb{E}}}
\newcommand\DD{{\mathbb{D}}}

\def\sgn{\operatorname{sgn}}

\def\id{\operatorname{id}}

\def\C{{\mathbb C}}
\def\Q{{\mathbb Q}}
\def\Z{{\mathbb Z}}
\def\N{{\mathbb N}}

\def\0{{\bar 0}}
\def\1{{\bar 1}}

\def\End{{\operatorname{End}}}

\newcommand{\sh}{\operatorname{sinh}}

\def\bid{\hbox{\boldmath{$1$}}}

\def\phi{{\varphi}}

\renewcommand{\Re}{\operatorname{Re}}
\input xy
\xyoption{all}
\usepackage{enumerate}
\usepackage{float}
\usepackage{graphicx}
\renewcommand{\mathcal}{\mathscr}
\newcommand{\bigsearrow}{\rotatebox[origin=b]{-45}{$\xrightarrow{\kern10mm}$}}
\newcommand{\bigswarrow}{\rotatebox[origin=b]{45}{$\xleftarrow{\kern10mm}$}}
\newcommand{\bignearrow}{\rotatebox[origin=b]{-315}{$\xrightarrow{\kern10mm}$}}
\begin{document}
\title{Laguerre semigroup and Dunkl operators}
\author{Salem Ben Sa{\" \i}d}
\address{S. Ben Sa{\"\i}d:  Universit\'e Henri Poincar\'e-Nancy 1, 
Institut Elie Cartan, D\'epartement de Math\'ematiques, 
B.P. 239, 54506 Vandoeuvre-Les-Nancy, Cedex, France} 
\email{Salem.BenSaid@iecn.u-nancy.fr}
\author{Toshiyuki Kobayashi}
\address{T. Kobayashi: the University of Tokyo,
IPMU, Graduate School of Mathematical Sciences, 
3-8-1 Komaba, Meguro, Tokyo, 153-8914, Japan}
\email{toshi@ms.u-tokyo.ac.jp}
\author{Bent {\O}rsted}
\address{B. {\O}rsted: University of Aarhus, Department
of Mathematical Sciences,  Building 530, Ny Munkegade, DK 8000,
Aarhus C, Denmark} \email{orsted@imf.au.dk}
\begin{abstract} 
We construct a two-parameter family of actions 
$\omega_{k,a}$ of the Lie algebra $\s\l(2,\R)$ by differential-difference operators
on $\R^N \setminus \{0\}$.
Here, $k$ is a multiplicity-function for the Dunkl
operators, 
and $a>0$ arises from the interpolation of the two $\s\l(2,\R)$ actions on the Weil representation of
$Mp(N,\R)$ and the minimal unitary representation of $O(N+1,2)$.
We prove that this action $\omega_{k,a}$ lifts 
to a unitary representation of the universal covering of
$SL(2,\R)$, 
and can even be extended to a holomorphic semigroup $\Omega_{k,a}$.
In the $k\equiv 0$ case,
our semigroup generalizes the Hermite semigroup
studied by R. Howe ($a=2$) and the Laguerre semigroup by
the second author with G. Mano ($a=1$).
One  boundary value of our semigroup $\Omega_{k,a}$ provides us with 
\textit{$(k,a)$-generalized Fourier transforms} $\mathcal{F}_{k,a}$,
which includes the Dunkl transform $\mathcal{D}_k$
($a=2)$ and  a new  unitary operator $\mathcal{H}_k$
($a=1$), namely a Dunkl--Hankel transform.
We establish the inversion formula, 
and a generalization of the Plancherel theorem, 
the Hecke
identity, the Bochner identity, and a Heisenberg uncertainty
relation for $\mathcal{F}_{k,a}$.
We also find  kernel functions for $\Omega_{k,a}$ and $\Fka$
for $a=1,2$ in terms of Bessel functions and the
Dunkl intertwining operator.
\end{abstract}
\maketitle
\tableofcontents

\section{Introduction}

The classical Fourier transform is one of the most
basic objects in analysis; it may be understood as
belonging to a one-parameter group of unitary 
operators on $L^2(\R^N)$, and this group may even be extended 
holomorphically to a semigroup (the \textit{Hermite semigroup}) $I(z)$
generated by the self-adjoint operator 
$\Delta - \Vert x \Vert^2$. 
This is a holomorphic semigroup of bounded
operators depending on a complex variable $z$ in the
complex right half-plane, viz. $I(z + w) = I(z)I(w).$
The structure of this semigroup and its properties may
be appreciated without any reference to representation
theory, whereas the link itself is rich as was
revealed beautifully by R. Howe \cite{H2} in connection with the
Schr\"{o}dinger model of the Weil representation.

Our primary aim of this article is to give a foundation of  the deformation theory
 of the classical situation, by constructing 
a generalization $\Fka$ of the Fourier transform, 
and the holomorphic semigroup $\mathcal{I}_{k,a}(z)$ with
infinitesimal generator
$\Vert x \Vert^{2-a} \Delta_k - \Vert x \Vert^a$, 
acting on a
concrete Hilbert space deforming $L^2(\R^N)$. 
Here $\Delta_k$ is the Dunkl Laplacian (a differential-difference operator).
We analyze these operators $\Fka$ and $\mathcal{I}_{k,a}(z)$ in
the context of integral operators as well as representation theory. 

The deformation parameters in our setting consist of a real parameter $a$ 
coming from the interpolation of the minimal unitary representations
of two different reductive groups by keeping smaller symmetries (see {\sc Diagram 1.4}), 
and a parameter $k$ coming from Dunkl's theory of
differential-difference operators associated to a finite Coxeter group; 
also the dimension
$N$ and the complex variable $z$ 
may be considered as a parameter of the theory. 

We point out, that already deformations with $k = 0$ are new
and interpolate the minimal representations of two reductive groups
$O_0(n+1,2)\sptilde$ and $Mp(n,\R)$. Notice that these unitary
representations are generated by the `unitary inversion operator'
($=\cal F_{0,a}$ with $a=1,2,$ up to a scalar multiplication) together
with an elementary action of the maximal parabolic subgroups (see
\cite{KM-intopq} and \cite[Introduction]{{KM-intopq-ams}}). 

This article establishes the foundation of these new operators. 
Our theorems on $(k,a)$-generalized Fourier transforms $\cal F_{k,a}$  include:
\begin{itemize}
\item[--] Plancherel and inversion formula (Theorems \ref{thm:Fkauni} and \ref{thm:5.6}),
\item[--]  Bochner-type theorem (Theorem \ref{thm:5.15}),
\item[--] Heisenberg's uncertainty relation (Theorem \ref{thm:5.20}),
\item[--] exchange of multiplication and differentiation (Theorem \ref{thm:Fka}).
\end{itemize}
 
We think of the results and the methods  here as opening potentially interesting studies such as: 
\begin{itemize}
\item[--] characterization of `Schwartz space' and Paley--Wiener type theorem,
\item[--] Strichartz estimates for Schr\"odinger  and wave equations,
\item[--] Brownian motions in a Weyl chamber  (cf. \cite{GY}),
\item[--] analogues of Clifford analysis for the Dirac operator (cf. \cite{OSS}),
\end{itemize}
working with deformations of classical operators.

In the diagram
below we have summarized some of the deformation
properties by indicating the limit behaviour of
the holomorphic semigroup $\mathcal{I}_{k,a}(z)$; it is seen how various previous
integral transforms fit in our picture. In particular
we obtain as special cases the Dunkl transform $\mathcal{D}_k$ \cite{D3}
($a=2$, $z=\frac{\pi i}{2}$ and $k$ arbitrary), 
the Hermite semigroup $I(z)$ \cite{xfol,H2} ($a=2$, $k\equiv0$ and $z$
arbitrary), 
and the
Laguerre semigroup \cite{KM05,KM}
($a=1$, $k\equiv 0$ and $z$ arbitrary). 
Our framework gives a new treatment  even on the theory of the
Dunkl transform.

The `boundary value' of the holomorphic semigroup
$\mathcal{I}_{k,a}(z)$ from $\Re z>0$ to the imaginary axis gives rise
to a one-parameter subgroup of unitary operators.
The underlying idea may be interpreted as a descendent of Sato's
hyperfunction theory \cite{Sato} and also that of the
Gelfand--Gindikin program \cite{GG,HN2,Ols,Sta} for unitary
representations of real reductive groups.
The specialization
$\mathcal{I}_{k,a}(\frac{\pi i}{2})$ will be our
$(k,a)$-generalized Fourier transform $\Fka$ (up to a phase factor),
which reduces to the Fourier transform
($a=2$ and $k\equiv0$),
the Dunkl transform $\mathcal{D}_k$ ($a=2$ and
$k$ arbitrary),
and the Hankel transform ($a=1$ and $k\equiv0$).

Yet another specialization is to take $N=1$.
This very special case contains (after some change of variables) the results on 
 the $L^2$-model of the highest weight
representations of the universal covering group of
 $SL(2, \R)$, which was obtained  by B. Kostant
\cite{Kos} and R. Rao \cite{RR} by letting
$\mathfrak{sl}_2$ act as
differential operators on the half-line (see Remark \ref{rem:3.6}). 

The secondary aim of this article is to contribute to the
theory of special functions, in particular orthogonal
polynomials; indeed we derive several new identities,
for example, 
the $(k,a)$-deformation of the classical Hecke identity
(Corollary \ref{coro:5.14}) where the Gaussian function and harmonic
polynomials in the classical setting are replaced respectively with
$\exp(-\frac{1}{a}\Vert x\Vert^a)$ and 
 polynomials annihilated by the Dunkl Laplacian.
Another example is the identity \eqref{eqn:4.37}, 
which expresses 
 an infinite sum
of products of Bessel functions and Gegenbauer functions
as a single Bessel function.
 
In the rest of  the  Introduction we describe a little more
the contents of this article.

In Sections \ref{subsec:1.1} and \ref{subsec:1.2},
without any reference to representation theory,
we discuss our holomorphic semigroup
$\mathcal{I}_{k,a}(z)$ and $(k,a)$-generalized
Fourier transforms $\Fka$ as a two-parameter deformation of the
classical objects, i.e.\ the Hermite semigroup and the Euclidean
Fourier transform.

In Section \ref{subsec:1.3},
we introduce the basic machinery of the present article, namely,
to construct triples of differential-difference operators
generating the Lie algebra of $SL(2, \R)$, and 
 see how they are integrated
 to unitary representations of the universal covering
group.

One further aspect of our constructions is the link to
minimal unitary representations.
For the specific two parameters $(a,k)=(1,0)$ and $(2,0)$,
 we are really working with representations of
much larger semisimple groups, and our deformation is
interpolating the representation spaces for 
the minimal representations of two different groups. 
We highlight these hidden symmetries in Section \ref{subsec:1.4}.

 Let us also note that there is in  our theory a natural 
 appearance of some symmetries of the double degeneration
 of the double affine Hecke algebra (sometimes
 called the {\textit {rational Cherednik algebra}}),
 see Section \ref{subsec:5.6}. 
 Here $a=2$ and $k$ arbitrary, and in particular,
 we recover the Dunkl transform.

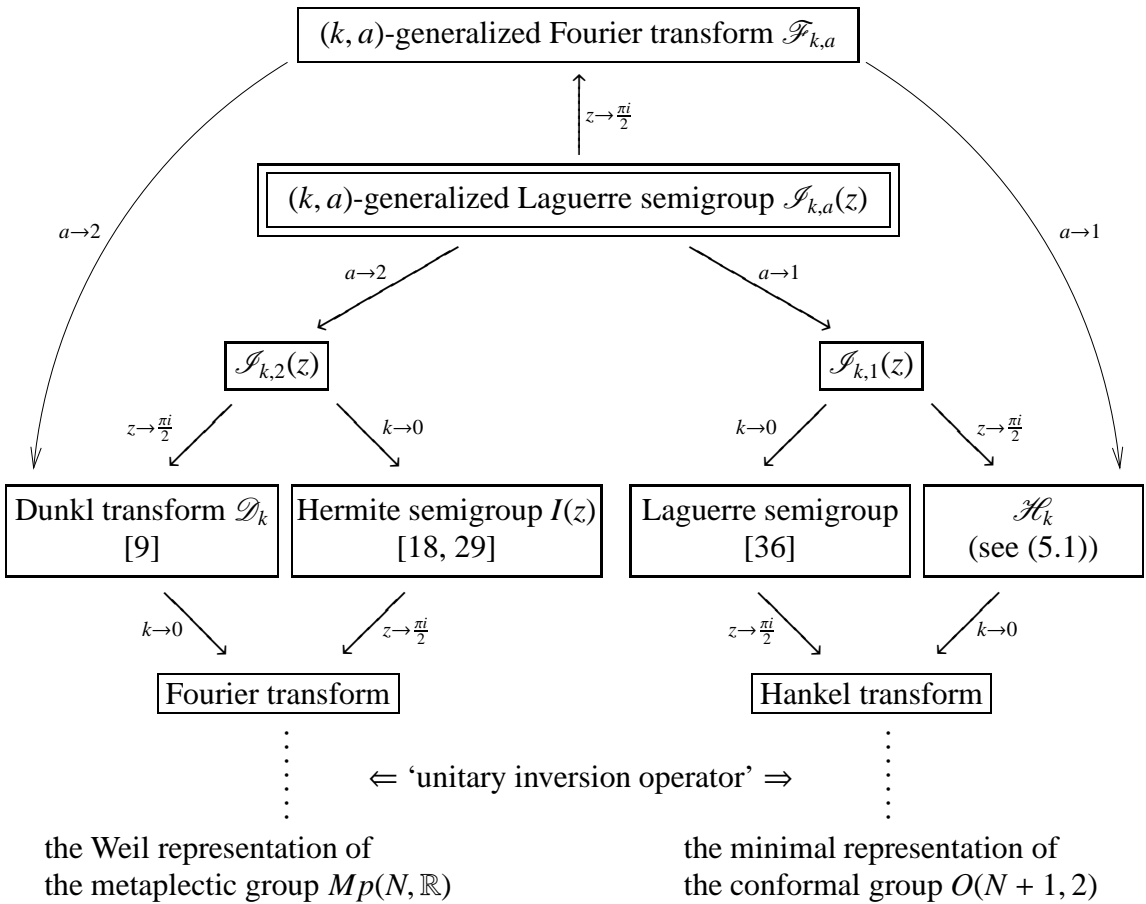
\begin{figure}[H]
\renewcommand{\figurename}{Diagram}
\setlength{\unitlength}{1mm}
\begin{picture}(160,120)
\put(26,106){\makebox(100,12){\fbox{\begin{tabular}{c}$(k,a)$-generalized
                              Fourier transform $\mathcal{F}_{k,a}$ 
                              \end{tabular}}}}
\put(7,85){$\scriptstyle a \to 2$}
\put(1,54){%
\setlength{\unitlength}{0.00088333in}

\begingroup\makeatletter\ifx\SetFigFont\undefined%
\gdef\SetFigFont#1#2#3#4#5{%
  \reset@font\fontsize{#1}{#2pt}%
  \fontfamily{#3}\fontseries{#4}\fontshape{#5}%
  \selectfont}%
\fi\endgroup%
{\renewcommand{\dashlinestretch}{30}
\begin{picture}(1530,2431)(0,-10)
\put(3172.000,-292.000){\arc{6328.507}{3.2365}{4.1638}}
\path(5.967,130.650)(22.000,8.000)(65.571,123.765)
\end{picture}
}
}
\put(140,85){$\scriptstyle a \to 1$}
\put(112,54){%
\setlength{\unitlength}{0.00088333in}
\begingroup\makeatletter\ifx\SetFigFont\undefined%
\gdef\SetFigFont#1#2#3#4#5{%
  \reset@font\fontsize{#1}{#2pt}%
  \fontfamily{#3}\fontseries{#4}\fontshape{#5}%
  \selectfont}%
\fi\endgroup%
{\renewcommand{\dashlinestretch}{30}
\begin{picture}(1529,2431)(0,-10)
\put(-1642.000,-292.000){\arc{6328.507}{5.2609}{6.1882}}
\path(1464.429,123.765)(1508.000,8.000)(1524.033,130.650)
\end{picture}
}
}
\put(75,100.5){$\Bigg\uparrow$}
\put(77,100.5){$\scriptstyle z \to \frac{\pi i}{2}$}
\put(51,84){\makebox(50,12){\fbox{\fbox{\begin{tabular}{c}$(k,a)$-generalized
                              Laguerre semigroup
                             $\mathcal{I}_{k,a}(z)$\end{tabular}}}}}
\put(45,79.5){$\scriptstyle a \to 2$}
\put(40,77.5){\rotatebox[origin=b]{30}{$\xleftarrow{\kern20mm}$}}
\put(90,77.5){\rotatebox[origin=b]{-30}{$\xrightarrow{\kern20mm}$}}
\put(100,79.5){$\scriptstyle a \to 1$}
\put(29,67){\fbox{$\mathcal{I}_{k,2}(z)$}}
\put(108,67){\fbox{$\mathcal{I}_{k,1}(z)$}}
\put(16,59){$\scriptstyle z \to \frac{\pi i}{2}$}
\put(20,58){\bigswarrow}
\put(43,58){\bigsearrow}
\put(50,59){$\scriptstyle k \to 0$}
\put(97,59){$\scriptstyle k \to 0$}
\put(99,58){\bigswarrow}
\put(122,58){\bigsearrow}
\put(129,59){$\scriptstyle z \to \frac{\pi i}{2}$}
\put(0,40){\framebox(36,12){\begin{tabular}{c}
                            Dunkl transform $\mathcal{D}_k$\\ 
                            \cite{D1}\end{tabular}}}
\put(38,40){\framebox(41,12){\begin{tabular}{c}
                            Hermite semigroup $I(z)$\\ 
                            \cite{xfol,H2}\end{tabular}}}
\put(83,40){\framebox(37,12){\begin{tabular}{c}
                              Laguerre semigroup\\ 
                              \cite{KM}\end{tabular}}}
\put(122,40){\framebox(29,12){\begin{tabular}{c}
                            $\mathcal{H}_k$\\
                            (see \eqref{eqn:HFk})\end{tabular}}}
\put(18,32){$\scriptstyle k \to 0$}
\put(20,33){\bigsearrow}
\put(43,33){\bigswarrow}
\put(50,32){$\scriptstyle z \to \frac{\pi i}{2}$}
\put(96,32){$\scriptstyle z \to \frac{\pi i}{2}$}
\put(99,33){\bigsearrow}
\put(122,33){\bigswarrow}
\put(129,32){$\scriptstyle k \to 0$}
\put(20,23){\fbox{Fourier transform}}
\put(99,23){\fbox{Hankel transform}}
\put(37,6.2){\rotatebox{90}{\hbox to 3.5em{\dotfill}}}
\put(48,12){$\Leftarrow$ `unitary inversion operator' $\Rightarrow$}
\put(117,6.2){\rotatebox{90}{\hbox to 3.5em{\dotfill}}}
\put(5,0){\parbox{15em}{the Weil representation of\\
                        the metaplectic group $Mp(N,\mathbb{R})$}}
\put(90,0){\parbox{15em}{the minimal representation of\\
                   the conformal group $O(N+1,2)$}}
\end{picture}
\label{fig:11}
\caption{Special values of holomorphic semigroup $\mathcal{I}_{k,a}(z)$}
\end{figure}

\subsection{Holomorphic semigroup $\mathcal{I}_{k,a}(z)$ with two
parameters $k$ and $a$}
\label{subsec:1.1}
\hfill\break
Dunkl operators are differential-difference operators associated to a
finite reflection group on the Euclidean space.
They were introduced by C. Dunkl \cite{D1}. 
This subject was motivated partly from harmonic
analysis on the tangent space of the Riemannian
 symmetric spaces, 
and resulted in a new theory of non-commutative harmonic analysis `without Lie groups'.
The Dunkl operators are also used as a tool for investigating an
algebraic integrability property for the Calogero--Moser quantum
problem related to root systems \cite{H}.
We refer to \cite{D08} for the up-to-date survey on various applications
of Dunkl operators.

Our holomorphic semigroup $\mathcal{I}_{k,a}(z)$ is built on Dunkl
operators.
To fix notation,
let 
\index{Cg@$\CC$}%
$\CC$ be the Coxeter group associated with a root system
\index{R@$\cal R$}%
$\mathcal{R}$ in $\R^N$.
For a $\mathfrak{C}$-invariant real function $k\equiv(k_\alpha)$
(\textit{multiplicity function}) on $\mathcal{R}$, 
we write $\Delta_k$ for the Dunkl Laplacian
 on $\R^N$
(see \eqref{eqn:2.3}).

We take $a>0$ to be a deformation parameter, 
and introduce the following differential-difference operator
\begin{equation}\label{eqn:kaLap}
\index{*Deltaka@$\Delta_{k,a}$|main}%
\Delta_{k,a}
 := \Vert x \Vert^{2-a} \Delta_k - \Vert x \Vert^a.
\end{equation}
Here, $\Vert x\Vert$ is the norm of the coordinate $x \in \R^N$,
and $\Vert x\Vert^a$ in the right-hand side of the formula 
stands for the multiplication operator by
$\Vert x\Vert^a$. 
Then, $\Delta_{k,a}$ is a symmetric operator on
the Hilbert space
$L^2(\mathbb{R}^N, \vartheta_{k,a}(x)dx)$
consisting of square integrable functions on $\R^N$ against
the measure $\vartheta_{k,a}(x)dx$, 
where the density function $\vartheta_{k,a}(x)$ on $\R^N$ is given by
\begin{equation}\label{eqn:01}
\index{*hthetakax@$\vartheta_{k,a}(x)$|main}%
\vartheta_{k,a}(x)
 := \|x\|^{a-2} \prod_{\alpha\in\mathcal{R}}
                   | \langle \alpha,x \rangle |^{k_\alpha}.
\end{equation}
Then $\vartheta_{k,a}(x)$ has a degree of homogeneity $a - 2 + 2 \langle k \rangle$,  
where $\langle k\rangle := \frac{1}{2} \sum_{\alpha\in\mathcal{R}} k_\alpha$
is the index of $k=(k_\alpha)$ (see \eqref{eqn:2.6}).

The \textit{$(k,a)$-generalized Laguerre semigroup} $\mathcal{I}_{k,a}(z)$ is 
defined to be the semigroup with infinitesimal generator 
$\frac{1}{a}\Delta_{k,a}$,
that is,
\begin{equation}\label{eqn:02}
\index{Ikaz@$\mathcal{I}_{k,a}(z)$|main}%
\mathcal{I}_{k,a}(z)
:=
 \exp\Bigl(\frac{z}{a}\Delta_{k,a}\Bigr),
\end{equation}
for $z \in \mathbb{C}$ such that $\operatorname{Re}z \ge 0$.
(Later, we shall use the notation
$\mathcal{I}_{k,a}(z) = \Omega_{k,a}(\gamma_z)$,
in connection with the Gelfand--Gindikin program.)

In the case $a=2$ and $k\equiv0$,
the density $\vartheta_{k,a}(x)$ reduces to $\vartheta_{0,2}(x) \equiv
1$ and we recover the classical setting where
\begin{align*}
&\Delta_{0,2} = \sum_{j=1}^N  \frac{\partial^2}{\partial x_j^2} 
                            - \sum_{j=1}^N x_j^2, 
 \ \text{the Hermite operator on $L^2(\R^N)$},
\\
&\mathcal{I}_{0,2}(z) = 
 \text{the Hermite semigroup $I(z)$ (\cite{xfol,H2}}).
\end{align*}

In this article,
we shall deal with a positive $a$ and a non-negative multiplicity function $k$
for simplicity,
though some of our results still hold for ``slightly-negative\rq\rq\ 
multiplicity functions (see Remark~\ref{rem:non-neg-k}).
We begin with:
\begin{thmalph}[see Corollary \ref{cor:kaLap}]\label{thm:A}
Suppose $a >0 $ and a non-negative multiplicity function $k$ satisfy
$
a+2\langle k\rangle+N-2>0.
$ Then,
\begin{enumerate}[\upshape 1)]
\item  
$\Delta_{k,a}$ extends to a self-adjoint operator on 
$L^2(\R^N,\vartheta_{k,a}(x)dx)$.
\item  
There is no continuous spectrum of $\Delta_{k,a}$.
\item  
All the discrete spectra are negative.
\end{enumerate}
\end{thmalph}
We also find all the discrete spectra explicitly in Corollary \ref{cor:kaLap}.

Turning to the $(k,a)$-generalized Laguerre semigroup
 $\mathcal{I}_{k,a}(z)$ (see \eqref{eqn:02}), we shall prove:
\begin{thmalph}[see Theorem \ref{thm:3.18}]\label{thm:B}
Retain the assumptions of Theorem \ref{thm:A}.
\begin{enumerate}[\upshape 1)]
\item  
$\mathcal{I}_{k,a}(z)$ is
a holomorphic semigroup in the complex right-half plane
$\{ z \in \C: \operatorname{Re}z > 0\}$ in the sense that 
$\mathcal{I}_{k,a}(z)$ is a Hilbert--Schmidt operator on
$L^2(\R^N,\vartheta_{k,a}(x)dx)$
satisfying
$$
\mathcal{I}_{k,a}(z_1) \circ \mathcal{I}_{k,a}(z_2)
= \mathcal{I}_{k,a}(z_1+z_2),
\quad (\operatorname{Re} z_1, \operatorname{Re} z_2 > 0),
$$
and that the scalar product $(\mathcal{I}_{k,a}(z) f, g)$ is a 
holomorphic function of $z$ for\/ $\operatorname{Re}z>0$,
 for any 
$f,g \in L^2(\R^N,\vartheta_{k,a}(x)dx)$.
\item  
$\mathcal{I}_{k,a}(z)$ is a one-parameter group of unitary operators
on the imaginary axis $\operatorname{Re}z=0$.
\end{enumerate}
\end{thmalph}

In Section \ref{sec:4.3},
we shall introduce a real analytic function
$\mathcal{I}(b,\nu;w;\cos\phi)$
in four variables defined on
$\{(b,\nu,w,\phi) \in \R_+\times\R\times\C\times\R/2\pi\Z: 1+b\nu>0\}$.
The special values
at $b=1,2$ are given by 
\begin{align}
&\cal J(1,\nu;w;t)=e^{wt},
\label{eqn:I1}
\\
&\cal J(2,\nu;w;t)=\Gamma(\nu+{1\over 2}) \widetilde I_{\nu-{1\over 2}}\Big({{w(1+t)^{1/2}}\over {\sqrt 2}}\Big).
\label{eqn:I2}
\end{align}
Here,
$\widetilde{I}_\lambda(z) = (\frac{z}{2})^{-\lambda} I_\lambda(z)$
is the (normalized) modified Bessel function of the first kind
(simply, $I$-Bessel function).
We notice that these are positive-valued functions of $t$ if $w\in \R$.

We then define the following continuous function of $t$ on the
interval $[-1,1]$ with parameters $r,s>0$ and $z\in \{z\in \C\;|\; \Re
z\geq 0\}\setminus i\pi \Z$ by
$$h_{k,a}(r,s;z;t)={{\exp\big(-{1\over a}(r^a+s^a)\coth (z)\big)}\over {\sh (z)^{{2\langle k\rangle +N+a-2}\over a}}} \cal J\left({2\over a}, {{2\langle k\rangle +N-2}\over 2}; {{2(rs)^{a\over 2}}\over {a \sh (z)}}; t\right),$$ 
where $\langle k\rangle  = \frac{1}{2} \sum_{\alpha\in\mathcal{R}} k_\alpha$  (see \eqref{eqn:2.6}).

For a function $h(t)$ of one variable,
let $(\widetilde{V}_k h)(x,y)$ be a $k$-deformation of the function
$h(\langle x,y\rangle)$ on $\R^N \times \R^N$.
(This $k$-deformation is defined by using the Dunkl intertwining
operator $V_k$, see \eqref{eqn:Vkhx}).

In the polar coordinates $x= r \omega$ and $y=s \eta$, 
 we set
$$\Lambda_{k,a}(x,y;z)=\widetilde V_k\big(h_{k,a}(r,s;z;\cdot)\big)(\omega, \eta).$$ 

For $a>0$ and a non-negative multiplicity function $k$,  
we introduce the following normalization constant 
\begin{equation}\label{cst}
c_{k,a}
:= (\int_{\R^N} \exp \left( -\frac{1}{a} \Vert x\Vert^a \right)
    \vartheta_{k,a} (x) dx )^{-1}.
\end{equation}
The constant $c_{k,a}$ can be expressed in terms  of the gamma function
owing to the work by Selberg, Macdonald, Heckman, Opdam \cite{Op}, and others
 (see Etingof \cite{E} for a uniform proof).  

Here is an integration formula of the holomorphic semigroup $\mathcal{I}_{k,a}(z).$
\begin{thmalph}[see Theorem \ref{thm:4.9}]
\label{thm:EE}
Suppose $a>0$ and $k$ is a non-negative multiplicity function.
Suppose $\Re z \ge 0$ and $z \notin i\pi\Z$.
Then, $\mathcal{I}_{k,a}(z) = \exp(\frac{z}{a}\Delta_{k,a})$ is given by 
\begin{equation}\label{eqn:LS}
\mathcal{I}_{k,a}(z)f(x)=c_{k,a} \int_{\R^N} f(y) \Lambda_{k,a}(x,y;z) \vartheta_{k,a}(y)dy.
\end{equation}
\end{thmalph}

The formula \eqref{eqn:LS} generalizes the $k\equiv0$ case;
see Kobayashi--Mano \cite{KM} for $(k,a)=(0,1)$,
and the Mehler kernel formula in Folland \cite{xfol} or Howe \cite{H2}
for $(k,a)=(0,2)$.

\subsection{$(k,a)$-generalized Fourier transforms
$\mathcal{F}_{k,a}$}
\label{subsec:1.2}
\hfill\break
As we mentioned in Theorem \ref{thm:B} 2),
the `boundary value' of the $(k,a)$-generalized 
Laguerre semigroup $\mathcal{I}_{k,a}(z)$
on the imaginary axis gives a one-parameter family of
unitary operators. 
The case $z=0$ gives the identity operator, namely, 
$\mathcal{I}_{k,a}(0) = \operatorname{id}$.
The particularly interesting case is when $z = \frac{\pi i}{2}$,
and we set 
$$
\index{Fka@$\Fka$}%
\Fka
 := c \, \mathcal{I}_{k,a} \Bigl(\frac{\pi i}{2}\Bigr)
= c \exp 
\Bigl( \frac{\pi i}{2a} (\Vert x\Vert^{2-a} \Delta_k - \Vert x\Vert^a) \Bigr)
$$
by multiplying the phase factor 
$
c = e^{i\frac{\pi}{2}(\frac{2\langle k \rangle + N+a-2}{a})}
$
(see \eqref{eqn:5.1}). 
Then, the unitary operator $\Fka$ for general $a$ and $k$ 
satisfies the following significant
properties:
\begin{thmalph}[see Proposition \ref{prop:3.13} and Theorem \ref{thm:Fka}]
\label{thm:D}
Suppose $a>0$ and $k$ is a non-negative multiplicity function such that $a+2\langle k\rangle+N-2>0$.
\begin{enumerate}[\upshape 1)]
\item 
$\Fka$ is a unitary operator on
$L^2(\R^N,\vartheta_{k,a}(x)dx)$.
\item 
$\Fka \circ E
 = - (E+N+2\langle k\rangle+a-2) \circ \Fka$.
\\
Here, $E = \sum_{j=1}^N x_j \partial_j$.
\item 
$\Fka \circ \Vert x\Vert^a
= - \Vert x\Vert^{2-a} \Delta_k \circ \Fka$,
\\
$\Fka \circ (\Vert x\Vert^{2-a} \Delta_k)
= -\Vert x\Vert^a \circ \Fka$.
\item 
$\Fka$ is of finite order if and only if $a\in\Q$.
Its order is $2p$ if $a$ is of the form $a=\frac{p}{q}$,
where $p$ and $q$ are positive integers that are relatively prime.
\end{enumerate}
\end{thmalph}

We call $\Fka$ a \textit{$(k,a)$-generalized Fourier transform} on $\R^N$.
We note that 
$\Fka$ reduces to the Euclidean Fourier transform $\mathcal{F}$ if $k\equiv0$ and $a=2$;
to the Hankel transform if $k\equiv0$ and $a=1$;
to the Dunkl transform $\mathcal{D}_k$ introduced by C. Dunkl himself
in \cite{D3} 
if $k>0$ and $a=2$.

For $a=2$,
our expressions of $\Fka$ amount to:
\begin{alignat*}{2}
&\mathcal{F} = e^{\frac{\pi i N}{4}} \exp\frac{\pi i}{4} (\Delta-\Vert x\Vert^2)
&\qquad&\text{(Fourier transform)},
\\
&
\index{Ds_k@$\protect\mathcal{D}_k$|main}%
\mathcal{D}_k = e^{\frac{\pi i(2\langle k\rangle+N)}{4}} \exp\frac{\pi i}{4} (\Delta_k-\Vert x\Vert^2)
&\qquad&\text{(Dunkl transform)}.
\end{alignat*}
For $a=1$ and $k \equiv 0$,
the unitary operator 
$$
\mathcal{F}_{0,1} = e^{\frac{\pi i(N-1)}{2}}
\exp \Bigl( \frac{\pi i}{2}\Vert x\Vert(\Delta-1)\Bigr)
$$
arises as the 
\textit{unitary inversion operator} of the Schr\"{o}dinger model of the minimal
representation of the conformal group $O(N+1,2)$
(see \cite{KM05,KM}).
Its Dunkl analogue, namely, the unitary operator $\Fka$ for $a=1$ and
$k>0$ seems also interesting, however,
it has never appeared in the literature, to the best of our knowledge.
The integral representation of
this unitary operator,
$$
\index{Hs_k@$\mathcal{H}_k$}%
\mathcal{H}_k
 := \mathcal{F}_{k,1} = 
e^{i\frac{\pi}{2}(2\langle k\rangle+N-1)}
\mathcal{I}_{k,1}\Bigl(\frac{\pi i}{2}\Bigr)
= e^{i\frac{\pi}{2}(2\langle k\rangle+N-1)}
 \exp\Bigl(\frac{\pi i}{2} \Vert x\Vert (\Delta_k-1)\Bigr),
$$
is given in terms of the Dunkl intertwining operator and the Bessel
 function due to the closed formula of $\mathcal{I}(b,\nu;w;t)$ at
 $b=2$ (see \eqref{eqn:I2}).

On the other hand,
our methods can be applied to general $k$ and $a$ in finding some
basic properties of the $(k,a)$-generalized Fourier transform $\Fka$
such as
the inversion formula, 
the Plancherel theorem,
the Hecke identity (Corollary \ref{coro:5.14}), 
 the Bochner identity (Theorem \ref{thm:5.15}),
and the following Heisenberg inequality (Theorem \ref{thm:5.20}):
\begin{thmalph}[Heisenberg type inequality]
\label{thm:E}
Let $\Vert \ \ \Vert_k$ denote by the norm on the Hilbert space
$L^2(\R^N,\vartheta_{k,a}(x)dx)$.
Then,
\[
\Bigl\| \, \Vert x\Vert^{\frac{a}{2}} f(x) \Bigr\|_k
\,
\Bigl\| \, \Vert \xi\Vert^{\frac{a}{2}} \Fka f(\xi) \Bigr\|_k
\ge
\frac{2\langle k\rangle+N+a-2}{2}
\, \Vert f(x)\Vert_k^2
\]
for any $f\in L^2(\R^N, \vartheta_{k,a}(x)dx)$.
The equality holds if and only if $f$ is a scalar multiple of\/ 
$\exp(-c\Vert x\Vert^a)$ for some $c>0$.
\end{thmalph}
This inequality was previously proved by R\"{o}sler \cite{R2} and
Shimeno \cite{Sh} for the $a=2$ case (i.e.\ the Dunkl transform $\mathcal{D}_k$). 
In physics terms we may think of the function where the equality holds
in Theorem \ref{thm:E} as a ground state;
indeed when $a=c=1$, $N=3$, and $k\equiv0$ it is exactly the wave function for the Hydrogen
atom with the lowest energy.

\subsection{$\mathfrak{sl}_2$-triple of
differential-difference operators}
\label{subsec:1.3}
\hfill\break
Over the last several decades, various works have been published that
develop applications of the representation theory of the special linear group $SL(2,\R)$.
We mention particularly the books of Lang \cite{L} and Howe--Tan
\cite{HT}, and the research papers of Vergne \cite{V} and Howe
\cite{H1}. These and other contributions show how the symmetries of
$\mathfrak{sl}_2$ can offer new perspectives on familiar topics from
inside and outside representation theory (character formulas, ergodic
theory, Fourier analysis, the Laplace equation, etc.).

The basic tool for the present article is also the $SL_2$ theory. 
We construct an $\mathfrak{sl}_2$-triple of
differential-difference operators with two parameters $k$ and $a$, 
and then apply representation theory of $\widetilde{SL(2,\R)}$,
the universal covering group of $SL(2,\R)$.
The resulting representation is a discretely decomposable unitary
representation in the sense of \cite{Kdisc},
which depends continuously on parameters $a$ and $k$.

To be more precise, 
we introduce the following differential-difference operators on
$\R^N \setminus \{0\}$ by
$$\EE_{k,a}^+:={i\over a} \Vert x\Vert^a,\qquad \EE_{k,a}^-:={i\over a}\Vert
x\Vert^{2-a} \Delta_k,\qquad \HH_{k,a}:={2\over a} \sum_{i=1}^N x_i
\partial_i+{{N+2\langle k\rangle+a-2}\over a}.$$
With these operators,
we have
$$
a\Delta_{k,a} = i \, (\EE_{k,a}^+ - \EE_{k,a}^-).
$$
The main point here is that our operator $\Delta_{k,a}$ can be
interpreted in the framework of the (infinite dimensional)
representation of the Lie algebra $\mathfrak{sl}(2,\R)$:
\begin{lemalph}[see Theorem \ref{thm:TDS}]\label{lem:E}
The differential-difference operators
$\{\HH_{k,a},\EE_{k,a}^+,\EE_{k,a}^-\}$ form an
$\mathfrak{sl}_2$-triple for any multiplicity-function $k$ and any
non-zero complex number $a$.
\end{lemalph}
In other words,
taking a basis of $\mathfrak{sl}(2,\R)$ as 
$$
{\bf e}^+ = \begin{pmatrix}0&1\\0&0\end{pmatrix}, \quad
{\bf e}^- = \begin{pmatrix}0&0\\1&0\end{pmatrix}, \quad
{\bf h} = \begin{pmatrix}1&0\\0&-1\end{pmatrix},
$$
we get a Lie algebra representation $\omega_{k,a}$ of
$\mathfrak{g}=\mathfrak{sl}(2,\R)$ with continuous parameters $k$ and
$a$ on functions on $\R^N$ by mapping
$$
{\bf h} \mapsto \HH_{k,a}, \quad
{\bf e}^+ \mapsto \EE_{k,a}^+, \quad
{\bf e}^- \mapsto \EE_{k,a}^-.
$$

The main result of Section \ref{sec:3} is to prove that the
representation $\omega_{k,a}$ of $\mathfrak{sl}(2,\R)$ lifts to the
universal covering group $\widetilde{SL(2,\R)}$:
\begin{thmalph}[see Theorem \ref{thm:3.12}]\label{thm:F}
If $a>0$ and $k$ is a non-negative multiplicity function such that $a+2\langle k\rangle+N-2>0$,
then $\omega_{k,a}$ lifts to a unitary representation of
$\widetilde{SL(2,\R)}$ on $L^2(\R^N,\vartheta_{k,a}(x)dx)$.
\end{thmalph}
Theorem \ref{thm:F} fits nicely into the framework of discretely decomposable unitary 
representations
\cite{Kdisc,Kaspm97}.  
In fact, we see in Theorem \ref{thm:3.16} that the Hilbert space
$L^2(\R^N,\vartheta_{k,a}(x)dx)$ decomposes discretely as a direct sum
of unitary representations of the direct product group
$\mathfrak{C} \times \widetilde{SL(2,\R)}$:
\begin{equation}\label{eqn:Hpideco}
L^2(\R^N,\vartheta_{k,a}(x)dx)
\simeq
\sideset{}{^\oplus}\sum_{m=0}^\infty\mathcal{H}_k^m(\R^N)_{ \big| S^{N-1}}
\otimes \pi \, \Bigl(\frac{2m+2\langle k\rangle+N-2}{a}\Bigr),
\end{equation}
where $\mathcal{H}_k^m(\R^N)$ stands for the representation of the
Coxeter group $\mathfrak{C}$ on the eigenspace of the Dunkl Laplacian
(the space of spherical $k$-harmonics of degree $m$) and $\pi(\nu)$ is
an irreducible unitary lowest weight representation of
$\widetilde{SL(2,\R)}$ of weight $\nu+1$ (see Fact \ref{fact:SL2rep}).
The unitary isomorphism \eqref{eqn:Hpideco} is constructed explicitly
by using Laguerre polynomials.

For general $N\ge 2$,
the right-hand side of \eqref{eqn:Hpideco} is an infinite sum.
For $N=1$,
\eqref{eqn:Hpideco} is reduced to the sum of
two terms ($m=0,1$).

The unitary representation of
$\widetilde{SL(2,\R)}$ on $L^2(\R^N,\vartheta_{k,a}(x)dx)$ extends
furthermore to a holomorphic semigroup of a complex three
dimensional semigroup (see Section \ref{sec:3.7}). 
Basic properties of the holomorphic semigroup
$\mathcal{I}_{k,a}(z)$ defined in \eqref{eqn:02} and the unitary
operator $\Fka$ can be read from the
`dictionary' of $\mathfrak{sl}(2,\R)$
as follows:
\begin{align*}
i\begin{pmatrix}0&1\\-1&0\end{pmatrix}
&\longleftrightarrow \frac{1}{a}\Delta_{k,a}
\\
\exp iz\begin{pmatrix}0&1\\-1&0\end{pmatrix}
&\longleftrightarrow \mathcal{I}_{k,a}(z) = \exp(\frac{z}{a}\Delta_{k,a})
\\
w_0 = \exp\frac{\pi}{2}\begin{pmatrix}0&-1\\1&0\end{pmatrix}
&\longleftrightarrow \Fka \ \text{(up to the phase factor)}
\\
\operatorname{Ad}(w_0){\bf e}^+ = {\bf e}^-
&\longleftrightarrow \Fka \circ \Vert x\Vert^a = -\Vert x\Vert^{2-a}\Delta_k\Fka
\\
\operatorname{Ad}(w_0){\bf e}^- = {\bf e}^+
&\longleftrightarrow \Fka \circ \Vert x\Vert^{2-a}\Delta_k = -\Vert x\Vert^a\Fka.
\end{align*}

\subsection{Hidden symmetries for $a=1$ and $2$}
\label{subsec:1.4}
\hfill\break
As we have seen in Section \ref{subsec:1.1},
one of the reasons that we find an explicit formula for the holomorphic semigroup
$\mathcal{I}_{k,a}(z)$ (and for the unitary operator $\Fka$) 
(see Section \ref{subsec:1.1}) is that there are large
 `hidden symmetries' on the Hilbert
space when $a=1$ or $2$.

We recall that our analysis is based on the fact that the Hilbert
space
$L^2(\R^N,\vartheta_{k,a}(x)dx)$ has a symmetry of the direct product
group $\CC \times \widetilde{SL(2,\R)}$ for all $k$ and $a$.
It turns out that
this symmetry becomes larger for special values of $k$ and $a$.
In this subsection,
we discuss these hidden symmetries.

First, in the case $k\equiv 0$,
the Dunkl Laplacian $\Delta_k$ becomes the Euclidean Laplacian
$\Delta$,
and consequently, not only the Coxeter group $\CC$ but also the whole
orthogonal group $O(N)$ commutes with $\Delta_k\equiv\Delta$.
Therefore, the Hilbert space
$L^2(\R^N,\vartheta_{0,a}(x)dx)$ is acted on by
$O(N)\times\widetilde{SL(2,\R)}$.
Namely, it has a larger symmetry
$$
\mathfrak{C} \times \widetilde{SL(2,\R)} \subset O(N)
\times \widetilde{SL(2,\R)}.
$$

Next, we observe that the Lie algebra of the direct product group 
$O(N)\times \widetilde{SL(2,\R)}$ may be seen as a subalgebra of two
different reductive Lie algebras
$\mathfrak{sp}(N,\R)$ and $\mathfrak{o}(N+1,2)$:
\begin{align*}
&\oo(N)\oplus\mathfrak{sl}(2,\R)
\simeq
\oo(N)\oplus\oo(1,2) \subset \oo(N+1,2)
\\
&\oo(N)\oplus\mathfrak{sl}(2,\R)
\simeq
\oo(N)\oplus\mathfrak{sp}(1,\R) \subset \mathfrak{sp}(N,\R)
\end{align*}
It turns out that they are the hidden symmetries of the Hilbert space
$L^2(\R^N,\vartheta_{0,a}(x)dx)$ for $a=1,2$, respectively.
To be more precise,
the conformal group $O(N+1,2)_0$ (or its double covering group if $N$
is even) acts on 
$L^2(\R^N,\vartheta_{0,1}(x)dx) = L^2(\R^N,\Vert x\Vert^{-1}dx)$
as an irreducible unitary representation,
while the metaplectic group $Mp(N,\R)$ (the double covering group of
the symplectic group $\mathit{Sp}(N,\R)$) acts on
$L^2(\R^N,\vartheta_{0,2}(x)dx) = L^2(\R^N,dx)$
as a unitary representation.

In summary,
we are dealing with the symmetries of the Hilbert space
$L^2(\R^N,\vartheta_{k,a}(x)dx)$ described below:
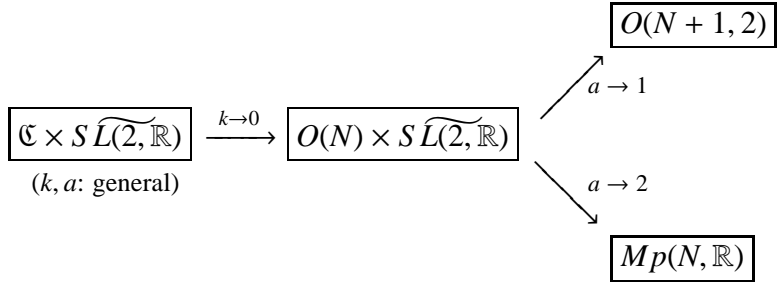
\begin{figure}[H]
\renewcommand{\figurename}{Diagram}
\renewcommand{\thefigure}{\thesubsection}
\begin{center}
\setlength{\unitlength}{1mm}
\begin{picture}(100,43)
\put(80,36){\fbox{$O(N+1,2)$}}
\put(77,28){\scriptsize $a\to1$}
\put(69,28){\bignearrow}
\put(0,21){\fbox{$\mathfrak{C}\times\widetilde{SL(2,\R)}$}}
\put(26,21){$\xrightarrow{\ k\to 0\ }$}
\put(37,21){\fbox{$O(N)\times\widetilde{SL(2,\R)}$}}
\put(3,15){\footnotesize ($k,a$: general)}
\put(69,14){\bigsearrow}
\put(80,5){\fbox{$Mp(N,\R)$}}
\put(77,15){\scriptsize $a\to2$}
\end{picture}
\end{center}
\label{fig:Hs}
\caption{Hidden symmetries in $L^2(\R^N,\vartheta_{k,a}(x)dx)$}
\end{figure}

For $a=2$, this unitary representation is nothing but the Weil
representation,
sometimes referred to as the Segal--Shale--Weil representation,
the metaplectic representation, 
or the oscillator representation,
and its realization on $L^2(\R^N)$ is called
the Schr\"{o}dinger model.

For $a=1$,
the unitary representation of the conformal group on 
$L^2(\R^N,\Vert x\Vert^{-1}dx)$ is irreducible and has a similar
nature to the Weil representation.
The similarity is illustrated by the fact that
both of these unitary representations are `minimal representations',
i.e.,
 their annihilator of the infinitesimal representations
are the Joseph ideal of the universal enveloping algebras,
and in particular, they attain the minimum of their Gelfand--Kirillov
dimensions.

In this sense, our continuous parameter $a>0$ interpolates two minimal
representations of different reductive groups by keeping smaller
symmetries
(i.e.\ the representations of $O(N)\times\widetilde{SL(2,\R)}$).
The $(k,a)$-generalized Fourier transform $\Fka$ plays a special
role in the global formula of the $L^2$-model of minimal
representations. 
In fact, the conformal group $O(N+1,2)$ is generated by a maximal parabolic subgroup (essentially, the affine conformal group for the Minkowski space $\R^{N,1}$) and the 
 inversion element $I_{N+1,2}={\rm diag}(1,\ldots,1,-1, -1).$ Likewise, the metaplectic group 
 $Mp(N,\R)$ is generated by the Siegel parabolic subgroup and the conformal inversion element. 
Since the maximal parabolic subgroup acts on the $L^2$-model
on the minimal representation, we can obtain
the global formula of the whole group
 if we determine the action of the inversion element.
For the Weil representation, this crucial action is nothing but the
Euclidean Fourier transform (up to the phase factor),
and it is the Hankel transform for the minimal representation of the
conformal group $O(N+1,2)$ (see \cite{KM05},
see also \cite{KM-intopq} and \cite[Introduction]{KM-intopq-ams} for some perspectives of this
direction in a more general setting).

A part of the results here has been announced in \cite{BKO-CRAS}
without proof.

\medbreak

\textbf{Notation}\enspace
$\N = \{ 0,1,2,\dotsc \}$,
$\N_+ = \{ 1,2,3,\dotsc \}$,
$\R_+ = \{ x \in \R \mid x>0 \}$,
and
$\R_{\ge0} = \{ t \in \R: t \ge 0 \}$.

\section{Preliminary results on Dunkl operators}\label{sec:2}
\subsection{Dunkl operators}\label{subsec:2.1}
\hfill\break
 Let $\langle
\cdot,\cdot\rangle$ be the standard Euclidean scalar product in
$\R^N$.
We shall use the same notation for its bilinear extension to $\C^N\times \C^N.$ For
$x\in \R^N,$ denote by $\Vert x\Vert=\langle x,x\rangle^{1/2}$.

For $\alpha\in \R^N\setminus \{0\},$ we write $r_\alpha$ for the
reflection with respect to the hyperplane $\langle \alpha\rangle^\perp$
orthogonal to $\alpha$ defined by
$$ r_\alpha(x):=x-2{{\langle \alpha,x\rangle}\over {\Vert \alpha
\Vert^2}}\alpha,\qquad x\in \R^N.$$

We say a finite  set  
$
\index{R@$\cal R$|main}%
\cal R 
$ in $\R^N \setminus \{0\}$ is a (reduced) \textit{root system} if: 
\begin{itemize}
\item[(R1)] $r_\alpha(\cal R)=\cal R$ for all $\alpha\in \cal R,$
\item[(R2)] $\cal R\cap \R\alpha =\{\pm \alpha\}$ for all $\alpha\in
\cal R$.
\end{itemize}
In this article,
we do not impose
 crystallographic conditions on the roots, and do not require that $\cal R$ spans $\R^N.$ 
However, we shall assume $\mathcal{R}$ is reduced, namely,
(R2) is satisfied.

The subgroup 
$
\index{Cg@$\CC$|main}%
\CC 
\subset O(N, \R)$    generated by the reflections  $\{r_\alpha\;|\;\alpha\in \cal R\}$ is called the finite Coxeter group  associated with  $\cal R.$  
The Weyl groups such as the symmetric group $\mathfrak{S}_N$ for the
type $A_{N-1}$ root system and the hyperoctahedral group for the type
$B_N$ root system are typical examples.
In addition,
$H_3, H_4$ (icosahedral groups) and $I_2(n)$ (symmetry group of the
regular $n$-gon) are also the Coxeter groups.
We refer to \cite{Hy}  for more details on the theory of Coxeter groups. 

\begin{defe}\label{def:mf}
A multiplicity function for $\CC$ is a function 
$k:\cal R\rightarrow \C$
which is constant on $\CC$-orbits.
\end{defe}
Setting $k_\alpha:=k(\alpha)$ for
$\alpha\in \cal R, $ we have $k_{h\alpha}=k_\alpha$ for all $h\in \CC$ 
from definition. 
We say $k$ is non-negative if $k_\alpha \ge 0$ for all $\alpha \in \mathcal{R}$.
The $\C$-vector space of multiplicity functions on $\cal R$ is
denoted by $\cal K.$  
The dimension of $\mathcal{K}$ is equal to the number of
$\CC\text{-orbits in}\; \cal R$.

For $\xi\in \C^N$ and $k\in \cal K,$   Dunkl \cite{D1}
introduced a family of first order differential-difference operators
\index{Txik@$T_\xi(k)$|main}%
$T_\xi(k)$ 
(\textit{Dunkl's operators}) by
\begin{equation}\label{eqn:2.1}
T_\xi(k)f(x):=\partial_\xi f(x)+\sum_{\alpha\in \cal
R^+}k_\alpha\langle \alpha, \xi\rangle{{f(x)-f(r_\alpha x)}\over
{\langle \alpha, x\rangle}},\qquad f\in C^1(\R^N). \end{equation} Here
$\partial_\xi$ denotes the directional derivative corresponding to
$\xi.$  Thanks to the $\CC$-invariance of the multiplicity function, this definition is independent of the choice
of the positive subsystem $\cal R^+.$ The  operators  $T_\xi(k)$ are 
homogeneous of degree $-1.$ Moreover, the Dunkl operators satisfy the
following properties (see \cite{D1}): 
\begin{itemize}
\item[(D1)] $L(h)\circ T_\xi(k)\circ L(h)^{-1} =T_{h\xi}(k)$  for all $h\in \CC,$
\item[(D2)] $T_\xi(k)T_\eta(k)=T_\eta(k)T_\xi(k)$ for all $\xi, \eta\in \R^N,$
\item[(D3)] $T_\xi(k) [fg ]=gT_\xi(k) f +fT_\xi(k)g$ if $f$ and $g$ are in $C^1(\R^N)$ and   at least one of them is $\CC$-invariant.
\end{itemize}
Here, we denote by $L(h)$ the left regular action of $h \in
\mathfrak{C}$ on the function space on $\R^N$:
\[
(L(h)f)(x):=f(h^{-1}\cdot x).
\]

\begin{rem}\label{rem:2.1}   
The Dunkl Laplacian arises as the radial part of the Laplacian on the
tangent space of a Riemannian symmetric spaces.
Let $\g$ be a real  semisimple Lie algebra with  Cartan decomposition $\g=\k\oplus \p.$  
We take a maximal abelian subspace 
$\a$ in $ \p,$ and let $\Sigma(\g, \a)$ be the set of   restricted roots,
and $m_\alpha$ the multiplicity of $\alpha\in \Sigma (\g,\a).$ 
We may consider $\Sigma(\g, \a)$  to  be a subset of $\a$ by means of the Killing form of $\g.$  
The Killing form endows $\p$ with a flat Riemannian symmetric
space structure, and we write $\Delta_{\p}$ for the (Euclidean)
Laplacian on $\p.$   Put $\cal R:=2\Sigma(\g, \a)$ and
$k_\alpha:=\frac{1}{2}\sum_{\beta\in\Sigma^+\cap\R\alpha}m_\beta$.  
We note that the root system $\mathcal{R}$ is not necessarily reduced.
Then the radial part of  $\Delta_{\p}$,
denoted by $Rad(\Delta_{\mathfrak{p}})$,
(see \cite[Proposition 3.13]{Helg})
is given by
$$
Rad(\Delta_{\p}) f =\Delta_k f
$$ 
for every  ${\CC}$-invariant function 
$f\in C^\infty(\a)$,
where $\Delta_k$ is the Dunkl Laplacian which will be defined in
\eqref{eqn:2.3}. 
\end{rem}
\begin{rem}\label{rem:non-neg-k}
Some of our results still hold for ``slightly-negative\rq\rq\ 
multiplicity functions.
For instance, when $k_\alpha=k$ for all $\alpha\in \cal R,$
 we may relax the assumption $k\geq 0$ by $k> -{1\over d_{\max}}$ where
$d_{\max}$ is the largest fundamental degree of the Coxeter group $\CC$
(see \cite[Theorem 3.1]{E}). However, for simplicity, we will
restrict ourselves to non-negative multiplicity functions
$k=(k_\alpha)_{\alpha\in \cal R}$.
\end{rem}

Let $\vartheta_k$ be the weight function on $\R^N$ defined by
\begin{equation}\label{eqn:2.5}
\index{*htheta_k(x)@$\vartheta_k(x)$|main}%
\vartheta_k(x) 
:=\prod_{\alpha\in \cal R^+} \vert \langle \alpha,
x\rangle\vert^{2k_\alpha},\qquad x\in \R^N. \end{equation} It is $\CC$-invariant
and homogeneous of degree $2\langle k\rangle,$ where the index
$\langle k\rangle$ of the
multiplicity function $k$ is defined as
\begin{equation}\label{eqn:2.6}
\index{k@$ \langle k\rangle$|main}%
\langle k\rangle 
:=\sum_{\alpha\in \cal R^+}k_\alpha
 = \frac{1}{2}\sum_{\alpha\in\mathcal{R}}k_\alpha.
\end{equation}
Let $dx$ be the Lebesgue measure on $\R^N$ with respect to the inner
product $\langle \ , \ \rangle$.
Then  the Dunkl operators are skew-symmetric with respect to the
measure $\vartheta_k(x)dx$ (see \cite{D1}).
In particular,
 if $f$ and $g$ are differentiable and one of them has
compact support, then 
\begin{equation}\label{eqn:2.7}
\int_{\R^N} (T_\xi(k)f)(x)g(x) \vartheta_k(x)dx=-\int_{\R^N }f(x) (T_\xi(k)g)(x)\vartheta_k(x)dx.\end{equation}

It  is shown in \cite{D2} that for any non-negative root multiplicity
function $k$ there is a unique linear isomorphism 
\index{V_k@$V_k$|main}%
$V_k$ 
(\textit{Dunkl's
intertwining operator}) on the space $\mathcal P(\R^N)$ of
polynomial functions on $\R^N$ such that:
\begin{itemize}
\item[(I1)] $V_k(\mathcal P_m(\R^N)) =\mathcal P_m(\R^N)  $ for all $m\in \N,$
\item[(I2)] ${V_k}_{ | \mathcal P_0(\R^N)}= \text  {id},$
\item[(I3)] $T_\xi(k) V_k= V_k \partial_\xi$ for all $\xi\in \R^N.$
\end{itemize}
Here, $\mathcal P_m(\R^N)$ denotes the space of 
homogeneous polynomials of degree $m.$  It is known that
 $V_k$ induces a homeomorphism of $C(\R^N)$ and also that of
$C^\infty(\R^N)$ (cf. \cite{T}).
See also \cite{DJO} for more results on $V_k$ for $\C$-valued
multiplicity functions on $\mathcal{R}$.

For arbitrary finite reflection group $\CC,$  and for any
non-negative multiplicity function $k$, 
R\"{o}sler \cite{R1} proved that there exists a unique positive Radon probability-measure 
\index{*mu_x^k@$\mu_x^k$|main}%
$\mu_x^k $ 
on $\R^N$ such that \begin{equation}\label{eqn:2.8}
 V_k f(x)=\int_{\R^N} f(\xi) d\mu_x^k(\xi). \end{equation}
The measure $\mu_x^k$ depends on $x\in\R^N$ and its
support is contained in the ball $B(\Vert x\Vert) := \{\xi\in \R^N\;|\;
\Vert \xi\Vert\leq \Vert x\Vert\}.$ 
Moreover, for any Borel set $S\subset \R^N,$  $g\in \CC $ and $r>0,$ 
the following invariant property holds:
$$
\mu_x^k(S) = \mu_{gx}^k(gS) = \mu_{rx}^k(rS).
$$
In view of the Laplace type representation \eqref{eqn:2.8}, 
Dunkl's intertwining operator  $V_k$ can be extended to a larger class of spaces.
For example,
let $B$ denote the closed unit ball in $\R^N$.
Then the support property of $\mu_x^k$ leads us to the following:
\begin{lem}\label{lem:VkB}
For any $R>0$,
$V_k$ induces a continuous endomorphism of
$C(B(R))$.
\end{lem}

\begin{proof}
Let $f \in C(B(R))$.
We extend $f$ to be a continuous function $\widetilde{f}$ on $\R^N$.
Then, $V_k\widetilde{f}$ is given by the integral
\[
V_k\widetilde{f}(x)
= \int_{\R^N} \widetilde{f}(\xi) d\mu_x^k(\xi).
\]
Suppose now $x \in B(R)$.
Then $\operatorname{Supp} \mu_x^k \subset B(\Vert x\Vert) \subset B(R)$.
Hence, $(V_k\tilde{f})|_{B(R)}$ is determined by the restriction
$f = \tilde{f}|_{B(R)}$.
Thus, the correspondence $f \mapsto (V_k \tilde{f})|_{B(R)}$ is
well-defined, and we get an induced
 linear map $V_k$:
$C(B(R)) \to C(B(R))$,
by using the same letter.

Next, suppose a sequence $f_j \in C(B(R))$ converges uniformly to
$f\in C(B(R))$ as $j \to \infty$.
Then we can extend $f_j$ to a continuous function
$\widetilde{f}_j$ on $\R^N$ such that $\widetilde{f}_j$ converges to
$\widetilde{f}$ on every compact set on $\R^N$.
Hence $V_k\widetilde{f}_j$ converges to $V_k\widetilde{f}$,
and so does $V_k f_j$ to $V_k f$.
\end{proof}
For a continuous function $h(t)$ of one variable,
we set
\[
h_y(\cdot) := h(\langle \cdot, y\rangle)
\quad
(y \in \R^N),
\]
and define
\begin{equation}\label{eqn:Vkhx}
(\widetilde{V}_k h)(x,y)
:= (V_k h_y)(x)
= \int_{\R^N} h(\langle \xi,y\rangle) d\mu_x^k(\xi).
\end{equation}
Then, $(\widetilde{V}_k h)(x,y)$ is a continuous function on
$(x,y) \in \R^N \times \R^N$.

We note that if $k\equiv0$ then
\[
(\widetilde{V}_0 h)(x,y) = h(\langle x,y \rangle).
\]

If $h(t)$ is defined only near the origin,
we can still get a continuous function
$(\widetilde{V}_k h)(x,y)$ as far as $|\langle x,y \rangle|$ is
sufficiently small.
To be more precise,
we prepare the following proposition for later purpose.
For simplicity, we write $B$ for the unit ball $B(1)$ in $\R^N$.
\begin{prop}\label{prop:VkB}
Suppose $h(t)$ is a continuous function on the closed interval\/ 
$[-1,1]$.
Then, $(\widetilde{V}_k h)(x,y)$ is a continuous function on 
$B \times B$.
Further, $\widetilde{V}_k h$ satisfies
\begin{align}\label{eqn:Vkhsup}
\Vert \widetilde{V}_kh \Vert_{L^\infty(B\times B)}
&\le \Vert h \Vert_{L^\infty([-1,1])}
\\
(\widetilde{V}_k h)(x,y)
&= (\widetilde{V}_k h)(y,x).
\label{eqn:VkBxy}
\end{align}
\end{prop}
\begin{proof}
We extend $h$ to a continuous function $\widetilde{h}$ on $\R$.
It follows from Lemma \ref{lem:VkB} that the values
$(\widetilde{V}_k \widetilde{h})(x,y)$ for $(x,y)$ satisfying 
$|\langle x,y\rangle| \le 1$ are determined by the restriction
$h = \widetilde{h}|_{[-1,1]}$.
Hence,
\[
(\widetilde{V}_k h)(x,y)
:= (\widetilde{V}_k \widetilde{h})(x,y),
\quad
(x,y) \in B \times B
\]
is well-defined.

Since $\mu_x^k$ is a probability measure,
we get an upper estimate
\eqref{eqn:Vkhsup} from the integral expression \eqref{eqn:Vkhx}.

By the Weierstrass theorem,
we can find a sequence of polynomials $h_j(t)$ $(j=1,2,\dotsc)$
such that $h_j(t)$ converges to $\widetilde{h}(t)$ uniformly on any
compact set of $\R$.
Then, $\widetilde{V}_k h_j$ converges to $\widetilde{V}_k h$ uniformly
on $B \times B$.
Thanks to \cite[Proposition 3.2]{D2},
we have
$(\widetilde{V}_k h_j)(x,y) = (\widetilde{V}_k h_j)(y,x)$.
Taking the limit as $j$ tends to infinity,
we get the equation \eqref{eqn:VkBxy}.
Hence, Proposition \ref{prop:VkB} is proved.
\end{proof}

Aside from the development of the general theory of the Dunkl
transform,
we note that explicit formulas for $V_k$ have been known
for only a few cases:  
$\CC=\Z_2^N$, $\CC=S_3 $, and the equal parameter case for the Weyl
group of $B_2$ 
(see \cite{D08} for the recent survey by C. Dunkl). 

\subsection{The Dunkl Laplacian}\label{subsec:DLap}
\hfill\break
Let $\{\xi_1,\ldots, \xi_N\}$ be an orthonormal basis of $(\R^N, \langle\cdot,\cdot\rangle).$ For the $j$-th
basis vector $\xi_j,$ we will use the abbreviation
$T_{\xi_j}(k)=T_j(k).$ 
The Dunkl--Laplace
operator, or simply, the Dunkl Laplacian, is defined as
\begin{equation}\label{eqn:2.2}
\index{*Delta_k@$\Delta_k$|main}%
\Delta_k
:=\sum_{j=1}^NT_{j}(k)^2. 
\end{equation}
The definition of $\Delta_k$ is independent of the choice of an
orthonormal basis of $\R^N$.
In fact, it is proved in \cite{D1} that $\Delta_k$ is expressed as
\begin{equation}\label{eqn:2.3}
\Delta_kf(x)=\Delta f(x)+ \sum_{\alpha\in \cal
R^+}k_\alpha\left\{ {{2\langle  \nabla f(x), \alpha\rangle}\over
{\langle \alpha, x\rangle}}-\Vert \alpha\Vert^2{{f(x)-f(r_\alpha x)}\over {\langle
\alpha, x\rangle^2}}\right\},\end{equation} where $\nabla$ denotes  the usual gradient operator. 

For $k \equiv 0$, the Dunkl--Laplace operator $\Delta_k$ reduces to the
Euclidean Laplacian $\Delta$,
which commutes with the action of $O(N)$.
For general $k$, 
it follows from (D1) and \eqref{eqn:2.2} that $\Delta_k$
commutes with  the action of the Coxeter group $\CC$, i.e. 
\begin{equation}\label{eqn:2.4}
L(h)\circ \Delta_k\circ L(h)^{-1}=\Delta_k,\qquad \forall h\in \CC.
\end{equation}

\begin{defe}\label{def:Hkm}
A $k$-harmonic polynomial of degree $m$ $(m\in \N)$
is a homogeneous polynomial $p$ on $\R^N$ of degree $m$ such that $\Delta_kp =0.$ 
\end{defe}

Denote by  
$
\index{Hskm(\R^N)@$\Hkm(\R^N)$|main}%
\Hkm(\R^N)
$ 
the space of $k$-harmonic polynomials  of degree $m$.
It is naturally a representation space of the Coxeter group $\CC$.

Let $d\sigma$ be the standard
measure on $S^{N-1}$, 
$\vartheta_{k}$  the
density given in
\eqref{eqn:2.5},  
and $d_k$  the normalizing constant defined by
\begin{equation}\label{eqn:dk}
\index{d_k@$d_k$|main}%
d_k 
:=\Big(\int_{S^{N-1}} \vartheta_k(\omega) d\sigma(\omega)\Big)^{-1}.
\end{equation}
We write $L^2(S^{N-1}, \vartheta_k(\omega)d\sigma(\omega))$ for the
Hilbert space with
the following  inner product $\langle \ , \ \rangle_k$ given by
$$\langle f,g\rangle_k :=d_k \int_{S^{N-1}} f(\omega) \overline{g(\omega)}
\vartheta_k(\omega) d\sigma(\omega).$$ 

For $k\equiv0$,
$d_k^{-1}$ is the volume of the unit sphere, namely,
\begin{equation}\label{eqn:d0}
d_0 = \frac{\Gamma(\frac{N}{2})}{2\pi^{\frac{N}{2}}}.
\end{equation}
Thanks to Selberg, Mehta,
  Macdonald \cite{bigMac}, Heckman, Opdam \cite{Op}, and others,
there is a closed form of $d_k$ in terms of Gamma functions when $k$ is a non-negative multiplicity function (see also \cite{E}).

As in the classical spherical harmonics (i.e.\ the $k\equiv0$ case), we have
(see \cite[page 37]{D88}): 
\begin{fact}\label{fact:kH}
~
\begin{enumerate}[\upshape 1)]
\item 
$\Hkm(\R^N)|_{S^{N-1}}$ $(m=0,1,2,\dotsc)$ are orthogonal to each
other with respect to $\langle \ ,\ \rangle_k$.
\item 
The Hilbert space $L^2(S^{N-1},\vartheta_k(\omega)d\sigma(\omega))$
decomposes as a direct Hilbert sum:
\begin{equation}\label{3.15}
 L^2\big(S^{N-1}, \vartheta_k(\omega) d\sigma(\omega)\big)=\sideset{}{^\oplus}\sum_{{m}\in \N} \Hkm(\R^N)|_{S^{N-1}}.\end{equation}
\end{enumerate}
\end{fact}

We pin down some basic formulae of
$\Delta_k$.
We write the Euler operator as
\begin{equation}\label{eqn:Eu}
\index{Ea@$E$|main}%
E := \sum_{j=1}^N x_j\partial_j.
\end{equation}
\begin{lem}\label{lem:DLap}
\begin{enumerate}[\upshape 1)]
\item 
The Dunkl Laplacian $\Delta_k$ is of degree $-2$, namely,
\begin{equation}\label{eqn:DLap2}
[E,\Delta_k] = -2\Delta_k
\end{equation}
\item 
\begin{equation}\label{3.12} 
\sum_{j=1}^N \big( x_jT_j(k)+T_j(k)x_j\big) = N+2\langle k\rangle+2 E .
\end{equation}
\item 
Suppose $\psi(r)$ is a $C^\infty$ function of one variable.
Then we have
\begin{equation}\label{eqn:Dkpsi}
[\Delta_k,\psi(\Vert x\Vert^a)]
= a^2\Vert x\Vert^{2a-2} \psi''(\Vert x\Vert^a)
  + a\Vert x\Vert^{a-2} \psi'(\Vert x\Vert^a) ((N+2\langle k\rangle+a-2)+2E).
\end{equation}
\end{enumerate}
\end{lem}

\begin{proof}
See \cite[Theorem 3.3]{H} for 1) and 2).

3) 
Take an arbitrary $C^\infty$ function $f$ on $\R^N$.
We recall from the definition \eqref{eqn:2.1} and (D3) that
\begin{align}
& T_j(k) g = \partial_j g,
\label{eqn:Tj g}
\\
& T_j(k)(fg)
= (T_j(k) f) g + f(\partial_j g),
\nonumber
\end{align}
if $g$ is a $\mathfrak{C}$-invariant function on $\R^N$.
In particular,
\begin{align*}
& T_j(k) \psi(\Vert x\Vert^a)
= ax_j \Vert x\Vert^{a-2} \psi'(\Vert x\Vert^a),
\\
& T_j(k) (f(x)\psi(\Vert x\Vert^a))
= (T_j(k) f(x)) \psi(\Vert x\Vert^a) + ax_j f(x) \Vert x\Vert^{a-2}
\psi'(\Vert x\Vert^a).
\end{align*}
Using (D3) again,
we get
\begin{align*}
 T_j(k)^2 (f(x) \psi(\Vert x\Vert^a))
={}
&(T_j(k)^2 f(x)) \psi(\Vert x\Vert^a) 
\\
&+ a\Vert x\Vert^{a-2} \psi'(\Vert x\Vert^a)
   (x_j(T_j(k)f(x)) + T_j(k)(x_j f(x)))
\\
&+ a f(x) x_j T_j(k) (\Vert x\Vert^{a-2} \psi'(\Vert x\Vert^a)).
\end{align*}
Taking the summation over $j$,
we arrive at
\begin{align*}
\Delta_k(f(x)\psi(\Vert x\Vert^a))
={}
& (\Delta_k f(x)) \psi(\Vert x\Vert^a)
 + a\Vert x\Vert^{a-2} \psi' (\Vert x\Vert^a)
   (2E+N+2\langle k\rangle) f(x)
\\
&+ a f(x) E (\Vert x\Vert^{a-2} \psi' (\Vert x\Vert^a )).
\end{align*}
Here, we have used the expression \eqref{eqn:2.2} of $\Delta_k$,
\eqref{3.12}, and \eqref{eqn:Tj g}.
Now, \eqref{eqn:Dkpsi} follows from the following observation:
 in the polar coordinate $x=r\omega$,
the Euler operator $E$ amounts to $r\frac{\partial}{\partial r}$,
and
$ 
r\frac{d}{dr} (r^{a-2} \psi'(r^a))
= (a-2) r^{a-2} \psi'(r^a) + ar^{2a-2} \psi'' (r^a)
$.
\end{proof}

To end this section,
we consider a `$(k,a)$-deformation' of the classical formula
\[
e^{\Vert x\Vert^2} \circ \Delta \circ e^{-\Vert x\Vert^2}
= \Delta + 4\Vert x\Vert^2 -2N - 4 E.
\]

\begin{lem}\label{lem:3.4}
For any $\nu\in\C$ and $a\ne0$,
we have
\begin{equation}\label{eqn:kLapc}
e^{\frac{\nu}{a}\Vert x\Vert^a} \circ \Vert x\Vert^{2-a} \Delta_k
\circ e^{-\frac{\nu}{a}\Vert x\Vert^a}
= \Vert x\Vert^{2-a}\Delta_k + \nu^2\Vert x\Vert^a -
  \nu( (N+2\langle k\rangle+a-2) + 2 E).
\end{equation}
\end{lem}
\begin{proof} 

The proof parallels to that of Lemma \ref{lem:DLap} 3).
By the property (D3) of the Dunkl operators,
we get
\[
T_j(k) (e^{\lambda\Vert x\Vert^a} h(x))
= (T_j(k) e^{\lambda\Vert x\Vert^a}) h(x)
  + e^{\lambda\Vert x\Vert^a} T_j(k) h(x).
\]
Then, substituting the formula
\begin{equation*}
T_j(k) e^{\lambda\Vert x\Vert^a}
= \partial_j e^{\lambda\Vert x\Vert^a}
= \lambda a x_j\Vert x\Vert^{a-2} e^{\lambda\Vert x\Vert^a},
\end{equation*}
we have
\begin{equation}\label{3.9}
e^{-\lambda\Vert x\Vert^a} \circ T_j(k) \circ e^{\lambda\Vert
x\Vert^a} h(x)
= \lambda a x_j \Vert x\Vert^{a-2} h(x) + T_j(k) h(x).
\end{equation} 
Iterating \eqref{3.9} and using
\[
T_j(k)\Vert x\Vert^{a-2}
= (a-2) x_j\Vert x\Vert^{a-4},
\]
we get
\begin{align*}
e^{-\lambda\Vert x\Vert^a} \circ T_j(k)^2 \circ e^{\lambda\Vert x\Vert^a}
&= (\lambda ax_j \Vert x\Vert^{a-2} + T_j(k))^2
\\
&= \lambda^2a^2x_j^2 \Vert x\Vert^{2a-2}
   + \lambda a\Vert x\Vert^{a-2} (x_jT_j(k)+T_j(k)x_j)
\\
&\quad +\lambda a(a-2) x_j^2\Vert x\Vert^{a-4} + T_j(k)^2.
\end{align*}
Summing them up over $j$, we have
\begin{equation}\label{3.10}
e^{-\lambda\Vert x\Vert^a} \circ \Delta_k \circ e^{\lambda\Vert x\Vert^a}
= \Delta_k + \lambda^2 a^2 \Vert x\Vert^{2a}
  + \lambda a \Vert x\Vert^{a-2} (a-2+ \sum_{j=1}^N
    (x_j T_j(k)+T_j(k)x_j))
\end{equation}
The substitution of \eqref{3.12} and $\lambda=-\frac{\nu}{a}$ 
to  \eqref{3.10} shows Lemma. 
\end{proof}

\section{The infinitesimal representation $\omega_{k,a}$ of  $\s\l(2,\R)$}
\label{sec:3}
\subsection{$\mathfrak{sl}_2$ triple of differential-difference operators}
\label{subsec:3.1} \hfill\\
In this subsection,
we construct a family of Lie algebras which are isomorphic to
$\mathfrak{sl}(2,\R)$ in the space of differential-difference
operators on $\R^N$.
This family is parametrized by a non-zero complex number $a$ and a
multiplicity function $k$ for the Coxeter group.

We take a basis for the Lie algebra $\s\l(2,\R)$ as
\begin{equation}\label{eqn:3.1}
\index{ey^+@${\bf e}^+$|main}%
{\bf e}^+ 
:=\begin{pmatrix}0&1\\0&0\end{pmatrix} , \qquad 
\index{ey^-@${\bf e}^-$|main}%
{\bf e}^- 
:=\begin{pmatrix} 0&0\\1&0\end{pmatrix} ,\qquad  
\index{hy@${\bf h}$|main}%
{\bf h} 
:=\begin{pmatrix} 1&0\\0&-1\end{pmatrix}.\end{equation}
The triple $\{{\bf e}^+, {\bf e}^-, {\bf h}\}$ satisfies the commutation relations 
\begin{equation}\label{eqn:3.2}
[{\bf e}^+, {\bf e}^-]={\bf h}, \qquad [{\bf h}, {\bf e}^+]=2{\bf e}^+,\qquad  [{\bf h}, {\bf e}^-]=-2{\bf e}^-. 
\end{equation}
\begin{defe}\label{def:sl2}
An $\mathfrak{sl}_2$ triple is 
a triple of non-zero elements in a Lie algebra 
satisfying the same relation with \eqref{eqn:3.2}.
\end{defe}

We recall from Section \ref{sec:2}
that $\Delta_k$ is the Dunkl Laplacian 
associated with 
a multiplicity function $k$ on the root system, 
and that $\langle k \rangle$ is the index defined in \eqref{eqn:2.6}. 
For  a non-zero complex parameter $a$,
we introduce the following differential-difference operators on $\R^N$:
\begin{equation}\label{eqn:3.3}
\index{Eb_{k,a}^+@$\EE_{k,a}^+$|main}%
\EE_{k,a}^+ 
:={i\over a} \Vert x\Vert^a,\qquad 
\index{Eb^-_{k,a}@$\EE^-_{k,a}$|main}%
\EE^-_{k,a} 
:={i\over a}\Vert x\Vert^{2-a} \Delta_k,\qquad 
\index{Hb_{k,a}@$\HH_{k,a}$|main}%
\HH_{k,a} 
:={{N+2\langle k\rangle +a-2}\over a}
+{2\over a} \sum_{i=1}^N x_i \partial_i .\end{equation}
The point of the definition is:
\begin{thm}\label{thm:TDS}
The operators $\EE^+_{k,a},$ $\EE^-_{k,a}$ and $\HH_{k,a}$ form an $\mathfrak{sl}_2$ triple for
any complex number $a \ne 0$ and any multiplicity function $k$.
\end{thm} 

\begin{proof}[Proof of Theorem \ref{thm:TDS}]
The operator $\EE_{k,a}^+$ is homogeneous of degree $a$,
and $\EE_{k,a}^-$ is of degree $(2-a)-2=-a$ by Lemma \ref{lem:DLap}
1).
Let $E = \sum_{j=1}^N x_j\partial_j$ be the Euler operator as in
\eqref{eqn:Eu}. 
Since $\HH_{k,a}$ is of the form $\frac{2}{a}E+\text{constant}$,
the identity 
$[\HH_{k,a},\EE_{k,a}^\pm] = \pm2\EE_{k,a}^\pm$
is now clear.

To see
$[\EE_{k,a}^+,\EE_{k,a}^-]=\HH_{k,a}$,
we apply Lemma \ref{lem:DLap} 3) to the function
$\psi(r) = r$.
Then we get 
\begin{equation}\label{eqn:3.4}
\Delta_k\circ \Vert x\Vert ^a-\Vert x\Vert^a\Delta_k=a(N+2\langle k\rangle+a-2)\Vert x\Vert^{a-2} +2a\Vert x\Vert^{a-2} E.
\end{equation}
Composing the multiplication operator $\Vert x\Vert^{2-a}$,
we have
$$\Vert x\Vert^{2-a}\Delta_k\circ \Vert x\Vert^a-\Vert x\Vert^{2} \Delta_k=a(N+2\langle k\rangle+a-2)+2a E.$$
In view of the definition \eqref{eqn:3.3},
this means $ [\EE_{k,a}^+,\EE_{k,a}^-] =\HH_{k,a}.$ 

Hence, Theorem \ref{thm:TDS} is proved.
\end{proof}

\begin{rem}\label{rem:3.3} 
Theorem \ref{thm:TDS} for particular cases was previously known.
\begin{itemize}
\item[(1)] For  $a=2$ and $k \equiv 0$, 
$\{ \EE_{k,a}^+, \EE_{k,a}^-, \HH_{k,a} \}$ is  the classical harmonic
$\mathfrak{sl}_2$ triple $\{{i\over 2} \Vert
x\Vert^2, {i\over 2} \Delta, {N\over 2}+\sum_i x_i\partial_i\}$.
This $\mathfrak{sl}_2$ triple was used in the analysis of
 the Schr\"{o}dinger model of the 
Weil representation of the metaplectic group $Mp(N, \R)$
(see Howe \cite{H2}, Howe--Tan \cite{HT}).
\item[(2)]
For $a=2$ and $k > 0$,
Theorem \ref{thm:TDS} was proved in Heckman \cite[Theorem 3.3]{H}.
\item[(3)] For $a=1$ and $k \equiv 0$, 
 $\{\EE^+_{k,a},\EE^-_{k,a},\HH_{k,a}\}$ is the $\mathfrak{sl}_2$ triple introduced
in Kobayashi and Mano \cite{KM05,KM} where the authors studied the
$L^2$-model of the minimal representation of the double covering group of $SO_0(N+1,2).$
(To be more precise, the formulas in \cite{KM} are given for the
$\mathfrak{sl}_2$ 
triple for $\{2\EE^+_{k,a},\frac{1}{2}\EE^-_{k,a},\HH_{k,a}\}$ in our notation.)
\item[(4)] 
For $k\equiv0$, the deformation parameter $a$ was also considered in Mano
\cite{Mano}. 
\end{itemize}
\end{rem}

The differential-difference operators \eqref{eqn:3.3}
stabilize $C^\infty(\R^N \setminus \{0\})$,
the space of (complex valued) smooth functions on $\R^N \setminus \{0\}$.
Thus, for each non-zero complex number $a$ and each 
multiplicity function $k$ on the root system,
we can define an $\R$-linear map
\begin{equation}\label{eqn:omegaka}
\index{*zomega_{k,a}@$\omega_{k,a}$|main}%
\omega_{k,a}
: \mathfrak{sl}(2,\R) \to \operatorname{End}(C^\infty(\R^N \setminus \{0\}))
\end{equation}
by setting 
\begin{equation}\label{eqn:3.5}
\omega_{k,a}({\bf h})=\HH_{k,a},\qquad \omega_{k,a}({\bf
e}^+)=\EE_{k,a}^+, \qquad \omega_{k,a}({\bf e}^-)=\EE_{k,a}^-.\end{equation}
Then, Theorem \ref{thm:TDS} implies that 
$\omega_{k,a}$ is a Lie algebra homomorphism. 

We denote by $U(\mathfrak{sl}(2,\C))$ the universal enveloping algebra of the complex
Lie algebra $\mathfrak{sl}(2,\C) \simeq \mathfrak{sl}(2,\R)
\otimes_\R \C$.
Then, we can extend \eqref{eqn:omegaka} to a $\C$-algebra homomorphism (by the
same symbol)
$$
\omega_{k,a}: U(\mathfrak{sl}(2,\C)) \to \operatorname{End}(C^\infty(\R^N \setminus \{0\})).
$$

We use the letter $L$ to denote by the left regular representation of the
Coxeter group $\CC$ on $C^\infty(\R^N \setminus \{0\})$. 
\begin{lem}\label{lem:Lomega}
The two actions $L$ of the Coxeter group $\CC$ 
and $\omega_{k,a}$ of the Lie algebra $\mathfrak{sl}(2,\R)$ commute. 
\end{lem}
\begin{proof}
Obviously, $L(h)$ commutes with the multiplication operator
$\EE_{k,a}^+ = \frac{i}{a}\Vert x\Vert^a$.
As we saw in \eqref{eqn:2.4},
$L(h)$ commutes with the Dunkl Laplacian. 
Hence, it commutes also with $\EE_{k,a}^-$.
Finally, the commutation relation  $ [\EE_{k,a}^+, \EE_{k,a}^-]=\HH_{k,a}$ implies  $L(h)\circ \HH_{k,a} =\HH_{k,a} \circ L(h) $.
\end{proof}

We consider the following unitary matrix 
\begin{equation}\label{eqn:CTc}
c:=\frac{1}{\sqrt{2}}\begin{pmatrix} -i& -1 \\  1&i\end{pmatrix}. 
\end{equation}
We set
$$ 
\s\u(1,1) := \{ X \in \mathfrak{sl}(2,\C):
X^* \begin{pmatrix}1&0\\0&-1\end{pmatrix}
    + \begin{pmatrix}1&0\\0&-1\end{pmatrix} X = 0\},
$$
 another real form of $\s\l(2,\C)$.
Then, $\operatorname{Ad}(c)$ induces a Lie algebra isomorphism (the Cayley
transform)
\[
\operatorname{Ad}(c): \mathfrak{sl}(2,\R)
\stackrel{\sim}{\to} \mathfrak{su}(1,1).
\]
We set
\begin{subequations}\label{3.8}
  \renewcommand{\theequation}{\theparentequation\ \alph{equation}}%
\begin{alignat}{4}
& \index{kb@$\protect{\bf k}$|main}%
{\bf k} &&:= \operatorname{Ad}(c){\bf h}
          &&= i \begin{pmatrix}0 &-1 \\ 1 &0 \end{pmatrix}
          &&= \frac{1}{i} \, ({\bf e}^+ - {\bf e}^-),
\label{3.8 a}
\\
& \index{nb+@$\protect{\bf n}^+$|main}%
{\bf n}^+  &&:= \operatorname{Ad}(c){\bf e}^+ 
           &&= \frac{1}{2} \begin{pmatrix} i &-1 \\ -1 &-i \end{pmatrix}
           && = \frac{1}{2i} (- {\bf h} + \frac{1}{i}{\bf e}^+ + \frac{1}{i}{\bf e^-}),
\\
& \index{nb-@$\protect{\bf n}^-$|main}%
{\bf n}^-  &&:= \operatorname{Ad}(c){\bf e}^- 
           &&= \frac{1}{2} \begin{pmatrix} -i &-1 \\ -1 &i \end{pmatrix}
           &&= \frac{1}{2i} ({\bf h} + \frac{1}{i}{\bf e}^+ + \frac{1}{i}{\bf e^-}).
\end{alignat}
\end{subequations}

Correspondingly to (\ref{3.8} a -- c), the Cayley transform of the operators \eqref{eqn:3.5}
amounts to: 
\begin{subequations}\label{eqn:knn} 
  \renewcommand{\theequation}{\theparentequation\ \alph{equation}}%
\begin{align}
\label{eqn:knn a}
\index{Hb_{k,a}@$\protect\widetilde{\protect\mathbb{H}}_{k,a}$|main}%
\widetilde{\HH}_{k,a}
&:=\omega_{k,a}({\bf k})={{\Vert x\Vert^a-\Vert x\Vert^{2-a}\Delta_k}\over a} 
  = -\frac{1}{a}\Delta_{k,a}, 
\\
\label{eqn:knn b}
\index{Eb_{k,a}^+@$\protect\widetilde{\protect\mathbb{E}}_{k,a}^+$|main}%
\widetilde{\EE}_{k,a}^+ 
&:=\omega_{k,a}({\bf n}^+)=i{{ 2E+(N+2\langle k\rangle+a-2)-\Vert x\Vert^{2-a}\Delta_k -{{\Vert x\Vert^a} }}\over {2a}},\\
\label{eqn:knn c}
\index{Eb_{k,a}^-@$\protect\widetilde{\protect\mathbb{E}}_{k,a}^-$|main}%
\widetilde{\EE}_{k,a}^- 
&:=\omega_{k,a}({\bf n}^-)=-i{{ 2E+(N+2\langle k\rangle+a-2)+\Vert x\Vert^{2-a}\Delta_k +{{\Vert x\Vert^a} }}\over {2a}}. 
\end{align}
\end{subequations}
Here,
$
E 
=\sum_{i=1}^N x_i\partial_i
$
is the Euler operator.

Since $\operatorname{Ad}(c)$ gives a Lie algebra isomorphism,
$\{\widetilde{\EE}_{k,a}^+,\widetilde{\EE}_{k,a}^-,\widetilde{\HH}_{k,a}\}$
also forms an $\mathfrak{sl}_2$ triple of differential-difference operators.
Putting $\nu=\pm1$ in Lemma \ref{lem:3.4},
we get another expression of the triple $\{\widetilde \EE_{k,a}^+, \widetilde \EE_{k,a}^-,
\widetilde \HH_{k,a}\}$ as follows:
\begin{lem}\label{lem:3.5}
Let\/ $\widetilde{\EE}_{k,a}^+$, $\widetilde{\EE}_{k,a}^-$, and
$\widetilde{\HH}_{k,a}$ be as in\/ {\upshape(\ref{eqn:knn} a, b, c)}.
Then, we have:
\begin{subequations}
  \renewcommand{\theequation}{\theparentequation\ \alph{equation}}%
\begin{align}
&\label{3.13 a}
\widetilde \EE_{k,a}^+=\omega_{k,a}({\bf n}^+)=-{i\over {2a}}e^{ {{\Vert x\Vert^a}\over a}}\circ\Vert x\Vert^{2-a}\Delta_k\circ e^{-  {{\Vert x\Vert^a}\over a}},\\
&\label{3.13 b}
\widetilde \EE_{k,a}^-=\omega_{k,a}({\bf n}^-)= -{i\over{2a}} e^{- {{\Vert x\Vert^a}\over a}}\circ\Vert x\Vert^{2-a}\Delta_k\circ e^{ {{\Vert x\Vert^a}\over a}},\\
&\label{3.13 c}
\widetilde \HH_{k,a}=\omega_{k,a}({\bf k})= e^{- {{\Vert x\Vert^a}\over a}}\circ\left( \HH_{k,a}-{{\Vert x\Vert^{2-a}\Delta_k}\over a}\right) \circ e^{  {{\Vert x\Vert^a}\over a}}
\\
&\phantom{\widetilde \HH_{k,a}=\omega_{k,a}({\bf k})}
= \frac{1}{a} e^{-\frac{\Vert x\Vert^a}{a}} \circ
  \left( (N+2\langle k\rangle+a-2)+2E-\Vert x\Vert^{2-a}\Delta_k \right)
  \circ e^{\frac{\Vert x\Vert^a}{a}}.
\nonumber
\end{align}
\end{subequations}
\end{lem}

\subsection{Differential-difference operators in the polar coordinate}
\label{subsec:polar}
~\newline
In this subsection, 
we rewrite the differential-difference operators introduced in Section
\ref{subsec:3.1} by means of the polar coordinate.

We set
\begin{equation}\label{3.17}
\index{*lambda_{k,a,m}@$\lambda_{k,a,m}$|main}%
\lambda_{k,a,m} 
:={{2m+2\langle k\rangle+N-2}\over a}.\end{equation}
We begin with the following lemma.
\begin{lem}\label{lem:3.8}
 Retain the notation of Section \ref{subsec:DLap}. 
For all $\psi \in C^\infty(\R_+)$ and $p\in \Hkm(\R^N),$  we have 
\begin{align}\label{3.19}
\HH_{k,a}\Big( p(x) \psi(\Vert x\Vert^a)\Big)
&=\Big\{ (\lambda_{k,a,m}+1)
\psi(\Vert x\Vert ^a )+2\Vert x\Vert^a\psi'(\Vert x\Vert^a)\Big\}
p(x), 
\\
\label{3.20}
\Delta_k\Big( p(x) \psi(\Vert x\Vert^a)\Big)
&=\Big\{a^2(\lambda_{k,a,m}+1)
 \Vert x\Vert ^{a-2} \psi'(\Vert x\Vert ^a)+a^2\Vert x\Vert^{2a-2} \psi''(\Vert x\Vert^a)\Big\} p(x).
\end{align}
\end{lem} 
\begin{proof} The first statement is straightforward because
the Euler operator $E$ is of the form $r\frac{\partial}{\partial r}$ in the polar coordinates
$x = r\omega$.
To see the second statement, 
we apply \eqref{eqn:Dkpsi} to $p(x)$.
Since $Ep=mp$ and $\Delta_k p=0$,
we get
the desired formula \eqref{3.20}.
\end{proof}

We consider the following linear operator:
\begin{equation}\label{eqn:Tkap}
T_a:
C^\infty(\R^N)\otimes C^\infty(\R_+) \to C^\infty(\R^N \setminus \{0\}),
\ \ (p,\psi) \mapsto p(x) \psi(\Vert x\Vert^a)
\end{equation}

\begin{lem}\label{prop:slreduction}
Via the linear map $T_a$,
the operators $\HH_{k,a}$, $\EE_{k,a}^+$, and $\EE_{k,a}^-$
(see \eqref{eqn:3.3}) take
 the following forms
on
$\Hkm(\R^N) \otimes C^\infty(\R_+)$:
\begin{subequations}
  \renewcommand{\theequation}{\theparentequation\ \alph{equation}}%
\begin{alignat}{2}
&\HH_{k,a} \circ T_a
&&= T_a \circ \Big(\operatorname{id} \otimes
\Bigl( 2r \frac{d}{dr} + (\lambda_{k,a,m}+1) \Bigr) \Big)
\\
&\EE_{k,a}^+ \circ T_a
&&= T_a \circ \Big(\operatorname{id} \otimes \frac{i}{a} r\Big)
\\
& \EE_{k,a}^- \circ T_a
&&= T_a \circ \Big(\operatorname{id} \otimes
ai \Bigl( r \frac{d^2}{dr^2} + (\lambda_{k,a,m}+1) \frac{d}{dr} \Bigr) \Big)
\end{alignat}
\end{subequations}
\end{lem}
\begin{proof}
Clear from Lemma \ref{lem:3.8} and the definition \eqref{eqn:3.3} of $\HH_{k,a}$,
$\EE_{k,a}^+$, and $\EE_{k,a}^-$.
\end{proof}
The point of Lemma \ref{prop:slreduction} is that the operators
$\HH_{k,a}$, $\EE_{k,a}^+$, and $\EE_{k,a}^-$ act only on the radial
part $\psi$ when applied to those
functions $p(x)\psi(\Vert x\Vert^a)$ for 
$p \in \Hkm(\R^N)$.

For $a>0$, we define an endomorphism of $C^\infty(\R_+)$ by
\[
U_a: C^\infty(\R_+) \stackrel{\sim}{\to} C^\infty(\R_+),
\  g(t) \mapsto (U_a g)(r)
:= \exp \Bigl( -\frac{1}{a} r \Bigr) g \Bigl( \frac{2}{a} r \Bigr).
\]
Clearly, $U_a$ is invertible.
Composing with $T_a$ (see \eqref{eqn:Tkap}),
we define the following
 linear operator $S_a$ by
\[
S_a :=T_a  \circ (\operatorname{id} \otimes U_a).
\]
That is, 
$
S_a: C^\infty(\R^N) \otimes C^\infty(\R_+) \to 
C^\infty(\R^N \setminus \{0\})
$
is given by
\begin{equation}\label{eqn:expa}
S_a(p\otimes g)(x)
:= p(x) \exp\Bigl(-\frac{1}{a} \Vert x\Vert^a\Bigr) 
g\Bigl(\frac{2}{a} \Vert x\Vert^a \Bigr).
\end{equation}
We set
\begin{equation}\label{eqn:LagP}
P_{t,\lambda} := t \frac{d^2}{dt^2} + (\lambda_{k,a,m}+1-t) \frac{d}{dt}.
\end{equation}
Here, $\lambda$ stands for $\lambda_{k,a,m}.$ Then
Lemma \ref{prop:slreduction} can be formulated as follows:
\begin{lem}\label{lem:pHE}
Via the map $S_a$, 
the operators $\HH_{k,a}$, $\EE_{k,a}^+$, and $\EE_{k,a}^-$ take
the following forms on $\Hkm(\R^N) \otimes C^\infty(\R_+)$:
\begin{alignat*}{2}
& \HH_{k,a} \circ S_a
&&= S_a \circ \Bigl(\operatorname{id} \otimes
 \Bigl(2t\frac{d}{dt}+(\lambda_{k,a,m}+1-t)\Bigr)\Bigr) ,
\\
& \EE_{k,a}^+ \circ S_a 
&&= S_a \circ \Bigl(\operatorname{id} \otimes
  \frac{i}{2} t\Bigr),
\\
& \EE_{k,a}^- \circ S_a 
&&= S_a \circ \Bigl(\operatorname{id} \otimes
  i\Bigl(2P_{t,\lambda} + \frac{t}{2} - \lambda_{k,a,m}-1\Bigr)\Bigr).
\end{alignat*}
\end{lem}
\begin{proof}
Immediate from Lemma \ref{prop:slreduction} and the following
relations:
\begin{alignat*}{2}
&U_a^{-1} \circ \frac{d}{dr} \circ U_a
&&= \frac{2}{a} \left( \frac{d}{dt} - \frac{1}{2} \right),
\\
&U_a^{-1} \circ r \circ U_a
&&= \frac{a}{2} t.
\end{alignat*}
\end{proof}

Similarly, by using (\ref{3.8} a--c),
the actions of $\widetilde{\HH}_{k,a}$,
$\widetilde{\EE}_{k,a}^+$,
and $\widetilde{\EE}_{k,a}^-$ (see (\ref{eqn:knn} a--c))
are given as follows:
\begin{lem}\label{lem:HEE}
Let $P_{t,\lambda}$ be as in \eqref{eqn:LagP}.
Then, through the linear map
 $S_a$ (see \eqref{eqn:expa}),
$\widetilde{\HH}_{k,a}$, $\widetilde{\EE}_{k,a}^+$, and $\widetilde{\EE}_{k,a}^-$
take the following forms on $\Hkm(\R^N) \otimes C^\infty(\R_+)$:
\begin{alignat*}{2}
&  \widetilde{\HH}_{k,a} \circ S_a
&&=S_a \circ \Big(\operatorname{id}\otimes \Big(-2P_{t,\lambda}+ \lambda_{k,a,m}+1 \Big)\Big),
\\
&  \widetilde{\EE}_{k,a}^+ \circ S_a
&&=S_a \circ \Bigl(\operatorname{id}\otimes \Big(-i\Bigl(P_{t,\lambda} - t\frac{d}{dt} + t - \lambda_{k,a,m} -1\Bigr)\Bigr)\Big),
\\
&  \widetilde{\EE}_{k,a}^- \circ S_a
&&=S_a \circ \Bigl(\operatorname{id}\otimes \Big(-i\Bigl(P_{t,\lambda} + t\frac{d}{dt}\Bigr)\Bigr)\Big).
\end{alignat*}
\end{lem}

\subsection{Laguerre polynomials revisited}
\label{sec:Lag}
~\newline
In this subsection,
after a brief summary on the (classical) Laguerre polynomials we give
a `non-standard' representation of them in terms of the one parameter group with infinitesimal
generator
$t \frac{d^2}{dt^2} + (\lambda+1) \frac{d}{dt}$
(see Proposition \ref{prop:homoLag}).

For a complex number $\lambda\in \C$ such that $\Re\lambda>-1$,
we write $L_\ell^{(\lambda)} $ for the Laguerre polynomial defined by
$$ 
\index{L_\ell^{(\lambda)}(t)@$L_\ell^{(\lambda)}(t)$|main}%
L_\ell^{(\lambda)}(t) 
:={{(\lambda+1)_\ell}\over {\ell!}} \sum_{j=0}^\ell {{(-\ell)_j}\over {(\lambda+1)_j}} {{t^j}\over {j!}}
= \sum_{j=0}^\ell
  \frac{(-1)^j \Gamma(\lambda+\ell+1)}
       {(\ell-j)! \Gamma(\lambda+j+1)}
  \, {{t^j}\over {j!}}.
$$
Here, $(a)_m:=a(a+1)\cdots(a+m-1)$ is the Pochhammer symbol.

We list some standard properties of Laguerre polynomials that we shall
use in this article.
\begin{fact}[{see \cite[\S6.5]{AAR}}]\label{fact:Lag}
Suppose $\Re\lambda>-1$.
\begin{enumerate}[\upshape 1)]
\item  
$L_\ell^{(\lambda)}(t)$ is the unique polynomial of degree $\ell$
satisfying the Laguerre differential equation
\begin{equation}\label{eqn:LODE}
\Bigl( t \frac{d^2}{dt^2} + (\lambda+1-t) \frac{d}{dt} + \ell\Bigr) f(t)
= 0
\end{equation}
and
\begin{equation}\label{eqn:Ltop}
f^{(\ell)}(0) = (-1)^{\ell}.
\end{equation}
\item  
(recurrence relation)
\begin{subequations}
  \renewcommand{\theequation}{\theparentequation\ \alph{equation}}%
\begin{align}
\label{3.22 a}
&(\ell+t\frac{d}{dt}-t+\lambda+1)L_\ell^{(\lambda)}(t)
= (\ell+1)L_{\ell+1}^{(\lambda)}(t),
\\
&(\ell-t\frac{d}{dt}) L_\ell^{(\lambda)}(t)
= (\ell+\lambda)L_{\ell-1}^{(\lambda)}(t).
\label{3.22 b}
\end{align}
\end{subequations}
\item  
(orthogonality relation)
\begin{equation}\label{3.18}
\int_0^\infty
L_\ell^{(\lambda)}(t)L_s^{(\lambda)}(t)t^\lambda e^{-t} dt=\delta_{\ell
s} {{\Gamma(\lambda+\ell+1)}\over {\Gamma(\ell+1)}}.
\end{equation}
\item  
(generating function)
\begin{equation} 
(1-r)^{-\lambda-1} \exp\Big({{rt}\over {r-1}}\Big) =\sum_{\ell=0}^\infty 
L_\ell^{(\lambda)}(t) r^\ell,\quad (|r|<1).
\end{equation}
\item  
$\{L_\ell^{(\lambda)}(t): \ell \in \N \}$
form an orthogonal basis in
$L^2(\R_+,t^\lambda e^{-t} dt)$ if $\lambda$ is real and $\lambda>-1$.
\end{enumerate}
\end{fact}

Finally,
we give a new representation of the Laguerre polynomial.

\begin{thm}\label{prop:homoLag}
For any $c\ne0$ and $\ell\in\N$,
\begin{equation}\label{eqn:expLag}
\exp \Bigl( -c\bigl(t\frac{d^2}{dt^2} + (\lambda+1) \frac{d}{dt}\bigr)\Bigr)
t^\ell
= (-c)^\ell \ell! L_\ell^{(\lambda)}\Bigl(\frac{t}{c}\Bigr).
\end{equation}
\end{thm}
Since the differential operator 
$$B_t := t \frac{d^2}{dt^2} + (\lambda+1) \frac{d}{dt}$$
is homogeneous of degree $-1$, namely,
$B_t = cB_x$ if $x=ct$,
it is sufficient to prove Theorem  \ref{prop:homoLag} in the case $c=1$.
We shall give two different proofs for this.

\begin{proof}[Proof 1]
We set
\[
A := t \frac{d}{dt} - \ell,
\quad
B := t \frac{d^2}{dt^2} + (\lambda+1) \frac{d}{dt}.
\]
It follows from $[A,B]=-B$ that
\[
AB^n = B^nA - nB^n
\]
for all $n \in \N$ by induction.
Then, by the Taylor expansion 
$e^{-B}= \sum_{n=0}^\infty \frac{(-1)^n}{n!} B^n$,
we get
\[
Ae^{-B} = e^{-B} A + B e^{-B}.
\]
Since $At^\ell = 0$,
we get $(B-A)(e^{-B}t^\ell) = 0$, namely, 
$e^{-B} t^\ell$ solves the Laguerre differential equation \eqref{eqn:LODE}.
On the other hand,
$e^{-B}t^\ell$ is clearly a polynomial of $t$ with top term $t^\ell$.
In view of \eqref{eqn:Ltop}, 
we have $e^{-B}t^\ell = (-1)^\ell \ell! L_\ell^{(\lambda)} (t)$.
\end{proof}

\begin{proof}[Proof 2]
A direct computation shows
\[
B t^\ell = \ell(\lambda+\ell) t^{\ell-1}.
\]
Therefore,
$B^j t^\ell = 0$ for $j>\ell$ and
\begin{align*}
e^{-B} t^\ell
&= \sum_{j=0}^\ell
   \frac{(-1)^j \ell(\ell-1)\cdots(\ell-j+1)(\lambda+\ell)(\lambda+\ell-1)
         \cdots (\lambda+\ell-j+1)}
        {j!}
  \, t^{\ell-j}
\\
&= \sum_{k=0}^\ell
   \frac{(-1)^{\ell+k} \ell! \Gamma(\lambda+\ell+1)}
        {(\ell-k)! \Gamma(\lambda+k+1)}
  \,  {{t^k}\over {k!}}
\\
&= (-1)^\ell \ell! L_\ell^{(\lambda)} (t).
\end{align*}
Hence, Theorem \ref{prop:homoLag} has been proved.
\end{proof}

\subsection{Construction of an orthonormal basis in
$L^2(\R^N,\vartheta_{k,a}(x)dx)$}
\label{sec:onb}
~\newline
We recall from \eqref{eqn:01} and \eqref{eqn:2.5} that the weight function
$\vartheta_{k,a}$ on $\R^N$ satisfies
$$
\vartheta_{k,a}(x) = \Vert x\Vert^{a-2}
\prod_{\alpha\in\mathcal{R}^+} \vert \langle \alpha,x \rangle \vert^{2k_\alpha}
= \Vert x \Vert^{a-2} \vartheta_k(x).
$$
Therefore, in the polar coordinates
$x=r\omega$ ($r>0$, $\omega\in S^{N-1}$),
we have
\begin{equation}\label{eqn:tpol}
\vartheta_{k,a}(x)dx
= r^{2\langle k\rangle+N+a-3} \vartheta_k(\omega)dr d\sigma(\omega),
\end{equation}
where $d\sigma(\omega)$ is the standard measure on the unit sphere.
Accordingly, we have a unitary isomorphism:
\begin{equation}\label{eqn:L2polar}
L^2(S^{N-1},\vartheta_k(\omega)d\sigma(\omega))
\mathbin{\widehat{\otimes}}
L^2(\R_+,r^{2\langle k\rangle+N+a-3} dr)
\stackrel{\sim}{\to}
L^2(\R^N, \vartheta_{k,a}(x)dx),
\end{equation}
where $\widehat{\otimes}$ stands for the Hilbert completion of the tensor
product space of two Hilbert spaces.

Combining \eqref{eqn:L2polar} with Fact \ref{fact:kH},
we get a direct sum decomposition of the Hilbert space:
\begin{equation}\label{eqn:L2RNsum}
\sideset{}{^\oplus}\sum_{m\in\N}
(\mathcal{H}_k^m(\R^N)|_{S^{N-1}})
\otimes L^2(\R_+,r^{2\langle k\rangle+N+a-3} dr)
\stackrel{\sim}{\to}
L^2(\R^N,\vartheta_{k,a}(x)dx).
\end{equation}

In this subsection, we demonstrate the irreducible decomposition theorem of the
$\mathfrak{sl}_2$ representation on (a dense
subspace of) $L^2(\R^N,\vartheta_{k,a}(x)dx)$ by using
\eqref{eqn:L2RNsum} and finding an orthogonal basis for
$L^2(\R_+,r^{2\langle k\rangle+N+a-3} dr)$.

For   $\ell,m\in \N$  and $p\in \Hkm(\R^N),$ we introduce the following
functions on $\R^N$:
\begin{equation}\label{eqn:PhiS}
\Phi_\ell^{(a)} (p,\cdot)
:= S_a (p \otimes L_\ell^{(\lambda_{k,a,m})}).
\end{equation}
Here, $S_a: C^\infty(\R^N) \otimes C^\infty(\R_+) \to C^\infty(\R^N \setminus \{0\})$ 
is a linear operator defined in \eqref{eqn:expa},
$\lambda_{k,a,m} = \frac{1}{a}(2m+2\langle k\rangle+N-2)$
(see \eqref{3.17}),
and $L_\ell^{(\lambda)}(t)$ is the Laguerre polynomial.
Hence, for
 $x=r\omega\in\R^N$ ($r>0$, $\omega\in S^{N-1}$),
we have
\begin{align}\label{3.16}
\index{*VPhi_{\ell}^{(a)}(p,x)@$\Phi_{ \ell}^{(a)}(p,x)$|main}%
\Phi_{ \ell}^{(a)}(p,x) 
={}&p(x) L_\ell^{(\lambda_{k,a,m})}\Big({2\over a}\Vert x\Vert^a\Big)
\exp\Bigl(-\frac{1}{a}\Vert x\Vert^a\Bigr)
\\
={}&p(\omega)r^mL_\ell^{(\lambda_{k,a,m})}
      \Bigl(\frac{2}{a}r^a\Bigr) \exp\Bigl(-\frac{1}{a}r^a\Bigr).
\nonumber
\end{align} 

We define the following vector space of functions on $\R^N$ by
\begin{equation}\label{eqn:Wka}
\index{W_{k,a}(\R^N)@$W_{k,a}(\R^N)$|main}%
W_{k,a}(\R^N) 
:= \mbox{$\C$-span} 
 \{\Phi_{\ell }^{(a)}(p,\cdot)\;|\; \ell\in \N, m\in \N, p\in
 \Hkm(\R^N)\}.
\end{equation}

\begin{prop}\label{prop:3 7}
Suppose  $k$ is a non-negative multiplicity function on the root system
$\mathcal{R}$ and $a>0$ such that
\begin{equation}\label{eqn:a2k}
a+2\langle k\rangle+N-2 > 0.
\end{equation}
Let $  \ell, s,{m}, {  n}\in \N,$  $p\in \Hkm(\R^N)$ and 
$q\in \mathcal H_{k}^n(\R^N)$.
\begin{enumerate}[\upshape 1)]
\item  
$\Phi_\ell^{(a)}(p,x) \in C(\R^N) \cap L^2(\R^N,\vartheta_{k,a}(x)dx)$.
\item  
$$
\int_{\R^N} \Phi_{\ell }^{(a)}(p,x)
\overline{\Phi_{s }^{(a)}(q,x)} \vartheta_{k,a}(x) dx 
=\delta_{m,n} \delta_{\ell,s} 
    {{a^{\lambda_{k,a,m}}\Gamma(\lambda_{k,a,m}+\ell+1)}\over {2^{1+\lambda_{k,a,m} }\Gamma(\ell+1)}}\int_{S^{N-1}} p(\omega) \overline{q(\omega)} \vartheta_k(\omega) d\sigma(\omega).
$$
\item  
$W_{k,a}(\R^N)$ is a dense subspace of
$L^2(\R^N,\vartheta_{k,a}(x)dx)$.
\end{enumerate}
\end{prop}
\begin{rem}
The special values of our
 functions $\Phi_\ell^{(a)}(p,x)$ have been used in various settings including:
\begin{alignat*}{2}
&a=2 &\qquad&\text{see \cite[\S3]{D08}},
\\
&k\equiv0, \  N=1 &\qquad&\text{see \cite{Kos}},
\\
&k\equiv0, \  a=1 &\qquad&\text{see \cite[\S3.2]{KM}}.
\end{alignat*}
\end{rem}
\begin{rem}
The condition \eqref{eqn:a2k} is automatically satisfied
for $a>0$ and a non-negative multiplicity $k$ if $N\ge2$.
\end{rem}
\begin{proof}

Our assumption \eqref{eqn:a2k} implies
\[
 \lambda_{k,a,m} > -1
\quad\text{for any $m\in\N$},
\]
and thus $\Phi_\ell^{(a)} (p,x)$ is continuous at $x=0$.
Therefore, it is a continuous function on
$x\in\R^N$ of exponential decay.
On the other hand,
we see from \eqref{eqn:tpol} that
the measure $\vartheta_{k,a}(x)dx$ 
 is locally integrable 
under our assumptions on $a$ and $k$.
Therefore, $\Phi_\ell^{(a)} (p,x) \in L^2(\R^N,\vartheta_{k,a}(x)dx)$.
Hence the first statement is proved.

To see the second and third statements,
we rewrite the left-hand side of the integral as
$$
\Big( \int_0^\infty L_\ell^{(\lambda_{k,a,m})}\big({2\over
a}r^a\big)L_s^{(\lambda_{k,a,n})}\big({2\over
a}r^a\big)\exp\Bigl(-\frac{2}{a}r^a\Bigr) r^{m+n+2\langle k\rangle +N+a-3} dr \Big)
\Big(\int_{S^{N-1}} p(\omega) \overline{q(\omega)} \vartheta_k(\omega)
d\sigma(\omega)\Big)$$   
in the polar coordinates $x=r\omega$.
Since $k$-harmonic polynomials of different degrees are orthogonal to
each other
(see Fact \ref{fact:kH}),  the integration over
$S^{N-1}$ vanishes if $m \not = n.$  

Suppose that $m= n$.
By changing the variable $t := \frac{2}{a} r^a$,
we see that the first integration amounts to
\begin{equation}\label{eqn:int1}
\frac{a^{\lambda_{k,a,m}}}{2^{1+\lambda_{k,a,m}}}
\int_0^\infty L_\ell^{(\lambda_{k,a,m})} (t) L_s^{(\lambda_{k,a,m})}(t)
 \,  t^{\lambda_{k,a,m}} e^{-t} dt.
\end{equation}
By the orthogonality relation \eqref{3.18},
we get
$$
\text{\eqref{eqn:int1}}
=
\delta_{\ell s}   {{a^{\lambda_{k,a,m}}\Gamma(\lambda_{k,a,m}+\ell+1)}\over {2^{1+\lambda_{k,a,m}}\Gamma(\ell+1)}}.$$
Hence, the second statement is proved.
The third statement follows from the completeness of the Laguerre
polynomials (see Fact \ref{fact:Lag} 4)).
\end{proof}
We pin down the following proposition which is already implicit in the
 proof of Proposition \ref{prop:3 7}:
\begin{prop}\label{prop:falr}
We fix $m\in\N$, $a>0$, and a multiplicity function $k$ satisfying
\[
2m+2\langle k\rangle+N+a-2 > 0.
\]
We set
\begin{equation}\label{eqn:4.4}
\index{fl_{\ell,m}^{(a)}(r)@$f_{\ell,m}^{(a)}(r)$|main}%
f_{\ell,m}^{(a)}(r) 
:= \Big( {{ 2^{\lambda_{k,a,m}+1} \Gamma(\ell+1)}\over
{a^{\lambda_{k,a,m}}\Gamma( \lambda_{k,a,m}+\ell+1)}}\Big)^{1/2}r^{m}
L_\ell^{(\lambda_{k,a,m})}\Big({2\over
a}r^a\Big)\exp(-\frac{1}{a}r^a)
\quad\text{for $\ell\in\N$}. 
\end{equation}
Then
$\{f_{\ell,m}^{(a)}(r): \ell\in\N\}$
forms an orthonormal basis in
$L^2(\R_+,r^{2\langle k\rangle+N+a-3}dr)$.
\end{prop}
\begin{rem}\label{rem:4.1}
Let $c_0, c_1, \ldots$ be a sequence of positive real numbers. Fix a
parameter $\alpha >0$. Dunkl \cite{D03} proves that the only possible orthogonal
sets
$\{
L_\ell^{(\alpha)}(c_\ell r) \exp( -{1\over 2} c_\ell r)
\}_{\ell=0}^\infty$
for the measure $r^{\alpha+ \mu} dr$ on $\R_+$, with $\mu \geq 0$, are the
two cases (1) $\mu=0$, $c_\ell=c_0$ for all $\ell$; (2) $\mu=1,$ $c_\ell=c_0
{{\alpha+1} \over { \alpha+2\ell +1} }.$
\end{rem}

For each   ${m}\in \N,$   
we take an
orthonormal basis 
$\{ h_j^{(m)}\}_{j\in J_{m}}$ 
of the space $\Hkm(\R^N)|_{S^{N-1}}$. 
Proposition \ref{prop:3 7} immediately yields the following statement. 
\begin{coro}\label{lem:5.17}
Suppose that $a>0$ and that the non-negative multiplicity function 
$k$ satisfies the inequality \eqref{eqn:a2k}.
 For $\ell,m\in \N $   and $j\in J_{m},$ 
we set
\[
\Phi_{\ell,m,j}^{(a)}(x)
:=  h_j^{(m)}
\Bigl( \frac{x}{\Vert x\Vert}\Bigr)
f_{\ell,m}^{(a)} (\Vert x\Vert).
\]
Then, 
the set $\Big\{\Phi_{\ell, {m}, j}^{(a)}\;|\; \ell\in \N, {m}\in \N,   j\in J_{m}\Big\}$ forms an orthonormal basis of $L^2\big(\R^N, \vartheta_{k,a}(x)dx\big).$
\end{coro}
\begin{rem}
A basis of $\Hkm(\R^N)$ is constructed in \cite[Corollary 5.1.13]{DX}.
\end{rem}

\subsection{$\mathfrak{sl}_2$ representation on $L^2(\R^N,
\vartheta_{k,a}(x)dx)$} 
\label{subsec:sl2}
~\newline
Now we are  ready to exhibit the action of the $\mathfrak{sl}_2$
triple $\{{\bf k}, {\bf n}^+, {\bf n}^-\}$ on the basis
$\Phi_\ell^{(a)}(p,\cdot)$ (see (\ref{3.8} a--c) and \eqref{3.16} for the definitions).
We recall from (\ref{eqn:knn} a--c) that
$\widetilde{\HH}_{k,a} = \omega_{k,a}({\bf k})$,
$\widetilde{\EE}_{k,a}^+ = \omega_{k,a}({\bf n}^+)$, and
$\widetilde{\EE}_{k,a}^- = \omega_{k,a}({\bf n}^-)$.

\begin{thm}\label{thm:3.9}
Let $W_{k,a}(\R^N)$ be the dense subspace of
$L^2(\R^N,\vartheta_{k,a}(x)dx)$ defined in \eqref{eqn:Wka}.
Then, $W_{k,a}(\R^N)$ is stable under the action of\/ $\s\l(2,\C)$.
More precisely, for each fixed $p \in \mathcal{H}_k^m(\R^N)$,
the action $\omega_{k,a}$  (see {\upshape(\ref{eqn:knn} a--c)}) is given as follows:
\begin{subequations}
  \renewcommand{\theequation}{\theparentequation\ \alph{equation}}%
\begin{align}
&\label{3.21 a}
\omega_{k,a}({\bf k})\Phi_{\ell }^{(a)}(p,x)=(2\ell+\lambda_{k,a,m}+1) \Phi_{\ell }^{(a)}(p,x),\\
&\label{3.21 b}
 \omega_{k,a}({\bf n}^+)\Phi_{\ell }^{(a)}(p,x)=i(\ell+1) \Phi_{\ell+1 }^{(a)}(p,x),\\
&\label{3.21 c}
\omega_{k,a}({\bf n}^-)\Phi_{\ell }^{(a)}(p,x)=i (\ell+\lambda_{k,a,m})  \Phi_{\ell-1 }^{(a)}(p,x), 
\end{align}
\end{subequations}
where $\Phi_\ell^{(a)} (p,x)$ is defined in \eqref{3.16} and
 $\lambda_{k,a,m}=(2m+2\langle k\rangle +N-2)/a$ (see
\eqref{3.17}).
We have used the convention $\Phi_{-1 }^{(a)}\equiv 0.$ 
\end{thm}
Theorem \ref{thm:3.9} may be visualized by
the diagram below.
We see that
for each fixed $k,a$, and $p \in \mathcal{H}_k^m(\R^N)$,
the operators $\omega_{k,a}({\bf n}^+)$ and $\omega_{k,a}({\bf n}^-)$ act as 
\textit{raising/lowering operators}.
\begin{figure}[H]
\renewcommand{\figurename}{Diagram}
\renewcommand{\thefigure}{\thesubsection}
$$
\xymatrix{
{\cdots}\kern-2.5em 
&{\circ} \ar@/^.7pc/[r] 
&{\circ} \ar@/^.7pc/[r]^{{\bf n}^+} \ar@/^.7pc/[l]
{\smash{\rlap{\raisebox{2.5ex}{$\kern-.7em\scriptstyle\ell-1$}}}}
&{\circ} \ar@/^.7pc/[r]^{{\bf n}^+} \ar@/^.7pc/[l]^{{\bf n}^-} 
{\smash{\rlap{\raisebox{2.5ex}{$\kern-.3em\scriptstyle\ell$}}}}
&{\circ} \ar@/^.7pc/[r] \ar@/^.7pc/[l]^{{\bf n}^-} 
{\smash{\rlap{\raisebox{2.5ex}{$\kern-.7em\scriptstyle\ell+1$}}}}
&{\circ} \ar@/^.7pc/[l] 
&\kern-2.5em {\cdots}
}
$$
\label{diagram:3.2}
\caption{}
\end{figure}
Here, the dots represent $\omega_{k,a}({\bf k})$ eigenvectors
$\Phi_\ell^{(a)}(p,x)$ arranged by increasing $\omega_{k,a}({\bf k})$
eigenvalues, from left to right.

\begin{proof}[Proof of Theorem \ref{thm:3.9}]
For simplicity, we 
 use the notation
$P_{t,\lambda} = t\frac{d^2}{dt^2}+(\lambda_{k,a,m}+1-t)\frac{d}{dt}$ as in
\eqref{eqn:LagP}, where $\lambda$ stands for $\lambda_{k,a,m}$. 
By the formula
$\Phi_\ell^{(a)}(p,\cdot) = S_a(p\otimes L_\ell^{(\lambda)})$
(see \eqref{eqn:PhiS}) and by Lemma \ref{lem:HEE},
it is sufficient to prove
\begin{subequations}\label{eqn:rec} 
  \renewcommand{\theequation}{\theparentequation\ \alph{equation}}%
\begin{align}
\label{eqn:rec a}
&(-2P_{t,\lambda} +(\lambda+1))L_\ell^{(\lambda)}
 = (2\ell+\lambda+1)L_\ell^{(\lambda)},
\\
\label{eqn:rec b}
&(-i(P_{t,\lambda}-t\frac{d}{dt}+t-\lambda-1))L_\ell^{(\lambda)}
 = i(\ell+1)L_{\ell+1}^{(\lambda)},
\\
\label{eqn:rec c}
&-i(P_{t,\lambda}+t\frac{d}{dt}) L_\ell^{(\lambda)}
 = i(\ell+\lambda)L_{\ell-1}^{(\lambda)}.
\end{align}
\end{subequations}
Since the Laguerre polynomial $L_\ell^{(\lambda)}(t)$ satisfies the
Laguerre differential equation
\[
P_{t,\lambda} L_\ell^{(\lambda)}(t) = -\ell L_\ell^{(\lambda)}(t)
\]
(see \eqref{eqn:LODE}),
the assertion \eqref{eqn:rec a} is now clear.
The assertions \eqref{eqn:rec b} and \eqref{eqn:rec c} are reduced to the recurrence
relations \eqref{3.22 a} and \eqref{3.22 b},
respectively.
\end{proof}

\begin{rem}
An alternative proof of \eqref{3.21 a} will be given in Section
\ref{sec:5.4} (see Remark \ref{rem:Thm3.19}).
\end{rem}

By using the orthonormal basis $\{f_{\ell,m}^{(a)}(r)\}$ (see
\eqref{eqn:4.4}), we normalize $\Phi_\ell^{(a)}(p,x)$ as
\begin{align}\label{eqn:nPhi}
\index{*VPhitilde_{\ell}^{(a)}(p,x)@$\protect\widetilde{\Phi}_{\ell}^{(a)}(p,x)$|main}%
\widetilde{\Phi}_{\ell}^{(a)}(p,x)
:={}&f_{\ell,m}^{(a)}(r)p(\omega)
\\
={}&\Big( {{ 2^{\lambda_{k,a,m}+1} \Gamma(\ell+1)}\over
{a^{\lambda_{k,a,m}}\Gamma( \lambda_{k,a,m}+\ell+1)}}\Big)^{\frac{1}{2}}
\Phi_\ell^{(a)}(p,x)
\nonumber
\end{align}
for $x=r\omega$ ($r>0$, $\omega\in S^{N-1}$).
Then, 
Theorem \ref{thm:3.9} is reformulated as follows:
\begin{thm}\label{thm:onbPhi}
For any $p\in \Hkm(\R^N)$,
we have
\begin{subequations}
  \renewcommand{\theequation}{\theparentequation\ \alph{equation}}%
\begin{align}
\omega_{k,a}({\bf k}) \widetilde{\Phi}_\ell^{(a)}(p,x)
&= (2\ell+\lambda_{k,a,m}+1) \widetilde{\Phi}_\ell^{(a)}(p,x),
\\
\omega_{k,a}({\bf n}^+) \widetilde{\Phi}_\ell^{(a)}(p,x)
&= i\sqrt{(\ell+1)(\lambda_{k,a,m}+\ell+1)} \, \widetilde{\Phi}_{\ell+1}^{(a)}(p,x),
\\
\omega_{k,a}({\bf n}^-) \widetilde{\Phi}_\ell^{(a)}(p,x)
&= i\sqrt{(\lambda_{k,a,m}+\ell)\ell} \, \widetilde{\Phi}_{\ell-1}^{(a)}(p,x).
\end{align}
\end{subequations}
\end{thm}

We recall that an
operator $T$ densely defined on a Hilbert space is called essentially
self-adjoint, if it is symmetric and its closure is a self-adjoint operator. 
\begin{coro}\label{cor:kaLap}
Let $a>0$ and $k$ be a  non-negative multiplicity function satisfying \eqref{eqn:a2k}.
\begin{enumerate}[\upshape 1)]
\item  
The differential-difference operator 
$\Delta_{k,a} = \Vert x\Vert^{2-a} \Delta_k - \Vert x\Vert^a$
is an essentially self-adjoint operator on
$L^2(\R^N,\vartheta_{k,a}(x)dx)$.
\item  
There is no continuous spectrum of $\Delta_{k,a}$.
\item  
The set of discrete spectra of $-\Delta_{k,a}$ is given by
\begin{alignat*}{3}
&\{2a\ell+2m+2\langle k\rangle+N-2+a
&&:
   \ell,m\in \N \}
&&\quad (N \ge 2),
\\
&\{2a\ell+2\langle k\rangle+a\pm1
&&:
   \ell\in \N \}
&&\quad (N=1).
\end{alignat*}
\end{enumerate}
\end{coro}

\begin{proof}
In light of the formula \eqref{eqn:knn a}
\[
\Delta_{k,a} = -a \omega_{k,a} ({\bf k}),
\]
the eigenvalues of $\Delta_{k,a}$ are read from Theorem \ref{thm:3.9}.
Since $W_{k,a}(\R^N)$ is dense in
$L^2(\R^N,\vartheta_{k,a}(x)dx)$
(see Proposition \ref{prop:3 7}),
the remaining statement of Corollary \ref{cor:kaLap} is straightforward
from the following fact.
\end{proof}

\begin{fact}\label{lem:3.11}
 Let $T$ be a symmetric operator on a Hilbert space 
$\mathcal H$ with domain $\DD(T)$, and let $\{f_n\}_n$ be a complete orthogonal set in $\mathcal H.$ If each $f_n\in \DD(T)$ and there exists $\mu_n\in \R$ such that $T f_n=\mu_n f_n,$ for every $n,$ then $T$ is essentially self-adjoint. 
\end{fact}

\begin{rem}
We shall see in Theorem \ref{thm:3.12} that the action of\/
$\mathfrak{sl}(2,\R)$ in Theorem \ref{thm:onbPhi} lifts to a unitary
representation of the universal covering group $\widetilde{SL(2,\R)}$
and that Corollary \ref{cor:kaLap}
1) is a special case of the general theory of discretely decomposable 
$(\mathfrak{g}_{\C},K)$-modules (see \cite{Kdisc,Kaspm97}).
\end{rem}

\subsection{Discretely decomposable representations}
\label{subsec:DDR}
\hfill\break
Theorem \ref{thm:3.9} asserts that $W_{k,a}(\R^N)$ is an
$\mathfrak{sl}(2,\C)$-invariant, dense subspace in
$L^2(\R^N,\vartheta_{k,a}(x)dx)$.
For $N>1$, this is a `huge' representation in
the sense that it contains an infinitely many inequivalent irreducible
representations of $\mathfrak{sl}(2,\C)$.

By a theorem of Harish-Chandra,
Lepowsky and Rader,
any irreducible,
infinitesimally unitary $(\g_\C,K)$-module is the underlying
$(\g_\C,K)$-module of a (unique) irreducible unitary representation of $G$
(see \cite[Theorem 0.6]{KN}).
This result was generalized to a discretely decomposable
 $(\g_\C,K)$-modules by the second-named author
(see \cite[Theorem 2.7]{Kaspm97}).

In this section,
we discuss the meaning of 
Theorem \ref{thm:3.9} from the point of view  of
discretely decomposable  representations.

We begin with a general setting.
Let $G$ be a semisimple Lie group,
and $K$ a maximal compact subgroup of $G$ (modulo the center of $G$).
We write $\mathfrak{g}$ for the Lie algebra of $G$,
and $\mathfrak{g}_\C$ for its complexification.
The following notion singles out an algebraic property of unitary
representations that split into irreducible representations
without continuous spectra.
\begin{defe}\label{def:DDR}
Let $(\varpi,X)$ be a $(\mathfrak{g}_\C,K)$-module.
\begin{enumerate}[\upshape (1)]
\item 
{\upshape(\cite[Part I, \S1]{Kdisc})}\enspace
We say $\varpi$ is $K$-admissible if\/
$\dim\operatorname{Hom}_K(\tau,\varpi) < \infty$ for any $\tau \in \widehat{K}$.
\item 
{\upshape (\cite[Part III, Definition 1.1]{Kdisc})}\enspace
We say $\varpi$ is a discretely decomposable if there exist a sequence of\/
$(\mathfrak{g}_\C,K)$-modules $X_j$ such that
\begin{align*}
&\{0\} = X_0 \subset X_1 \subset X_2 \subset \cdots,
\quad
X = \bigcup_{j=0}^\infty X_j,
\\
&\text{$X_j/X_{j-1}$ is of finite length as a
$(\mathfrak{g}_\C,K)$-module for $j=1,2,\dotsc.$}
\end{align*}
\item 
We say $\varpi$ is infinitesimally unitarizable if there exists a
Hermitian inner product $(\ ,\ )$ on $X$ such that
\[
(\varpi(Y)u,v) = -(u,\varpi(Y)v)
\quad\text{for any $Y \in \g$, and any $u,v\in X$}.
\]
\end{enumerate}
\end{defe}
We collect some basic results on discretely decomposable
 $(\g_\C,K)$-modules: 
\begin{fact}[see {\cite{Kdisc,Kaspm97}}]\label{fact:deco}
Let $(\varpi,X)$ be a $(\g_C,K)$-module.
\begin{enumerate}[\upshape 1)]
\item  
If $\varpi$ is $K$-admissible,
then $\varpi$ is discretely decomposable as a $(\mathfrak{g}_\C,K)$-modules.
\item  
Suppose $\varpi$ is discretely decomposable as a
$(\mathfrak{g}_\C,K)$-module. 
If $\varpi$ is infinitesimally unitarizable, then $\varpi$ is isomorphic to
an algebraic direct sum of irreducible $(\mathfrak{g}_\C,K)$-modules.
\item  
Any discretely decomposable, infinitesimally unitary
$(\g_\C,K)$-module is the underlying $(\g_\C,K)$-module of a unitary
representation of $G$.
Furthermore, such a unitary representation is unique.
\end{enumerate}
\end{fact}

We shall apply this concept to the specific situation where
$\mathfrak{g}=\mathfrak{sl}(2,\R)$ and $G$ is the universal covering
group $\widetilde{SL(2,\R)}$ of $SL(2,\R)$.

We recall from \eqref{3.8 a} that
$$
{\bf k} = i \begin{pmatrix}0&-1\\1&0\end{pmatrix}
= i \, ({\bf e}^- - {\bf e}^+) \in \mathfrak{sl}(2,\C).
$$
Let $\k :=\R({\bf e}^--{\bf e}^+)=i\,\R {\bf k}$ and $K$ be the subgroup
of $G$ with Lie algebra  $\k.$  Since $G$ is
taken to be simply connected, the exponential map 
$$
\R\rightarrow K, \quad  t\mapsto \Exp(it \/ {\bf k})
$$ 
is a diffeomorphism.  

For $z \in i \, \R$, we set
\begin{equation}\label{eqn:gammaz}
\index{*cgamma_z@$\gamma_z$|main}%
\gamma_z 
:= \Exp(-z{\bf k}) 
= \Exp \begin{pmatrix}0 &iz \\ -iz &0 \end{pmatrix} \in K.
\end{equation}
Since $\{ {\bf k}, {\bf n}^+, {\bf n}^- \}$ forms an $\mathfrak{sl}_2$
triple, we have
\[
\operatorname{Ad}(\gamma_z){\bf n}^+ = e^{-2z}{\bf n}^+,
\quad
\operatorname{Ad}(\gamma_z){\bf n}^- = e^{2z}{\bf n}^-.
\]
Then it is easy to see that the subgroup
\begin{equation}\label{eqn:ZSL}
\index{Ca(G)@$C(G)$|main}%
C(G) 
:= \{ \gamma_{n\pi i} : n \in \Z \}  \simeq \Z
\end{equation}
coincides with the center of $G$.

Next, we give a parametrization of one-dimensional
representations of $K \simeq \R$ as
\begin{equation}\label{eqn:Khat}
\widehat{K} \simeq \C, 
\quad
\chi_\mu \leftrightarrow \mu
\end{equation}
by the formula
$
\chi_\mu (\gamma_z)
= e^{-\mu z}
$
or equivalently,
$
d\chi_\mu({\bf k}) = \mu
$.

We shall call $\chi_\mu$ simply as the $K$-type $\mu$.

Let $(\varpi,X)$ be a $(\mathfrak{g}_\C,K)$-module.
A non-zero vector $v \in X$ is a \textit{lowest weight vector} of weight $\mu \in
\C$ if $v$ satisfies
\begin{equation*}
 \varpi({\bf n}^-)v = 0,
\quad\text{and}\quad
\varpi({\bf k})v = \mu v.
\end{equation*}
We say $(\varpi,V)$ is a \textit{lowest weight module} of weight $\mu$ if $V$ is
generated by such $v$.
For each $\lambda \in \C$,
there exists   a unique  irreducible lowest weight
$(\mathfrak{g}_\C,K)$-module, to be denoted by $\pi_K(\lambda)$,
of weight $\lambda+1$.

With this normalization,
we pin down the following well-known properties of the $(\g_\C,K)$-module 
\index{*pi_K(\lambda)@$\pi_K(\lambda)$|main}%
$\pi_K(\lambda)$ for $\g = \s\l(2,\R)$:
\begin{fact}\label{fact:SL2rep}
~
\begin{enumerate}[\upshape 1)]
\item   
For a real $\lambda$ with $\lambda \ge -1$,
there exists   a  unique unitary
representation, denoted by 
\index{*pi(\lambda)@$\pi(\lambda)$|main}%
$\pi(\lambda)$, of $G = \widetilde{SL(2,\R)}$
such that its underlying $(\mathfrak{g}_\C,K)$-module is isomorphic to
$\pi_K(\lambda)$. 
\item   
$\pi(-1)$ is the trivial one-dimensional representation.
\item   
For $\lambda>0$, 
$\pi(\lambda)$ is a relative discrete series representation, namely,
its matrix coefficients are square integrable over $G$ modulo its
center $C(G)$.
\item   
$\pi (\frac{1}{2}) \oplus \pi (-\frac{1}{2})$ is the Weil
representation of $Mp(1,\R)$,
the two fold covering group of $SL(2,\R)$.
\item   
$\gamma_{ \pi i}$ 
(see \eqref{eqn:gammaz}) acts on $\pi(\lambda)$ as scalar
$e^{-\pi i(\lambda+1)}$.
\item   
$\pi(\lambda)$ is well-defined as a unitary representation of
$Mp(1,\R)$ if $\lambda \in \frac{1}{2}\Z$,
of $SL(2,\R)$ if $\lambda \in \Z$,
and of $PSL(2,\R)$ if $\lambda \in 2\Z+1$.
\item   
For $\lambda \ne -1,-3,-5,\dots,\pi_K(\lambda)$ is an infinite
dimensional representation.
For $\lambda>-1$, we fix a $G$-invariant inner product on the representation space of $\pi(\lambda)$. Then  we can find an orthonormal basis $\{v_\ell: \ell\in\N\}$ such that
\begin{alignat*}{2}
&\pi_K(\lambda)({\bf k})v_\ell
&&= (2\ell+\lambda+1)v_\ell,
\\
&\pi_K(\lambda)({\bf n}^+)v_\ell
&&= i\sqrt{(\ell+1)(\lambda+\ell+1)} \, v_{\ell+1},
\\
&\pi_K(\lambda)({\bf n}^-)v_\ell
&&= i\sqrt{(\lambda+\ell)\ell} \, v_{\ell-1}.
\end{alignat*}
Here, we set $v_{-1}=\{0\}$.
In particular,
$\pi_K(\lambda)$ has the $K$-type
$\{ \lambda+1, \lambda+3, \lambda+5, \dotsc \}$ with respect to the parametrization  \eqref{eqn:Khat}.
\item  
For $\lambda=-m$ $(m=1,2,\dotsc)$,
$\pi_K(\lambda)$ is an $m$-dimensional irreducible representation of 
$\s\l(2,\C)$.
\end{enumerate}
\end{fact}
By using Fact \ref{fact:SL2rep},
we can read from the formulas in 
Theorem \ref{thm:onbPhi} the following statement:
\begin{thm}\label{thm:Wka}
Suppose $a$ is a non-zero complex number and
$k$ is  a non-negative root  multiplicity function
 satisfying the inequality
\eqref{eqn:a2k}, i.e.\ 
$a+2\langle k\rangle+N-2 > 0$.
\begin{enumerate}[\upshape 1)]
\item 
$(\omega_{k,a}, W_{k,a}(\R^N))$ is a $\CC \times
(\mathfrak{g}_\C,K)$-module. 
\item 
As a $(\mathfrak{g}_\C,K)$-module,
$\omega_{k,a}$ is $K$-admissible and hence discretely decomposable
(see Definition \ref{def:DDR}).
\item 
$(\omega_{k,a}, W_{k,a}(\R^N))$ is decomposed into the direct sum of
 $\CC \times (\mathfrak{g}_\C, K)$-modules as follows:
\begin{equation}\label{eqn:Wkadeco}
W_{k,a}(\R^N) \simeq \bigoplus_{m=0}^\infty \mathcal{H}_k^m (\R^N)_{\big| S^{N-1}} \otimes
\pi_K ( \lambda_{k,a,m} ).
\end{equation}
\end{enumerate}
Here, 
 $\lambda_{k,a,m} = \frac{2m+2\langle k \rangle+N-2}{a}$
(see \eqref{3.17}). The Coxeter group $\CC$ acts on the first factor, and the Lie algebra $\s\l(2,\R)$ acts on the second factor of each summand in \eqref{eqn:Wkadeco}.
\end{thm}

\begin{proof}[Proof of Theorem \ref{thm:Wka}]
We fix a non-zero
$p \in \mathcal{H}_k^m(\R^N)$.
Then, 
it follows from Theorem \ref{thm:3.9} and Fact \ref{fact:SL2rep} 
that for 
$\mathfrak{sl}(2,\R)$ acts on the vector space
$$
\text{$\C$-span} \{ \Phi_\ell^{(a)}(p,\cdot): \ell \in \N \}
$$
as an irreducible lowest weight module $\pi_K(\lambda_{k,a,m})$.
By \eqref{eqn:L2RNsum},
we get the isomorphism \eqref{eqn:Wkadeco} as
$(\g_\C,K)$-modules.

On the other hand, the Coxeter group $\CC$ leaves
$\mathcal{H}_k^m(\R^N)$ invariant.
Furthermore, as we saw in
Lemma \ref{lem:Lomega}, the action of
$\CC$ and $\mathfrak{sl}(2,\R)$ commute with each other.
Hence, the first and third statements are proved.

It follows from the decomposition formula
\eqref{eqn:Wkadeco} that $\omega_{k,a}$ is $K$-admissible because the
$K$-type of an individual $\pi_K(\lambda)$ is of the form
$\{ \lambda+1, \lambda+3, \dotsc \}$ by Fact \ref{fact:SL2rep},
$\lambda_{k,a,m}$ increases as $m$ increases,  and
$\dim\mathcal{H}_k^m(\R^N)<\infty$.
Hence, the second statement is also proved.
\end{proof}

For $f,g\in L^2\big(\R^N, \vartheta_{k,a}(x)dx\big),$ we write its
inner product as 
\begin{equation}\label{3.23}
\index{1@\par\indexspace$\langle\kern-.2em\langle \ ,\ \rangle\kern-.2em\rangle _k$|main}%
\langle\!\langle f,g \rangle\!\rangle _k 
:=\int_{\R^N} f(x) \overline{g(x)} 
 \, \vartheta_{k,a}(x)dx.
\end{equation} 

\begin{prop}\label{prop:iu}
Suppose that $a>0$ and that $k$ is a non-negative multiplicity function such that $a+2\langle k\rangle+N-2 > 0$.
Then, the representation $\omega_{k,a}$ of\/ $\mathfrak{sl}(2,\R)$ on
$W_{k,a}(\R^N)$ is infinitesimally unitary with respect to 
the inner product
$\langle\!\langle \ , \ \rangle\!\rangle_k$, namely,
$$
\langle\!\langle \omega_{k,a}(X)f,g\rangle\!\rangle_k
= -\langle\!\langle f,\omega_{k,a}(X)g\rangle\!\rangle_k
$$
for any $X\in\mathfrak{sl}(2,\R)$ and
$f,g\in W_{k,a}(\R^N)$.
\end{prop}

\begin{proof}
As we saw in \eqref{eqn:2.7} that the Dunkl operators are
skew-symmetric with respect to the measure $\vartheta_k(x)dx$.
In view of the definitions of $\Delta_k$
(see \eqref{eqn:2.2}) and
$\EE_{k,a}^- = \frac{i}{a} \Vert x\Vert^{2-a} \Delta_k$
(see \eqref{eqn:3.3}),
we see that $\EE_{k,a}^-$ is a skew-symmetric operator with respect to
the inner product $\langle\!\langle \cdot,\cdot\rangle\!\rangle_k$.
Likewise for $\EE_{k,a}^+$.
Further, the commutation relation   $\HH_{k,a}=[\EE_{k,a}^+, \EE_{k,a}^-]$ shows  that $\HH_{k,a}$ is also skew-symmetric. 
Thus, for all $X\in  \s\l(2,\R),$  $\omega_{k,a}(X)$ is skew-symmetric. 
\end{proof}

\subsection{The integrability of  the representation
$\omega_{k,a}$}\label{subsec:3.2}
\hfill\\
Applying the general result on discretely decomposable representations
(see Fact \ref{fact:deco}) to our specific setting where $G$ is the
universal covering group of $SL(2,\R)$,
we get the following two theorems:
\begin{thm}\label{thm:3.12}
Suppose $a>0$ and $k$ is a  non-negative multiplicity function 
satisfying
\begin{equation}\label{eqn:unitary}
a+2\langle k\rangle+N-2 > 0.
\end{equation}
Then the infinitesimal representation $\omega_{k,a}$ of $\s\l(2,\R)$
lifts to a unique unitary representation,
to be denoted by 
\index{*ZOmegaka@$\Omega_{k,a}$|main}%
$\Omega_{k,a} $, 
of $G$ on the
Hilbert space $L^2\big(\R^N, \vartheta_{k,a}(x)dx\big) $.
In particular, we have
$$\omega_{k,a}(X)={{\rm d} \over {{\rm d} t}}\Big|_{t=0} \Omega_{k,a}(\Exp(tX)),\qquad X\in \g,$$
on $W_{k,a}(\R^N)$,
the dense subspace \eqref{eqn:Wka} of $L^2(\R^N,\vartheta_{k,a}(x)dx)$. 
Here, we have written $\Exp$ for the exponential map of the Lie
algebra $\s\l(2,\R)$ into $G.$ 
\end{thm} 

\begin{thm}\label{thm:3.16} 
Retain the assumption of Theorem \ref{thm:3.12}.
Then, as a representation of the direct product group
 $\CC\times G$, the
unitary representation $L^2\big(\R^N, \vartheta_{k,a}(x)dx\big)$ 
 decomposes
discretely as 
\begin{equation}\label{eqn:Hmkpi}
L^2\big(\R^N, \vartheta_{k,a}(x)dx\big)
= \sideset{}{^\oplus}\sum_{m=0}^\infty \Bigl(\Hkm(\R^N)|_{S^{N-1}}\Bigr)\otimes \pi(\lambda_{k,a,m}).
\end{equation}
Here, we recall that $\CC$ is the Coxeter group of the root system,
$G$ is the universal covering group of $SL(2,\R)$,
and
$\lambda_{k,a,m}=\frac{2m+2\langle k\rangle+N-2}{a}$ (see \eqref{3.17}).
The decomposition \eqref{eqn:Hmkpi} of the Hilbert space
$L^2(\R^N,\vartheta_{k,a}(x)dx)$
is given by the formula \eqref{eqn:L2RNsum} and $\pi(\lambda)$ $(\lambda>-1)$ is
the irreducible unitary representation of $\widetilde{SL(2,\R)}$
described in Fact \ref{fact:SL2rep}.
In particular,
the summands are mutually orthogonal with respect to the inner product \eqref{3.23} on $L^2\big(\R^N, \vartheta_{k,a}(x)dx\big).$
\end{thm}

\begin{rem}[The $N=1$ case]\label{rem:3.6}
 In \cite{Kos} Kostant exhibits a family of representations with
 continuous parameter of $\s\l(2,\R)$ by  second order differential operators on $(0,\infty)$. He uses Nelson's result\/ \cite{Nels} to study the exponentiation of
 such representations. See also \cite{RR}. 
 
In the $N=1$ case, the decomposition  in Theorem \ref{thm:Wka} (and hence in Theorem \ref{thm:3.16}) is reduced to a finite sum 
because $\Hkm(\R^N) = 0$ if $m \ge 2$ and $N=1$. Indeed, there are two summands according to even ($m=0$) and odd ($m=1$) functions. We observe that the difference operator  in \eqref{eqn:2.3} vanishes on even functions of one variable, so the Dunkl Laplacian $\Delta_k$ collapses to the differential operator ${{d^2}\over {dx^2}} +\frac{2k}{x}{d\over {dx}}.$ Thus, our generators \eqref{eqn:3.3} acting on even functions on $ (0,\infty)$ take the form 
$$\HH_{k,a}={2\over a} x{d\over {dx}} +{{2k+a-1}\over a},\qquad
 \EE_{k,a}^+= {i\over a} x^a,\qquad \EE_{k,a}^-={i\over a} x^{2-a}\Big( {{d^2}\over {dx^2}} +\frac{2k}{x}{d\over {dx}}\Big).$$

We may compare these with the generators in Kostant's paper \cite{Kos}, where his generators on $(0,\infty)$ are
$$iy,\qquad  2y{d\over {dy}}+1,\qquad  i\Big( y{{d^2}\over {dy^2}}+{d\over {dy}}-{{r^2}\over {4y}}\Big),
 $$ which by the substitution $y={1\over a} x^a$ and $\varphi(y)=x^{k-{1\over 2}} \varphi(x)$  become our operators $\{ \HH_{k,a}, \EE_{k,a}^+, \EE_{k,a}^-\}$ with $r={{2k-1}\over a}.$ Note that our generators acting on odd functions do not appear in Kostant's picture. 
 \end{rem}

\begin{rem}\label{rem:a2kN}
The assumption $a+2\langle k\rangle+N-2>0$ implies
\[
\lambda_{k,a,m} > -1
\quad\text{for any $m\in\N$},
\]
whence there exists an irreducible, infinite dimensional unitary representation
$\pi(\lambda_{k,a,m})$ of $G$ such that its underlying
$(\g_\C,K)$-module is isomorphic to $\pi_K(\lambda_{k,a,m})$ by Fact
\ref{fact:SL2rep} (1).
\end{rem}
By the explicit construction of the direct summand in Theorem
\ref{thm:3.9},
we have
\begin{coro}\label{cor:minK}
As a representation of $\widetilde{SL(2,\R)}$, 
minimal $K$-types of the irreducible summands in \eqref{eqn:Hmkpi} are
given by
\[
h(x) \exp \Bigl( -\frac{1}{a}\Vert x\Vert^a\Bigr),
\quad h \in \Hkm(\R^N).
\]
\end{coro}

As we have seen that $\omega_{k,a}$ lifts to the unitary representation
$\Omega_{k,a}$ of the universal covering group
$G = \widetilde{SL(2,\R)}$ for any $k$ and $a$ with certain positivity
\eqref{eqn:unitary}.
On the other hand,
if $k$ and $a$ satisfies a certain rational condition (see below),
then $\Omega_{k,a}$ is well-defined for some finite covering groups of
$PSL(2,\R)$.
This representation theoretic observation gives an explicit formula of the
order of the $(k,a)$-generalized Fourier transform $\Fka$ (see Section
\ref{sec:5}). We pin down a precise statement here.
\begin{prop}\label{prop:3.13}   
Retain the notation of Theorem \ref{thm:3.12}, and recall the
definition of the index $\langle k\rangle$ from \eqref{eqn:2.6}. 
Then the unitary representation $\Omega_{k,a}$ of the universal
covering group $G$ of $SL(2,\R)$ is well-defined also as a representation of some finite
covering group of $PSL(2,\R)$ if and only if both $a$ and $\langle k\rangle$
are rational numbers.
\end{prop}

\begin{proof}
It follows from Fact \ref{fact:SL2rep} 5) that the central element 
$\gamma_{n\pi i} \in C(G)$ acts on $\pi(\lambda_{k,a,m})$ by the
scalar
\[
e^{-\pi in(\lambda_{k,a,m}+1)}
= \exp \Bigl( -\frac{n}{a}\cdot 2\pi m i\Bigr)
  \exp \Bigl( -\frac{N+2\langle k\rangle+a-2}{a} n \pi i \Bigr).
\]
This equals $1$ for all $m$ if and only if
\[
\frac{n}{a} \in \Z
\quad\text{and}\quad
\frac{n(2\langle k\rangle+N-2+a)}{a} \in 2\Z.
\]
It is easy to see that there exists a non-zero integer $n$ satisfying these two conditions if
and only if both $a$ and $\langle k\rangle$ are rational numbers.
For such $n$,
$\Omega_{k,a}$ is well-defined for $G/n\Z$.
Hence, Proposition \ref{prop:3.13} is proved.
\end{proof}

We recall from \eqref{eqn:ZSL} that we have identified the center
$C(G)$ of the simply-connected Lie group $G = \widetilde{SL(2,\R)}$
with the integer group $\Z$.
Then, we have
$$
PSL(2,\R) \simeq G/\Z,
\quad
SL(2,\R) \simeq G/2\Z,
\quad
Mp(1,\R) \simeq G/4\Z.
$$
As a special case of Proposition \ref{prop:3.13} and its proof,
we have:
\begin{rem}\label{rem:3.14}
Let $\Omega_{k,a}$ be the unitary representation of the universal
covering group $G$.
\begin{itemize}
\item[(1)] Suppose $a=2$. 
\begin{itemize}
\item[(a)] $\Omega_{k,2}$ descends to $SL(2,\R)$ if and only if $2\langle k\rangle +N$ is an even integer. 
\item[(b)] $\Omega_{k,2}$ descends to $Mp(1,\R)$ if and only if $2\langle k\rangle +N$ is an integer. 
\end{itemize}
This compares well with the Schr\"odinger model on  $L^2(\R^N)$ of the
Weil representation $ \Omega_{0,2}$ of the metaplectic group
$Mp(N,\R)$ and its restriction to a subgroup locally isomorphic to $SL(2,\R)$
 (cf. \cite{Weil} and \cite{HT}). 
\end{itemize}
\begin{itemize}
\item[(2)] Suppose $a=1$. 
\begin{itemize}
\item[(a)]
$\Omega_{k,1}$ descends to $PSL(2,\R)$ if and only if $2\langle
k\rangle+N$ is an odd integer.
\item[(b)] $\Omega_{k,1}$ descends to $SL(2,\R)$ if and only if $2\langle k\rangle  $ is an integer. 
\item[(c)] $\Omega_{k,1}$ descends to $Mp(1,\R)$ if and only if $4\langle k\rangle  $ is an integer. 
\end{itemize}
The case $k\equiv0$ corresponds to the Schr\"odinger model on  $L^2(\R^N, {{dx}\over {\Vert x\Vert}})$ of  the  minimal  representation 
 $\Omega_{0,1}$ of the conformal group and its restriction to a
 subgroup locally isomorphic to $SL(2,\R)$  (cf. \cite{KO}). 
\end{itemize}
 \end{rem}
\begin{rem}\label{rem:3.15}
C. Dunkl reminded us of that the parity
 condition of $2\langle k\rangle +N$ appeared also in a different
 context, i.e., in \cite[Lemma 5.1]{DJO}, where the authors investigated 
a sufficient condition on $k$ for the existence and uniqueness of expanding a homogenous polynomial in terms of $k$-harmonics.
\end{rem}

\subsection{Connection with the Gelfand--Gindikin program}\label{sec:3.7}
\hfill\\
We consider the following closed cone in
$\mathfrak{g}=\mathfrak{sl}(2,\R)$ defined by
\[
W := \left\{ \begin{pmatrix}a&b\\c&-a\end{pmatrix}:
     a^2+bc\le0, \  b\ge c \right\}.
\]
Then, $W$ is $SL(2,\R)$-invariant
 and is
expressed as
\begin{equation*}
W = i \operatorname{Ad}(SL(2,\R)) \R_{\ge0} {\bf k}.
\end{equation*}

We write
$\exp_\C: \g_\C \to SL(2,\C)$
for the exponential map. 
Its restriction to $iW$ is an injective map,
and we define the following subset $\Gamma(W)$ of $SL(2,\C)$ by
$$
\index{*CGamma(W)@$\Gamma(W)$|main}%
\Gamma(W)
:=SL(2,\R)\exp_{\C}(iW).
$$ 
Since $W$ is $SL(2,\R)$-invariant, 
$\Gamma(W)$ becomes a semigroup, sometimes referred to as the \textit{Olshanski
semigroup}.

Denote by $\widetilde{\Gamma(W)}$ the universal covering semigroup of $\Gamma(W),$ and write  
$$
\Exp:\g+i  W\rightarrow \widetilde {\Gamma( W)}
$$ 
for the lifting of $\exp_\C{_{|\g+i  W}}: \g+i  W\rightarrow   \Gamma( W).$ Then $\widetilde {\Gamma(  W)}=  \widetilde {SL(2,\R)} \Exp(iW) $ and the polar map
$$  
\widetilde {SL(2,\R)}\times  W\rightarrow \widetilde {\Gamma( W)}, \;(g,X)\mapsto g\Exp(iX)
$$
is a homeomorphism. 

Since $W$ is an $\operatorname{Ad}(SL(2,\R))$-invariant cone,
$\Gamma(W)$
 is invariant under the action of $SL(2,\R)$ from the left and
right.
Thus, the semigroup $\Gamma(W)$ is written also as
\[
\Gamma(W) = SL(2,\R) \exp_{\C} (- \R_{\ge0} {\bf k}) SL(2,\R).
\]
Its interior is given by
\[
\Gamma(W^0) = SL(2,\R) \exp(- \R_+{\bf k}) SL(2,\R).
\]
See \cite[Theorem  7.25]{HN}.
Accordingly, we have
\[
\widetilde{\Gamma(W)} = \widetilde{SL(2,\R)}
\exp_{\C} (-\R_{\ge0} {\bf k}) \widetilde{SL(2,\R)}.
\]

By Theorem \ref{thm:3.9}, $\Omega_{k,a}$ 
 is a discretely decomposable unitary representation of $\widetilde{SL(2,\R)}$ on
$L^2(\R^N,\vartheta_{k,a}(x)dx)$.
It has a lowest weight $(2\langle k\rangle+N+a-2)/a.$ 
It then follows  from \cite[Theorem  B]{HN2} that $\Omega_{k,a}$
 extends to a representation of the Olshanski semigroup
$\widetilde{\Gamma(W)},$ denoted by the same symbol $\Omega_{k,a}$, such that:
\begin{itemize}
\item[(P1)] $\Omega_{k,a}:\widetilde{\Gamma(W)}\rightarrow \mathcal B(L^2)$ is strongly continuous semigroup homomorphism.
\item[(P2)] For all $ f\in L^2(\R^N, \vartheta_{k,a}(x)dx),$ the map
$\gamma\mapsto \langle\!\langle \Omega_{k,a}(\gamma)f,
f\rangle\!\rangle_k$ is holomorphic in the interior of  $\widetilde{\Gamma(W)}$. 
\item[(P3)] $\Omega_{k,a}(\gamma)^*=\Omega_{k,a}(\gamma^\sharp),$ where $\gamma^\sharp=\Exp(iX) g^{-1}$ for $\gamma=g\Exp(iX).$
\end{itemize}

Here, we have denoted by $\mathcal{B}(L^2)$ the space of bounded
operators on $L^2(\R^N,\vartheta_{k,a}(x)dx)$.

\begin{rem}
The Gelfand--Gindikin program \cite{GG} seeks for the understanding of
a `family of irreducible representations' by using complex geometric
methods.
This program has been particularly developed for lowest weight
representations by
Olshanski \cite{Ols} and Stanton \cite{Sta}, 
Hilgert, Neeb \cite{HN}, and some others.
The study of our holomorphic semigroup $\Omega_{k,a}$ by using the Olshanski semigroup
$\widetilde{\Gamma(W)}$ may be regarded as a descendant of this program.
\end{rem}

Henceforth we will use the notation 
$
\index{Cb^+@$\C^+$|main}%
\C^+ 
:=\{z\in \C\;|\; \Re(z)\geq 0\}$ and 
$
\index{Cb^{++}@$\C^{++}$|main}%
\C^{++} 
:=\{z\in \C\;|\; \Re(z)> 0\} .$

For
$z\in \C^+$, we extend the one-parameter subgroup $\gamma_z$ $(z\in i\R)$
(see \eqref{eqn:gammaz}) holomorphically as
\begin{equation}\label{3.27} 
\index{*cgamma_z@$\gamma_z$|main}%
\gamma_z 
:= \Exp(-z {\bf k})
= \Exp(iz\begin{pmatrix}0&1\\-1&0\end{pmatrix})\in \widetilde {\Gamma (  W)}.
\end{equation} 
Then, the operators $\Omega_{k,a}(\gamma_z) $ have the following property:
\begin{align*}
&\Omega_{k,a}(\gamma_{z_1})  \,  \Omega_{k,a}(\gamma_{z_2})
=\Omega_{k,a}(\gamma_{z_1+z_2}),\qquad \forall z_1,z_2\in \C^+ ,
\\
&\Omega_{k,a}(\gamma_z)^*
= \Omega(\gamma_{\bar{z}}),
\qquad z \in \C^+,
\\
&\Omega_{k,a}(\gamma_0) =\id.
\end{align*}
Following the formulation of \cite[Proposition 3.6.1]{KM} (the case
$k\equiv0$, $a=1$),
we summarize basic properties of the holomorphic representation 
$\Omega_{k,a}$ of the semigroup $\widetilde{\Gamma(  W)}.$
\begin{thm}\label{thm:3.18} 
Suppose $a>0$ and $k$ is a non-negative 
multiplicity function on the root system
satisfying \eqref{eqn:a2k}, i.e.\ $a+2\langle k\rangle+N-2 > 0$.
\begin{itemize}
\item[(1)] The map $$\widetilde {\Gamma( W)}\times L^2\big(\R^N, \vartheta_{k,a}(x)dx\big)\longrightarrow L^2\big(\R^N, \vartheta_{k,a}(x)dx\big),\quad (\gamma,f)\mapsto \Omega_{k,a}(\gamma)f$$ is continuous.
\item[(2)] 
For any 
$p \in \mathcal{H}_k^{(m)}(\R^N)$ and $\ell\in \N$,
 $\Phi_\ell^{(a)}(p,\cdot)$ (see \eqref{3.16}) is an
eigenfunction of the operator $\Omega_{k,a}(\gamma_z)
 = \exp(\omega_{k,a}(-z{\bf k}))$:
$$ \Omega_{k,a}(\gamma_z)\Phi_{\ell }^{(a)}
(p,x)=e^{-z(\lambda_{k,a,m}+1+2\ell)} \Phi_{\ell }^{(a)}(p,x) ,$$
where  $\lambda_{k,a,m}=\frac{1}{a}(2m+2\langle k\rangle +N-2)$ (see \eqref{3.17}).
\item[(3)]  The operator norm $\Vert \Omega_{k,a}(\gamma_z)\Vert_{\rm
op}$  is 
$\exp(-\frac{1}{a}(2\langle k\rangle+N+a-2)\Re z)$.
\item[(4)]  If\/ $\operatorname{Re}(z)>0$,  then $\Omega_{k,a}(\gamma_z)$ is a Hilbert--Schmidt operator.
\item[(5)]   If\/ $\operatorname{Re}(z)=0,$ then $\Omega_{k,a}(\gamma_z)$ is a unitary operator.
\item[(6)]  The representation $\Omega_{k,a}$ is faithful on
$\widetilde{\Gamma(W)}$ if at least one of $a$ or $\langle k\rangle$
is irrational,
  and on $ \widetilde{\Gamma(W)}/ D$ for some discrete abelian kernel
$D$ if both $a$ and $\langle k\rangle$ are rational.
\end{itemize}
\end{thm}
\begin{proof}
The second statement follows from \eqref{3.21 a}.
The fifth statement is a special case of Theorem \ref{thm:3.12}.
The proof of the other statements is parallel to that of
\cite[Proposition 3.6.1]{KM}, and we omit it.
\end{proof}

\section{The integral representation of the holomorphic semigroup $\Omega_{k,a}(\gamma_z)$}\label{sec:4}

We have seen in Theorem \ref{thm:3.18} that $\Omega_{k,a}(\gamma_z)$
is a Hilbert--Schmidt operator for $\Re z >0$ and is a unitary operator
for $\Re z=0$.
By the Schwartz kernel theorem,
the operator $\Omega_{k,a}(\gamma_z)$ can be expressed by means of a
distribution kernel $\Lambda_{k,a}(x,y;z)$.
If we adopt Gelfand's notation on
 a generalized functions,
 we may write the operator $\Omega_{k,a}(\gamma_z)$ on
$L^2(\R^N,\vartheta_{k,a}(x)dx)$ as an `integral transform'
against the measure $\vartheta_{k,a}(x)dx$:
\begin{equation}\label{eqn:Lambda}
\Omega_{k,a}(\gamma_z)f(x)
= c_{k,a} \int_{\R^N} \Lambda_{k,a}(x,y;z)
f(y) \vartheta_{k,a}(y)dy.
\end{equation}
Here, we have normalized  the kernel $\Lambda_{k,a}(x,y;z)$ by the constant $c_{k,a}$ that will be defined in \eqref{eqn:cka}.
In light of the unitary isomorphism
\[
L^2(\R^N,\vartheta_{k,a}(x)dx) \stackrel{\sim}{\to} L^2(\R^N,dx), \ 
f(x) \mapsto f(x)\vartheta_{k,a}(x)^{\frac{1}{2}}.
\]
We see that
$\Lambda_{k,a}(x,y;z)\vartheta_{k,a}(x)^{\frac{1}{2}}\vartheta_{k,a}(y)^{\frac{1}{2}}$
is a tempered distribution of $(x,y)\in\R^N \times \R^N$.

The goal of this section is to find the kernel $\Lambda_{k,a}(x,y;z)$.
The main result of this section is Theorem \ref{thm:4.9}.

\subsection{Integral representation for the radial part of $\Omega_{k,a}(\gamma_z)$}\label{sec:4.1}\hfill\\
By Lemma \ref{prop:slreduction},
the $\mathfrak{sl}_2$-action $\omega_{k,a}$ on 
$C^\infty(\R^N\setminus\{0\})$ 
(see \eqref{eqn:omegaka} for definition) can be described in
a simple form on each $k$-spherical component $\Hkm(\R^N)$, namely,
it can be expressed as the action only
on the radial direction.
Accordingly, we can define the `radial part' of the 
holomorphic semigroup $\Omega_{k,a}^{(m)}(\gamma_z) $ (see
\eqref{eqn:4.3} below for definition) on $L^2(\R_+, r^{2\langle k\rangle+N+a-3}dr)$.
The main result of this subsection is the integral formula for
$\Omega_{k,a}^{(m)}(\gamma_z)$, 
which will be given in Theorems \ref{thm:4.2} and \ref{thm:4.4}.

\subsubsection{Radial part of holomorphic semigroup}
\hfill\break
Recall that $\Hkm(\R^N)$ is the space of $k$-harmonic polynomials of
degree $m\in\N$.
Let
$$\alpha_{k,{a}}^{(m)}:
\Hkm(\R^N)|_{S^{N-1}}\otimes L^2(\R_+, r^{2\langle k\rangle+N+a-3}dr) \longrightarrow L^2\big(\R^N, \vartheta_{k,a}(x)dx\big)$$  
be a linear map defined by  
$$\alpha_{k,{a}}^{(m)}( p \otimes f)(x)=p\big({x\over {\Vert x\Vert}}\big)f(\Vert x\Vert)
\quad\text{for $p \in \mathcal{H}_k^m(\R^N)|_{S^{N-1}}$ and 
$f\in L^2(\R_+, r^{2\langle k\rangle+N+a-3}dr)$} .
$$  
Summing up $\alpha_{k,a}^{(m)}$, we get
a direct sum decomposition of the Hilbert space: 
\begin{equation}\label{eqn:4.1}
L^2\big(\R^N, \vartheta_{k,a}(x)dx\big)=\sideset{}{^\oplus}\sum_{{m}\in \N} 
 \Hkm(\R^N)|_{S^{N-1}}\otimes L^2\big(\R_+, r^{2\langle k\rangle+N+a-3}dr\big). 
\end{equation}
It follows from Theorem \ref{thm:3.16} that 
the unitary representation $\Omega_{k,a}$ of $\widetilde{SL(2,\R)}$ on
the Hilbert space
$L^2(\R^N,\vartheta_{k,a}(x)dx)$ 
induces a family of unitary operators,
to be denoted by $\Omega_{k,a}^{(m)}(\gamma_z)$ $(z\in i\mathbb{R}, m\in\N)$, 
 on $L^2(\R_+, r^{2\langle k\rangle+N+a-3}dr) $ such that
\begin{equation}\label{eqn:4.2}
\alpha_{k,{a}}^{(m)}\Big(p \otimes \Omega_{k, a}^{(m)}(\gamma_z)(f)\Big)=\Omega_{k,a}(\gamma_z)\Big(\alpha_{k,{a}}^{(m)}(p\otimes f)\Big).\end{equation}

As is Theorem \ref{thm:3.18} for $\Omega_{k,a}(\gamma_z)$,
the unitary operator $\Omega_{k, a}^{(m)}(\gamma_z)$ extends to a
holomorphic semigroup of Hilbert--Schmidt operators on
$L^2(\R_+,r^{2\langle k\rangle+N+a-3} dr)$
for $\Re(z)>0 .$  
Further, 
there exists a unique kernel $\Lambda_{k,a}^{(m)}(r,s;z)$ for each
$z$ and $m\in\N$ such that
\begin{equation}\label{eqn:4.3}
 \Omega_{k, a}^{(m)}(\gamma_z) f(r) =\int_0^\infty f(s)
\index{*Lambda_{k,a}^{(m)}(r,s;z)@$\Lambda_{k,a}^{(m)}(r,s;z)$|main}%
\Lambda_{k,a}^{(m)}(r,s;z) 
s^{2\langle k\rangle+N+a-3}ds,\end{equation}
holds for any $f\in L^2(\R_+,r^{2\langle k\rangle+N+a-3}dr)$.

According to the direct sum \eqref{eqn:4.1}, the semigroup
$\Omega_{k,a}(\gamma_z)$ is decomposed as follows:    
\begin{equation}\label{eqn:Omdeco}
\Omega_{k,a}(\gamma_z)=\sideset{}{^\oplus}\sum_{{m}\in \N}
{\rm id}_{ |{ \cal H_k^m(\R^N)}} \otimes \Omega_{k, a}^{(m)}(\gamma_z) .
\end{equation}

Comparing the integral expressions \eqref{eqn:Lambda} and
\eqref{eqn:4.3} of $\Omega_{k,a}(\gamma_z)$ and
$\Omega_{k,a}^{(m)}(\gamma_z)$ respectively,
we see that the kernels $\Lambda_{k,a}(x,y;z)$ and $\Lambda_{k,a}^{(m)}(r,s;z)$
satisfy the following identities:
\begin{equation}\label{eqn:Lam1}
c_{k,a} \int_{\R^N} \Lambda_{k,a}(x,y;z) 
p\Bigl(\frac{y}{\Vert y\Vert}\Bigr) f(\Vert y\Vert)
\vartheta_{k,a}(y) dy
= p \Bigl(\frac{x}{\Vert x\Vert}\Bigr)
  \int_0^\infty f(s) \Lambda_{k,a}^{(m)} (r,s;z)
  s^{2\langle k\rangle+N+a-3} ds
\end{equation}
for any $p \in \mathcal{H}_k^{(m)} (\R^N)$ and 
$f \in L^2(\R_+,r^{2\langle k\rangle+N+a-3}dr)$.

In light of the following formula for the measures
\[
\vartheta_{k,a}(y)dy
= \vartheta_k(\eta)s^{2\langle k\rangle+N+a-3} d\sigma(\eta)ds
\]
with respect to polar coordinates $y=s\eta$,
we see that \eqref{eqn:Lam1} is equivalent to
\begin{equation}\label{eqn:Lam2}
c_{k,a} \int_{S^{N-1}} \Lambda_{k,a} (r\omega,s\eta;z)p(\eta)
\vartheta_k(\eta)d\sigma(\eta)
= p(\omega) \Lambda_{k,a}^{(m)} (r,s;z).
\end{equation}

Therefore, 
the distribution $\Lambda_{k,a}$ is determined by the set of functions
$\Lambda_{k,a}^{(m)}$ $(m\in\N)$ as follows:
\begin{prop}\label{pro:Lm}
Fix $z\in\C$ with $\operatorname{Re}z \ge 0$.
Then, the distribution $\Lambda_{k,a}(x,y;z)$ on $\R^N\times \R^N$ is
characterized by the condition \eqref{eqn:Lam2} for any
$p\in\mathcal{H}_k^m(\R^N)$ and any $m\in\N$.
\end{prop}
The relation between $\Lambda_{k,a}$ and $\Lambda_{k,a}^{(m)}$
$(m\in\N)$ will be discussed again in Theorem \ref{thm:4.14} by means
of the `Poisson kernel'.

\subsubsection{The case $\Re(z)>0$}\hfill\\
Suppose $\Re(z)>0$.
Then, $\Omega_{k, a}^{(m)} (\gamma_z)$  is a
Hilbert--Schmidt operator on $L^2(\R_+, r^{2\langle
k\rangle+N+a-3}dr),$ and consequently, the kernel $
\Lambda_{k,a}^{(m)}(\cdot,\cdot\;;z)$ is square integrable function
with respect to the measure $(rs)^{2\langle k\rangle+N+a-3}dr \, ds$.

We shall find a closed formula for $\Lambda_{k,a}^{(m)}(r,s;z)$.
Let us fix $m\in\N$ (as well as $k$ and $a$) once and for all. 
We have given in Proposition \ref{prop:falr} 
 an explicit orthonormal basis
$\{f_{\ell,m}^{(a)}(r): \ell\in\N\}$ of 
$L^2(\R_+,r^{2\langle k\rangle+N+a-3}dr)$.
On the other hand, it follows from Theorem \ref{thm:3.18} (2) that
$\widetilde{\Phi}_\ell^{(a)}(p,x) = f_{\ell,m}^{(a)}(r)p(\omega)$
(see \eqref{eqn:nPhi}) is an eigenfunction of the Hilbert--Schmidt
operator $\Omega_{k,a}(\gamma_z)$:
\[
\Omega_{k,a}(\gamma_z) \widetilde{\Phi}_\ell^{(a)} (p,x)
= e^{-z(2\ell+\lambda_{k,a,m}+1)} \widetilde{\Phi}_\ell^{(a)}(p,x).
\]

Using the identity \eqref{eqn:4.2},  we deduce that 
\begin{equation}\label{eqn:4.5}
 \Omega_{k,a}^{(m)}(\gamma_z) f_{\ell,m}^{(a)}(r)=e^{-z(2\ell+\lambda_{k,a,m}+1)} f_{\ell,m}^{(a)}(r), \end{equation}  where the constant 
$\lambda_{k,a,m}$ is defined in  \eqref{3.17}.  Hence, the kernel
$\Lambda_{k,a}^{(m)}(r,s;z)$ in \eqref{eqn:4.3} is given by the
following series expansion:
\begin{equation*}
\Lambda_{k,a}^{(m)}(r,s;z)=\sum_{\ell=0}^\infty f_{\ell,m}^{(a)}(r) f_{\ell,m}^{(a)}(s) e^{-z(\lambda_{k,a,m}+1+2\ell)}.
\end{equation*}
In view of the definition \eqref{eqn:4.4} of $f_{\ell,m}^{(a)}(r)$,
$\Lambda_{k,a}^{(m)}(r,s;z)$ amounts to
\begin{equation*}
{{e^{-z(\lambda_{k,a,m}+1) }(rs)^{m} e^{-{1\over a}(r^a+s^a)}}\over {a^{\lambda_{k,a,m}} 2^{-(\lambda_{k,a,m}+1)}}}\sum_{\ell=0}^\infty {{\Gamma(\ell+1)}\over {\Gamma(\lambda_{k,a,m}+\ell+1)}} e^{-2\ell z} 
L_\ell ^{(\lambda_{k,a,m})}\Big({2\over a} r^a\Big) L_\ell ^{(\lambda_{k,a,m})}\Big({2\over a} s^a\Big). 
\end{equation*}

In order to compute this series expansion,
we recall some basic identities of Bessel functions.
Let $I_\lambda$ be the $I$-Bessel function defined by
\[
\index{Ia_\lambda(w)@$I_\lambda(w)$|main}%
I_\lambda(w)
:= e^{-\frac{\pi}{2}\lambda i}
   J_\lambda(e^{\frac{\pi}{2}i} w).
\]

It is also convenient to introduce the normalized $I$-Bessel function by
\begin{align}\label{eqn:4.15}
\index{Iat_\lambda(w)@$\protect\widetilde I_\lambda(w)$|main}%
\widetilde I_\lambda(w) 
:={}&\Big({w\over 2}\Big)^{-\lambda}I_\lambda(w)=
\sum_{\ell =0}^\infty {{w^{2\ell}}\over {2^{2\ell} \ell ! \Gamma(\lambda+\ell+1)}}
\\
={}&\frac{1}{\sqrt{\pi} \, \Gamma(\lambda+\frac{1}{2})}
    \int_{-1}^1 e^{wt} (1-t^2)^{\lambda-\frac{1}{2}} dt.
\label{eqn:Itildeint}
\end{align} 
We note that $\widetilde{I}_\lambda(w)$ is an entire function of $w\in\C$
satisfying
\begin{equation*}
\widetilde{I}_\lambda(0)
= \frac{1}{\sqrt{\pi}\, \Gamma(\lambda+\frac{1}{2})}.
\end{equation*}
Now, we can use
 the following Hille--Hardy identity \cite[(6.2.25)]{AAR}
\begin{align*}
\sum_{\kappa=0}^\infty  {{\Gamma(\kappa+1)}\over {\Gamma(\lambda+\kappa+1)}} 
L_\kappa ^{(\lambda)}(u) L_\kappa ^{(\lambda)}(v)w^\kappa
&={1\over
{1-w}} \exp\Big(-{{(u+v)w}\over {1-w}}\Big)
(uvw)^{-\frac{\lambda}{2}}I_\lambda\Big({{2\sqrt {uv w}}\over
{1-w}}\Big)
\\
&= \frac{1}{(1-w)^{\lambda+1}}
   \exp \left( -\frac{(u+v)w}{1-w} \right)
   \widetilde{I}_\lambda
   \left( \frac{2\sqrt{uvw}}{1-w} \right).
\end{align*}
Here the left-hand side converges for $|w|<1$.
Hence, we get a closed formula for
$\Lambda_{k,a}^{(m)}(r,s;z)$:
\begin{align}\label{eqn:4.6}
\Lambda_{k, a}^{(m)}(r,s;z)
&={{(rs)^{-\langle k\rangle-\frac{N}{2}+1}\over {\sinh(z)}}}e^{-{1\over a}(r^a+s^a)\coth(z)}I_{\lambda_{k,a,m}}\left(
{2\over a}{ {(rs)^{\frac{a}{2}}} \over {\sinh(z)}}\right)
\\
&= \frac{(rs)^m}{a^{\lambda_{k,a,m}}(\sinh(z))^{\lambda_{k,a,m}+1}}
   e^{-\frac{1}{a}(r^a+s^a)\coth(z)}
   \widetilde{I}_{\lambda_{k,a,m}} 
   \left( \frac{2}{a} \frac{(rs)^{\frac{a}{2}}}{\sinh(z)} \right).
\nonumber
\end{align}

Next, let us give an upper estimate of the kernel function
$\Lambda_{k,a}^{(m)}(r,s;z)$.
For this,
we recall from \cite[\S4.2]{KM} the following elementary lemma.
\begin{lem}\label{lem:abz}
For $z = x+ iy$, we set
\begin{align*}
&\alpha(z) := \frac{\sinh(2x)}{\cosh(2x)-\cos(2y)},
\\
&\beta(z) := \frac{\cos(y)}{\cosh(x)}.
\end{align*}
Then, we have
\begin{enumerate}[\upshape 1)]
\item  
\begin{alignat}{2}
&\Re \coth(z) &&= \alpha(z),
\label{eqn:az}
\\
&\Re \frac{1}{\sinh(z)} &&= \alpha(z)\beta(z).
\label{eqn:abz}
\end{alignat}
\item  
If $z \in \C^+ \setminus i\pi\Z$,
then we have $\cosh(2x)-\cos(2y)>0$,
and
\[
\alpha(z) \ge 0 \quad\text{and}\quad |\beta(z)| < 1.
\]
\item  
If\/ $\Re z > 0$,
then
$
\alpha(z) > 0
$.
\end{enumerate}
\end{lem}
We set
\begin{equation}\label{eqn:Ckamz}
C(k,a,m;z)
:= \frac{1}{a^{\lambda_{k,a,m}} \Gamma(\lambda_{k,a,m}+1)
            \,\mathopen|\sinh(z)|^{\lambda_{k,a,m}+1}}.
\end{equation}
With these notations, we have:
\begin{lem}\label{lem:Lm}
For $z\in\C^+ \setminus i\pi\Z$,
the kernel function $\Lambda_{k,a}^{(m)}(r,s;z)$ has the following
upper estimate:
\begin{equation}\label{eqn:Lkamz}
|\Lambda_{k,a}^{(m)} (r,s;z)|
\le C(k,a,m;z)(rs)^m
\exp \Bigl( -\frac{1}{a}(r^a+s^a) \alpha(z) (1-|\beta(z)|)\Bigr).
\end{equation}
\end{lem}

\begin{proof}
By the following upper estimate of the $I$-Bessel function 
(see \cite[Lemma 8.5.1]{KM}) 
\begin{equation}\label{eqn:4.9}
\vert \widetilde{I}_\nu(w)\vert\leq \Gamma(\nu+1)^{-1}
e^{\vert \Re(w)\vert}
\quad\text{for $\nu\geq -\frac{1}{2}$  and $w\in \C$}
\end{equation} 
that we get
\begin{equation}\label{eqn:4.10}
\vert \Lambda_{k,a}^{(m)}(r,s;z)\vert\leq 
C(k,a,m;z)(rs)^m 
\exp \Bigl( -\frac{1}{a} (r^a+s^a) \Bigl( \Re \coth(z) -
     \Bigl| \Re \frac{1}{\sinh(z)} \Bigr| \Bigr) \Bigr).
\end{equation} 
Here, we have used  $r^a+s^a\geq 2(rs)^{\frac{a}{2}}$. 
Then, the substitution of \eqref{eqn:az} and \eqref{eqn:abz} shows
Lemma.
\end{proof}

We are ready to complete the proof of the following:
\begin{thm}\label{thm:4.2}
 Let $\gamma_z = \Exp\Bigl(iz\begin{pmatrix}0&1\\-1&0\end{pmatrix}\Bigr)$ 
be an element of $\widetilde{\Gamma(W)}$ (see \eqref{3.27}),
and $\Lambda_{k,a}^{(m)}(r,s;z)$ the function defined by
\eqref{eqn:4.6}.
Assume
 $m\in \N$ and $a>0$ satisfy
\begin{equation}\label{eqn:mkNa}
2m+2\langle k\rangle+N+a-2 > 0.
\end{equation}
Then, for $z\in\C^{++}$,
the Hilbert--Schmidt operator $\Omega_{k,a}^{(m)}(\gamma_z) $ on $L^2(\R_+, r^{2\langle k\rangle+N+a-3}dr)$ is given by 
\begin{equation}\label{eqn:4.7}
 \Omega_{k,a}^{(m)}(\gamma_z) f(r)= \int_0^\infty \Lambda_{k,a}^{(m)}( r,s;z) f(s)s^{2\langle k\rangle+N+a-3}ds.\end{equation}
The integral in \eqref{eqn:4.7}   converges absolutely for $f\in L^2(\R_+, s^{2\langle k\rangle+N+a-3}ds).$ 
\end{thm}

\begin{proof}
We have already proved the formula \eqref{eqn:4.6} for
$\Lambda_{k,a}^{(m)}(r,s;z)$.
The convergence of the integral \eqref{eqn:4.7} is deduced from
the Cauchy--Schwarz inequality because
 $\Lambda_{k,a}^{(m)}(r,\cdot\,;z) \in L^2(\R_+, s^{2\langle
k\rangle+N+a-3}ds)$ for all $z\in \C^{++}$ if \eqref{eqn:mkNa} is
fulfilled.
\end{proof}

\subsubsection{The case $\Re(z)=0$}\hfill\\
The operator $\Omega_{k,a}^{(m)}(\gamma_z)$ is unitary if
$\Re(z)=0$.
In this subsection,
we discuss its distribution kernel.

We note that the substitution of $z=i\mu$  into
\eqref{eqn:4.6} makes sense as far as $\mu \notin \pi\Z$, 
and we get the following formula
\begin{equation}\label{eqn:4.8}
\Lambda_{k,a}^{(m)} (r,s;i\mu)
= \exp \Bigl( -\frac{\pi i}{2} (\lambda_{k,a,m}+1)\Bigr)
\frac{(rs)^{-\langle k\rangle-\frac{N}{2}+1}}{\sin\mu}
\exp \Bigl( \frac{i}{a} (r^a+s^a) \cot(\mu) \Bigr)
J_{\lambda_{k,a,m}} 
\left( \frac{2}{a} \frac{(rs)^{\frac{a}{2}}}{\sin\mu} \right).
\end{equation}
Here, we have used the relation
$I_\lambda \Bigl(\frac{z}{i}\Bigr) = e^{-\frac{i\lambda\pi}{2}} J_\lambda(z)$.

In this subsection, we shall prove:
\begin{thm}\label{thm:4.4}
 Retain the notation and the assumption \eqref{eqn:mkNa} as in Theorem \ref{thm:4.2}. For $ \mu\in  \R\setminus \pi \Z ,$ the unitary operator $\Omega_{k,a}^{(m)}(\gamma_{i\mu})$ on $L^2(\R_+, r^{2\langle k\rangle+N+a-3}dr)$ is given by \begin{equation}\label{eqn:4.13}
 \Omega_{k,a}^{(m)}(\gamma_{i\mu}) f(r)= \int_0^\infty f(s) \Lambda_{k,a}^{(m)}( r,s;i\mu) s^{2\langle k\rangle+N+a-3}ds.\end{equation}
The integral in the right-hand side 
\eqref{eqn:4.13} converges absolutely for all  $f$ in
 the dense subspace, in $L^2(\R_+, r^{2\langle k\rangle+N+a-3}dr),$
 spanned by the functions $\big\{f_{\ell,m}^{(a)}\big\}_{\ell\in \N}$ 
(see \eqref{eqn:4.4} for definition).
\end{thm}
\begin{proof}
Let   $\epsilon >0 $ and $\mu\in \R\setminus \pi \Z.$ By Theorem \ref{thm:4.2} we have \begin{equation}\label{eqn:4.12}
\Omega_{k,a}^{(m)}(\gamma_{\epsilon+{i\mu}}) f(r)= \int_0^\infty  f(s)\Lambda_{k,a}^{(m)}( r,s;\epsilon +i\mu)s^{2\langle k\rangle+N+a-3}ds.\end{equation}
As $\epsilon \rightarrow  0$ the left-hand side converges to
$\Omega_{k,a}^{(m)}(\gamma_{i\mu})$ by Theorem
\ref{thm:3.18} (1). 

On the other hand,  the addition formula 
$${\rm csch}(\epsilon+i\mu)={{{\rm csch}(\epsilon){\rm csch}(i\mu)}\over {\coth(\epsilon)+\coth(i\mu)}}$$     gives 
\begin{equation}\label{eqn:shxy}
\vert {\rm csch}( \epsilon+i\mu)\vert <\vert {\rm csch}(i\mu)\vert.
\end{equation}
Hence, it follows from Lemma \ref{lem:Lm} that we have
$$\vert \Lambda_{k,a}^{(m)}(r,s; \epsilon+i\mu)\vert 
\leq C(k,a,m;i\mu)(rs)^m 
\leq C {{(rs)^{m}}\over {\vert \sin (\mu)\vert^{\lambda_{k,a,m}+1}}},$$ 
for some constant $C$.
In view of \eqref{eqn:4.4}, we get
$$\big\vert \Lambda_{k, a}^{(m)}(r,s;  \epsilon+i\mu) f_{\ell,m}^{(a)}(s)\big\vert \leq C' {{r^{m}\exp(-\frac{1}{a}s^a)}\over {\vert\sin(\mu)\vert^{\lambda_{k,a,m}+1}}} s^{2m} \Big\vert L_\ell^{(\lambda_{k,a,m})}\Big({2\over a}s^a\Big)\Big\vert.$$ 
Now, we can use the dominated  convergence theorem to deduce that the right-hand side of \eqref{eqn:4.12} goes to 
$$ \int_0^\infty f (s)\Lambda_{k,a}^{(m)}(r,s; i\mu)  s^{2\langle k\rangle +N+a-3}ds $$ as $\epsilon \rightarrow 0.$ 
Hence, Theorem has been proved.
\end{proof}

As a corollary of Theorem \ref{thm:4.4},
 we obtain  representation theoretic proofs of the following two
 classical integral formulas of Bessel functions:
\begin{coro}\label{cor:Web}
~
 \begin{itemize}
\item[(1)] 
(Weber's second exponential integral, \cite[6.615]{GR}).
$$\int_0^\infty
 e^{-\delta T }  J_\nu(2\alpha\sqrt  T) J_\nu(2\beta \sqrt T) dT={1\over {\delta}} e^{-{1\over {\delta}}(\alpha^2+\beta^2)} I_\nu\Big( {{2\alpha\beta}\over {\delta}}\Big),$$ where $\vert \arg(\delta)\vert <{\pi\over 2}$ and $\nu\geq 0.$
\item[(2)]
(See \cite[7.421.4]{GR}).
\begin{equation}\label{eqn:web2}
\int_0^\infty e^{-\delta T } L_\ell ^{(\nu)}(\alpha T ) J_\nu(\beta \sqrt T) T^{\frac{\nu}{2}} dT={{(\delta-\alpha)^\ell \beta^\nu}\over {
2^{\nu} \delta^{\nu +\ell+1} }}
e^{ -{{\beta^2}\over {4\delta}}} L_\ell^{(
\nu)}\Big( {{\alpha\beta^2}\over {4\delta(\alpha-\delta)}}\Big),
\end{equation}
for $\Re(\delta)>0$ and\/ $\Re(\nu)>0.$ 
\end{itemize}
\end{coro}

\begin{proof}[Proof (Sketch)]
(1) 
The semigroup  law $\Omega_{k,a}^{(m)} (\gamma_{z_1}) \Omega_{k,a}^{(m)} (\gamma_{z_2}) =\Omega_{k,a}^{(m)} (\gamma_{z_1+z_2})$ yields
\begin{equation}\label{eqn:Web}
\int_0^\infty
\Lambda_{k,a}^{(m)}(r,s;z_1)\Lambda_{k,a}^{(m)}(s, r';z_2) s^{2\langle
k\rangle+N+a-3}ds =\Lambda_{k,a}^{(m)}(r, r';z_1+z_2).
\end{equation}
Using the
expression  \eqref{eqn:4.8} of $\Lambda_{k,a}^{(m)},$  we get the
identity (1).

(2)
The identity \eqref{eqn:4.5} (in terms of group theory, this comes from the
$K$-type formula \eqref{eqn:rec a}) can be restated as
$$\int_0^\infty \Lambda_{k,a}^{(m)}(r,s; z) f_{\ell,m}^{(a)}(s) s^{2\langle k\rangle+N+a-3}ds= e^{-z(2\ell+\lambda_{k,a,m}+1)} f_{\ell,m}^{(a)}(r),$$
in terms of the integral kernels
by Theorem \ref{thm:4.2}.
 After some simplifications and by putting  constants together, we get 
the identity \eqref{eqn:web2}.
\end{proof}

\begin{rem}\label{rem:4.5}
\begin{enumerate}[\upshape 1)]
\item  
The operator $\pi(\lambda_{k,a,m})(\gamma_z)$ acts on the irreducible
representation $\pi(\lambda_{k,a,m})$ of $\widetilde{SL(2,\R)}$
as a scalar
multiplication if $z\in i\pi\Z$ i.e.\ if $\gamma_z$ belongs to the
center (see Fact \ref{fact:SL2rep}).
Correspondingly, the kernel function
 $\Lambda_{k,a}^{(m)}(r,s;\epsilon+i\mu)$
approaches to a scalar multiple of Dirac's delta function as $\epsilon$ goes to $0$ if
$\mu\in\pi\Z$. 
\item  
Of particular interest is another case where 
$\mu \in \pi(\Z+\frac{1}{2})$.
For simplicity, let
$\mu=\frac{\pi}{2}$.
Then,
the formula \eqref{eqn:4.8} collapses to
\[
\Lambda_{k,a}^{(m)}\Bigl(r,s;\frac{\pi i}{2}\Bigr)
= \exp \Bigl( -\frac{\pi i}{2} (\lambda_{k,a,m}+1)\Bigr)
  (rs)^{-\langle k\rangle -\frac{N}{2}+1}
  J_{\lambda_{k,a,m}} \Bigl(\frac{2}{a}(rs)^{\frac{a}{2}}\Bigr).
\]
For $a=1,2$, we have
\begin{alignat*}{3}
&\Lambda_{k,1}^{(m)}(r,s;{\pi i\over 2})
&&= e^{-i{\pi\over 2}(2m +2\langle k\rangle+N-1)}(rs)^{-\langle k\rangle-{N\over 2}+1} J_{2m +2\langle k\rangle +N-2}(2\sqrt{rs})
&\quad&(a=1),
\\
&\Lambda_{k,2}^{(m)}(r,s;{\pi i\over 2})
&&= e^{-i{\pi\over 2}( m+ \langle k\rangle+{N\over 2})}(rs)^{-\langle
    k\rangle-{N\over 2}+1} J_{ m + \langle k\rangle +{N\over 2}-1}( {rs})
&&(a=2).
\end{alignat*}
\end{enumerate}
\end{rem}

We shall discuss the unitary operator $\Omega(\gamma_{{\pi i}\over 2})=\lim_{\epsilon \rightarrow 0}\Omega(\gamma_{\epsilon+{\pi i}\over 2})$ 
 in full detail,
 which we call the
$(k,a)$-\textit{generalized Fourier transform} $\Fka$ (up to a phase factor)
in Section \ref{sec:5}.

\subsection{Gegenbauer transform}
\hfill\break
In this section, we summarize some basic properties of the Gegenbauer
polynomials and the corresponding integral transforms.
\subsubsection{The Gegenbauer polynomial}
\hfill\break
The Gegenbauer polynomial $C_m^\nu(t)$ of degree $m$ is defined by the
 generating function
\begin{equation}\label{eqn:Ggen}
(1-2rt+r^2)^{-\nu}
= \sum_{m=0}^\infty C_m^\nu (t)r^m.
\end{equation}
To be more explicit, it is given as 
\begin{equation}\label{eqn:4.14}
\index{Ca_m^\nu(t)@$C_m^\nu(t)$|main}%
C_m^\nu(t) 
=  \big(-{1\over 2}\big)^m{{(2\nu)_m}\over {m ! (\nu+{1\over 2})_m}}
(1-t^2)^{-\nu+\frac{1}{2}}{{d^m}\over
{dt^m}} (1-t^2)^{m+\nu-\frac{1}{2}}.
\end{equation}

If we put $t=\cos\theta$,
and expand
\[
(1-2r\cos\theta+r^2)^{-\nu}
= ((1-r e^{i\theta})(1-r e^{-i\theta}))^{-\nu}
\]
by the binomial theorem,
then we have
\begin{equation}\label{eqn:Cpol}
C_m^\nu(\cos\theta)
= \sum_{k=0}^m
  \frac{(\nu)_k (\nu)_{m-k}}{k! (m-k)!} \cos(m-2k)\theta.
\end{equation}
Then, the following fact is readily seen.
\begin{fact}\label{fact:Gegen1}
\begin{enumerate}[\upshape 1)]
\item  
$C_m^\nu(t)$ is a polynomial of $t$ of degree $m$,
and is also a polynomial in parameter $\nu$.
\item  
$C_m^0(t) \equiv 0$ \  for any $m\ge1$.
\item  
$C_0^\nu(t) \equiv 1$, \  $C_1^\nu(t) = 2\nu t$.
\item  
$C_m^\nu(1) = \dfrac{\Gamma(m+2\nu)}{m! \Gamma(2\nu)}$.
\end{enumerate}
\end{fact}

In this subsection, we prove:
\begin{lem}\label{lem:Gbdd}
Fix $\nu\in\R$.
Then there exists a constant $B(\nu) >0$ such that
\begin{equation}\label{eqn:Cmnew}
\sup_{-1\le t\le 1} 
\left| \frac{1}{\nu} C_m^\nu(t) \right|
\le B(\nu) m^{2\nu-1}
\quad\text{for any $m\in\N_+$}.
\end{equation}
\end{lem}

\begin{rem}\label{rem:Gbdd}
\begin{enumerate}[\upshape 1)]
\item  
For $\nu>0$, it is known that the upper bound of
$|C_m^\nu(t)|$ is attained at $t=1$,
namely, 
\[
\sup_{-1\le t\le 1} \left| C_m^\nu(t) \right|
= C_m^\nu(1)
= \frac{\Gamma(m+2\nu)}{m! \Gamma(2\nu)}.
\]
This can be verified easily by \eqref{eqn:Cpol}, or alternatively, 
by the following integral expression for $\nu>0$:
\begin{equation}\label{eqn:Geint}
C_m^\nu(t)
= \frac{\Gamma(\nu+\frac{1}{2})\Gamma(m+2\nu)}
       {\sqrt{\pi}m!\Gamma(\nu)\Gamma(2\nu)}
  \int_0^\pi (t+\sqrt{t^2-1} \cos\theta)^m
  \sin^{2\nu-1}\theta \, d\theta.
\end{equation}
\item  
For $\nu=0$, the left-hand side of \eqref{eqn:Cmnew} is interpreted as
the $L^\infty$-norm of\/ $\lim_{\nu\to0} \frac{1}{\nu} C_m^\nu(t)$,
which is a polynomial of $t$ by Fact \ref{fact:Gegen1} 2).
\item  
Our proof below works also for all $\nu\in\C$.
\end{enumerate}
\end{rem}

Before proving Lemma \ref{lem:Gbdd},
we prepare the following estimate:
\begin{claim}\label{clm:Gam}
Let $\lambda\in\R$.
Then there exists a constant $A(\lambda)>0$ such that
\[
\left| \frac{\Gamma(\lambda+k)}{\Gamma(\lambda)k!} \right|
\le A(\lambda)k^{\lambda-1}
\quad\text{for any $k\in\N$}.
\]
\end{claim}

\begin{proof}
We recall Stirling's asymptotic formula of the Gamma function:
\[
\Gamma(x) \sim \Gamma_0(x)
\quad\text{as $x\to\infty$},
\]
where $\Gamma_0(x) := \sqrt{2\pi} x^{x-1} e^{-x}$.
In light of the following ratio:
\[
\frac{\Gamma_0(k+\lambda)}{\Gamma_0(k+1)}
= k^{\lambda-1} 
\Biggl( \frac{1+\frac{\lambda}{k}}{1+\frac{1}{k}} \Biggr)^{k+\frac{1}{2}}
\left( 1+\frac{\lambda}{k} \right)^{\lambda-1} e^{1-\lambda},
\]
we get
\[
\lim_{k\to\infty}\frac{\Gamma(\lambda+k)}{k! k^{\lambda-1}} 
=1.
\]
Thus, Claim follows.
\end{proof}

\begin{proof}[Proof of Lemma \ref{lem:Gbdd}]
By \eqref{eqn:Cpol},
we have
\begin{align*}
\sup_{-1\le t\le 1} \vert C_m^\nu(t)\vert
&\le \sum_{k=0}^m
 \left|
 \frac{\Gamma(\nu+k)}{k! \Gamma(\nu)}
 \frac{\Gamma(m-k+\nu)}{(m-k)! \Gamma(\nu)}
 \right|.
\\
\intertext{We note that there is no pole in the Gamma factors in the
right-hand side. We now use Claim \ref{clm:Gam},
and get}
&\le \sum_{k=0}^m A(\nu)^2 k^{\nu-1} (m-k)^{\nu-1}
\\
&\le A(\nu)^2 (m+1) \left( \frac{m}{2} \right)^{2\nu-2}.
\end{align*}
Hence, \eqref{eqn:Cmnew} is proved for $\nu\ne 0$.

For $\nu=0$, we use \eqref{eqn:Cpol}, and get
\begin{equation}\label{eqn:Cm0}
\lim_{\nu\to0} \frac{1}{\nu} C_m^\nu(\cos\theta)
= \frac{1}{m} (\cos(n-2m)\theta+\cos n\theta).
\end{equation}
Hence, the inequality \eqref{eqn:Cmnew} also holds for $\nu=0$.
Thus, we have proved Lemma \ref{lem:Gbdd}.
\end{proof}

\subsubsection{The Gegenbauer transform}
\hfill\break
We summarize $L^2$-properties of Gegenbauer polynomials in a way
that we shall use later.   

\begin{fact}[see \cite{AAR}, {\cite[Chapter 15]{DB}}]\label{fact:Gegen}
Suppose $\nu>-\frac{1}{2}$.
\begin{enumerate}[\upshape 1)]
\item  
$\{ C_m^\nu(t): m\in\N \}$
is an orthogonal basis in the Hilbert space
$\mathcal{H}_\nu :=
 L^2((-1,1),(1-t^2)^{\nu-\frac{1}{2}} dt)$.
\item  
$\displaystyle \int_{-1}^1 C_m^\nu(t)^2(1-t^2)^{\nu-\frac{1}{2}} dt
 =  \frac{\pi\Gamma(2\nu+m)}
         {2^{2\nu-1}\Gamma(m+1)(m+\nu)\Gamma(\nu)^2}
$.
\item  
We set a normalized constant $b_{\nu,m}$ by
\begin{equation}\label{eqn:balC}
b_{\nu,m}
:= \frac{2^{2\nu-1}\Gamma(m+1)\Gamma(\nu)\Gamma(\nu+1)}
        {\pi\Gamma(m+2\nu)}.
\end{equation}
Then, the Gegenbauer transform defined by
\begin{equation}\label{eqn:Ialh}
\index{Cs_{\nu,m}@$\mathcal{C}_{\nu,m}$|main}%
\mathcal{C}_{\nu,m}(h)
:= b_{\nu,m} \int_{-1}^1 h(t) C_m^\nu(t)
  (1-t^2)^{\nu-\frac{1}{2}} dt,
\quad\text{for $h\in \mathcal{H}_\nu$}
\end{equation}
has the following inversion formula.
\begin{equation}\label{eqn:invCt}
h = \frac{1}{\nu} \sum_{m=0}^\infty (m+\nu)
    \mathcal{C}_{\nu,m}(h) C_m^\nu(t).
\end{equation}
\end{enumerate}
\end{fact}

The orthonality relation in Fact \ref{fact:Gegen} can be restated in
terms of the Gegenbauer transform as follows:
\begin{equation}\label{eqn:bal}
\mathcal{C}_{\nu,m}(C_n^\nu) 
=
\begin{cases}
 \dfrac{\nu}{m+\nu}
 &(n=m),
\\[1.5ex]
 0
 &(n\ne m).
\end{cases}
\end{equation}

The Gegenbauer transform $\mathcal{C}_{\nu,m}$ arises also in 
 a Dunkl analogue of 
the classical Funk--Hecke formula for spherical harmonics as follows.
\begin{fact} [see {\cite[Theorem 2.1]{X1}}]\label{fact:4.8}
 Let  $h$ be a continuous function  on $[-1,1]$ and $p\in \Hkm(\R^N).$  
Then, we have
  \begin{equation}\label{eqn:4.19}
 d_k  \int_{S^{N-1}} (\widetilde{V}_kh)(\omega,\eta) p(\eta)
  \vartheta_k(\eta)d\sigma(\eta)
= \mathcal{C}_{\langle k\rangle+\frac{N-2}{2},m} (h) p(\omega),\qquad
  \omega\in S^{N-1}.\end{equation} 
Here, $d_k$ is the constant defined in \eqref{eqn:dk},
and $\widetilde{V}_kh$ is defined
  in \eqref{eqn:Vkhx} by using the Dunkl intertwining operator $V_k$.
For $k\equiv 0,$
$(\widetilde{V}_kh)(\eta,\omega) = h(\langle\eta,\omega\rangle)$ and
  $\vartheta_k(\eta)\equiv1$,  
so that the identity \eqref{eqn:4.19} collapses  to the original Funk--Hecke formula. 
\end{fact}

\subsubsection{Explicit formulas of Gegenbauer transforms}
\hfill\break
In this subsection,
we present two explicit formulas of Gegenbauer transforms
$\mathcal{C}_{\nu,m}$ (see \eqref{eqn:Ialh}).
These results will be used in describing the kernel distributions
$\Lambda_{k,a}(x,y;z)$ for $a=1,2$ (see Theorem \ref{thm:a12}).
\begin{lem}\label{lem:key}
~
\begin{enumerate}[\upshape 1)]
\item  
$\mathcal{C}_{\nu,m} \Bigl(\widetilde{I}_{\nu-\frac{1}{2}}
\Bigl(\alpha(1+t)^{\frac{1}{2}}\Bigr) \Bigr)
= 2^{2\nu-m} \pi^{-\frac{1}{2}} \alpha^{2m} \Gamma(\nu+1)
  \widetilde{I}_{2m+2\nu} (\sqrt{2} \, \alpha)
= \frac{\alpha^{2m}\Gamma(2\nu+1)}{2^m\Gamma(\nu+\frac{1}{2})}
  \, \widetilde{I}_{2m+2\nu} (\sqrt{2} \alpha)$.
\item  
$\mathcal{C}_{\nu,m} (e^{\alpha t})
= 2^{-m} \alpha^m \Gamma(\nu+1) \widetilde{I}_{\nu+m}(\alpha)$.
\end{enumerate}
\end{lem}
This lemma is an immediate consequence of the following integral
formulas and the duplication formula of the Gamma function: 
\begin{equation}\label{eqn:dup}
\Gamma(2\nu)
= 2^{2\nu-1} \pi^{-\frac{1}{2}} \Gamma(\nu) \Gamma\Bigl(\nu+\frac{1}{2}\Bigr).
\end{equation}
\begin{lem}\label{lem:4.7}
 For $\alpha,\nu\in \C$ such that $\Re(\nu)>0,$ the following two integral formulas hold: 
\begin{itemize}
\item[(1)] 
\begin{equation}\label{eqn:4.16}
  \int_{-1}^1\widetilde I_{\nu-{1\over 2}} \big(\alpha(1+t)^{\frac12}\big)C_{m }^{\nu}(t) (1-t^2)^{ \nu-{1\over 2}}dt=
{2\sqrt \pi {\Gamma(m +2\nu)}\over {\Gamma(m +1)\Gamma\big(\nu \big)}}\Big({{ \alpha}\over {\sqrt 2}}\Big)^{2m}\widetilde I_{2m+2\nu}\big(\sqrt 2\alpha\big).\end{equation}  
\item[(2)] \begin{equation}\label{eqn:4.17}
 \int_{-1}^1e^{ \alpha t} C_m^{ \nu}(t) (1-t^2)^{\nu-\frac12} dt= {{2^{ -2\nu-m+1}\pi  \Gamma(2\nu+m)}\over {\Gamma(m+1) \Gamma(\nu)}}\alpha^{m} \widetilde I_{\nu+m}(\alpha). \end{equation}  
\end{itemize}
\end{lem} 

\begin{proof}[Proof of Lemma \ref{lem:4.7}] 
(1) 
This identity was proved in \cite[Lemma 8.5.2]{KM}.

(2) The  integral formula \eqref{eqn:4.17} is well known (see for instance \cite[page 570]{Velin}). However, for the convenience of the readers, we give a simple proof. Using \eqref{eqn:4.14}, we have 
\begin{eqnarray*}
 \int_{-1}^1 e^{ \alpha t} C_m^\nu(t) (1-t^2)^{\nu-\frac12} dt 
& =&{{\Gamma(\nu+{1\over 2})}\over {\Gamma(2\nu)}}{{\Gamma(m+2\nu)}\over {2^m m!\Gamma(m+\nu+{1\over 2})}}\int_{-1}^1 {{d^m}\over {dt^m}}(e^{ \alpha t})(1-t^2)^{m+\nu-{1\over 2}} dt \\
& =&{{\Gamma(\nu+{1\over 2})}\over {\Gamma(2\nu)}}{{\Gamma(m+2\nu)}\over {2^m m!\Gamma(m+\nu+{1\over 2})}} \alpha^m \int_{-1}^1   e^{ \alpha t} (1-t^2)^{m+\nu-{1\over 2}} dt. 
\end{eqnarray*} Now, \eqref{eqn:4.17}  follows from the integral
representation \eqref{eqn:Itildeint} of $\widetilde{I}_\nu(w)$.
\end{proof}

\begin{rem}[expansion formula]\label{rem:expa}
Applying the inversion formula of the Gegenbauer transform (see Fact
\ref{fact:Gegen}),
we get the following expansion formulas from Lemma \ref{lem:key}:
\begin{gather}\label{eqn:ewt}
e^{wt}=\Gamma(\nu) \sum_{m=0}^\infty (\nu+m) \Big({w\over 2}\Big)^{m}\widetilde I_{\nu+m}(w) C_m^\nu(t),\qquad \Re(\nu)>0
\\
\widetilde I_{\nu-{1/2}}\Big({{w(1+t)^{1/2}}\over {\sqrt 2}}\Big)={{2^{2\nu} \Gamma(\nu) }\over \sqrt \pi}\sum_{m=0}^\infty   (\nu+m)  \Big({w\over 2}\Big)^{2m} \widetilde I_{2m+2\nu}( w)  C_m^{\nu}(t),\qquad \Re(\nu)>0 
\label{eqn:4.37}
\end{gather} 
or equivalently,
\begin{equation*}
 w^{2\nu} \widetilde J_{\nu-{1/2}}\Big({{w(1+t)^{1/2}}\over {\sqrt 2}}\Big) 
=  {{2^{4\nu} \Gamma(\nu)}\over {\sqrt \pi}} \sum_{m=0}^\infty (-1)^m ( \nu+ m) J_{2\nu+2m}(w) C_m^\nu(t).
\end{equation*}
The first formula \eqref{eqn:ewt} is 
Gegenbauer's expansion (see for instance \cite[7.13(14)]{W}),
whereas the second expansion formula \eqref{eqn:4.37} was proved in Kobayashi--Mano
\cite[Proposition 5.7.1]{KM}.
\end{rem}

\subsection{Integral representation for
$\Omega_{k,a}(\gamma_z)$}
\label{sec:4.3}
\hfill\\
In this section,
we find the integral kernel $\Lambda_{k,a}(x,y;z)$ of the operator
$\Omega_{k,a}(\gamma_z)$ for $z\in\C^+ \setminus i\pi\Z$.
The main result is Theorem \ref{thm:4.9}.

\subsubsection{The function $\mathcal{I}(b,\nu;w;t)$}\hfill\\
In this subsection,
we introduce a function $\mathcal{I}(b,\nu;w;t)$ of four
variables,
and study its basic properties.

Let $\widetilde{I}_\lambda(w) = (\frac{w}{2})^{-\lambda}I_\lambda(w)$
 be the normalized $I$-Bessel function
(see \eqref{eqn:4.15}),
and $C_m^\nu(t)$ the Gegenbauer polynomial.
Consider the following infinite sum:
\begin{equation}\label{eqn:Ibw}
\mathcal{I}(b,\nu;w;t)
= \frac{\Gamma(b\nu+1)}{\nu}
  \sum_{m=0}^\infty (m+\nu)
  \left( \frac{w}{2} \right)^{bm}
  \widetilde{I}_{b(m+\nu)}(w) C_m^\nu(t).
\end{equation}
We note that $\nu=0$ is not a singularity in the summand because
$C_m^0(t)\equiv 0$ for $m \ge 1$ (see Fact \ref{fact:Gegen1} 2);
see also \eqref{eqn:Cm0}).
In this subsection, we prove:
\begin{lem}\label{lem:Ib0}
\begin{enumerate}[\upshape 1)]
\item  
The summation \eqref{eqn:Ibw} converges absolutely and uniformly on
any compact subset of 
\begin{equation}\label{eqn:1bnu}
U := \{(b,\nu,w,t) \in \R_+\times \R \times \C \times [-1,1]:
1+b\nu > 0 \}.
\end{equation}
In particular,
$\mathcal{I}(b,\nu;w;t)$ is a continuous function on $U$.
\item  
(Special value at $w=0$)
\begin{equation}\label{eqn:Ib0}
\mathcal{I}(b,\nu;0;t) \equiv 1.
\end{equation}
\item  
(Gegenbauer transform) \ 
For $\nu > -\frac{1}{2}$,
\[
\mathcal{C}_{\nu,m} (\mathcal{I}(b,\nu;w;\cdot))
= \Gamma(1+b\nu) \left(\frac{w}{2}\right)^{bm}
  \widetilde{I}_{b(m+\nu)} (w),
\quad \text{for $m\in\N$}.
\]
\end{enumerate}
\end{lem}

\begin{proof}
1)
It is sufficient to show that for a sufficiently large $m_0$
the summation over $m$ $(\ge m_0)$
converges absolutely and uniformly
on any compact set of $U$.
We recall from \eqref{eqn:4.9} and \eqref{eqn:Cmnew} that
\begin{align*}
|\widetilde{I}_\lambda(w)|
&\le \frac{e^{|\Re w|}}
          {\Gamma(\lambda+1)},
\\
\left| \frac{1}{\nu} C_m^\nu(t) \right|
&\le B(\nu) m^{2\nu-1}
\quad\text{for any $m\ge1$}.
\end{align*}
Then,
\begin{align*}
\frac{1}{\nu} \sum_{m=m_0}^\infty
 \Biggl| (m+\nu) \left(\frac{w}{2}\right)^{bm}
  \widetilde{I}_{b(m+\nu)} (w) C_m^\nu(t) \Biggr|
&\le \sum_{m=m_0}^\infty \left| 
     \frac{(m+\nu)B(\nu)w^{bm} e^{|\Re w|}m^{2\nu-1}}
          {2^{bm} \Gamma(bm+b\nu+1)} \right|
\\
&= B(\nu) e^{\vert\Re w\vert}
   \sum_{m=m_0}^\infty
   \frac{(1+\frac{\nu}{m}) m^{2\nu}}
        {\Gamma(bm+b\nu+1)}
   \biggl( \left|\frac{w}{2}\right|^b \biggr)^m.
\end{align*}
Since $b>0$, 
$\Gamma(bm+b\nu+1)$ grows faster than any other term in each summand
as $m$ goes to
infinity,
and consequently, 
the last sum converges.
Furthermore, the convergence is uniform on any compact set of
parameters $(b,\nu,w)$.
Hence, we have proved the first assertion.

2)
Since $b>0$, 
the summand in \eqref{eqn:Ibw} vanishes at $w=0$ for any $m>0$,
and therefore
\begin{align*}
\mathcal{I}(b,\nu;0;t)
&= \frac{\Gamma(b\nu+1)}{\nu} \cdot \nu \cdot
   \widetilde{I}_{b\nu}(0) \cdot C_0^\nu (t)
\\
&=1.
\end{align*}
Thus, the second assertion is proved.

3)
This is an immediate consequence of Fact \ref{fact:Gegen} on the
Gegenbauer transform $\mathcal{C}_{\nu,m}$.
\end{proof}

\begin{exa}\label{ex:Ibnwt}
The special values at
$b=1,2$ are given by
\begin{alignat}{2}
&\mathcal{I}(1,\nu;w;t) 
&&= e^{wt},
\label{eqn:I1t}
\\
&\mathcal{I}(2,\nu;w;t) 
&&= \Gamma\Bigl(\nu+\frac{1}{2}\Bigr) \, \widetilde{I}_{\nu-\frac{1}{2}}
  \Biggl( \frac{w(1+t)^{\frac{1}{2}}}{\sqrt{2}} \Biggr).
\label{eqn:I2t}
\end{alignat}
\end{exa}

\begin{proof}[Proof of Example \ref{ex:Ibnwt}]
First, let us prove the identity \eqref{eqn:I1t}.
By Lemma \ref{lem:Ib0} 3), we have
\[
\mathcal{C}_{\nu,m}(\mathcal{I}(b,\nu;w;\cdot))
= \Gamma(1+b\nu) \left(\frac{w}{2}\right)^{bm} \widetilde{I}_{b(m+\nu)} (w),
\quad\text{for all $m\in\N$}.
\]
By Lemma \ref{lem:key} 2), we have
\begin{equation*}
\mathcal{C}_{\nu,m}
( e^{wt} )
= \Gamma(1+\nu)
\left(\frac{w}{2}\right)^{m} \widetilde{I}_{\nu+m}(w).
\end{equation*}
This shows that
\[
\mathcal{C}_{\nu,m}\text{(left-hand side)}
=
\mathcal{C}_{\nu,m}\text{(right-hand side)}
\quad
\text{for all $m\in\N$}
\]
with regard to the identity \eqref{eqn:I1t}.
Since the Gegenbauer polynomials form a complete orthogonal basis in the
Hilbert space $L^2((-1,1),(1-t^2)^{\nu-\frac{1}{2}}dt)$ 
(see Fact \ref{fact:Gegen}),
we have proved \eqref{eqn:I1t}.

The proof for the identity \eqref{eqn:I2t} goes similarly by using
Lemma \ref{lem:key} 1).
Thus, Example \ref{ex:Ibnwt} has been shown.
\end{proof}

\subsubsection{The normalization constant}\hfill\\
For $a>0$ and a multiplicity function $k$ on the root system $\mathcal{R}$, 
we introduce the following normalization constant
\begin{equation}\label{eqn:cka}
\index{cx_{k,a}@$c_{k,a}$|main}%
c_{k,a} 
:=\Big(\int_{\R^N} \exp\Bigl(-\frac{1}{a}\Vert x\Vert^a\Bigr) \vartheta_{k,a}(x)  dx\Big)^{-1}
\end{equation}  
where $\vartheta_{k,a}$ is the density defined in \eqref{eqn:01}.
Using the polar coordinates, we have  
\begin{align*}
c_{k,a}^{-1}
&= \int_0^\infty \int_{S^{N-1}}
   \exp \Bigl(-\frac{1}{a} r^a\Bigr)
   r^{2\langle k\rangle+N+a-3} \vartheta_k(\omega) d\sigma(\omega)dr
\\
&= d_k^{-1} \int_0^\infty e^{-t} (at)^{\frac{N+2\langle k\rangle-2}{a}} dt.
\end{align*}
Here, $d_k^{-1}$ is the $k$-deformation of the volume of the unit
sphere (see \eqref{eqn:dk}).
For a non-negative multiplicity function $k$,
 the integral converges if $2\langle k\rangle+N+a-2>0$,
and we get
\begin{equation}\label{eqn:cka-dk}
c_{k,a}= a^{-({{2\langle k\rangle+N-2}\over a})} \Gamma({{2\langle k\rangle+N+a-2}\over a})^{-1}d_k.
\end{equation}
For $k\equiv 0$,
we have $d_0^{-1} = \frac{2\pi^{\frac{N}{2}}}{\Gamma(\frac{N}{2})}$
(see \eqref{eqn:d0}),
and in particular,
\[
c_{0,1} = \frac{\Gamma(\frac{N}{2})}{2\pi^{\frac{N}{2}}\Gamma(N-1)}
= \frac{1}{(4\pi)^{\frac{N-1}{2}}\Gamma(\frac{N-1}{2})},
\quad
c_{0,2} = \frac{1}{(2\pi)^{\frac{N}{2}}}.
\]

\subsubsection{Definition of $h_{k,a}(r,s;z;t)$ and $\Lambda_{k,a}(x,y;z)$}
\hfill\break
We now
introduce the following continuous function of $t$ on the interval
$[-1,1]$ with
parameters
$r,s > 0$, and $z \in \C^+ \setminus i \pi \Z$:
\begin{equation}\label{eqn:hka}
h_{k,a}(r,s;z;t)
:=
\frac{\exp\bigl(-\frac{1}{a}(r^a+s^a)\coth(z)\bigr)}
     {\sinh(z)^{\frac{2\langle k\rangle+N+a-2}{a}}}
\mathcal{I} 
\left( \frac{2}{a},\frac{2\langle k\rangle+N-2}{2};
       \frac{2(rs)^{\frac{a}{2}}}{a\sinh(z)}; t \right).
\end{equation}

We observe that,
for $\mu \in \R \setminus \pi\Z$,
the substitution $z=i\mu$ into \eqref{eqn:hka} yields:
\begin{equation}\label{eqn:hkaim}
h_{k,a}(r,s;i\mu;t)
= \frac{\exp(\frac{i}{a}(r^a+s^a)\cot(\mu))}
       {e^{\frac{2\langle k\rangle+N+a-2}{2a}\pi i}\sin(\mu)^{\frac{2\langle k\rangle+N+a-2}{a}}}
  \mathcal{I}
  \left( \frac{2}{a}, \frac{2\langle k\rangle+N-2}{2};
         \frac{2(rs)^{\frac{a}{2}}}{ai\sin(\mu)}; t
  \right).
\end{equation}
We recall from \eqref{eqn:4.2} that
$\Lambda_{k,a}^{(m)}(r,s;z)$ is the integral kernel of the operator
$\Omega_{k,a}^{(m)}(\gamma_z)$ on
$L^2(\R_+,r^{2\langle k\rangle+N+a-3} dr)$.
Up to a constant factor (independent of $m$),
the Gegenbauer transform of $h_{k,a}$ coincides with
$\Lambda_{k,a}^{(m)}(r,s;z)$:
\begin{lem}\label{lem:Ghka}
Suppose $2\langle k\rangle+N>1$.
Then, for every $m\in\N$, we have
\begin{align}\label{eqn:Ghka}
\mathcal{C}_{\langle k\rangle+\frac{N}{2}-1,m}
(h_{k,a}(r,s;z;\cdot))
&= \frac{d_k}{c_{k,a}} \Lambda_{k,a}^{(m)} (r,s;z)
\\
&= a^{\frac{2\langle k\rangle+N-2}{a}}
   \Gamma \left( \frac{2\langle k\rangle+N+a-2}{a} \right)
   \Lambda_{k,a}^{(m)} (r,s;z).
\nonumber
\end{align}
\end{lem}

\begin{proof}
We observe 
\[
1+b\nu = \frac{2\langle k\rangle+a+N-2}{a}
\quad\text{and}\quad
\lambda_{k,a,m} = \frac{2}{a} (m+\nu)
\]
if $b=\frac{2}{a}$ and $\nu=\frac{2\langle k\rangle+N-2}{2}$.
Then, Lemma \ref{lem:Ghka} follows from Lemma \ref{lem:Ib0} 3)
and the definition \eqref{eqn:4.6}  of
$\Lambda_{k,a}^{(m)} (r,s;z)$.
\end{proof}

We are ready to define the following function on
$\R^N\times\R^N\times(\C^{+}\setminus i\pi\Z)$ by 
\begin{equation}\label{eqn:4.22}
\Lambda_{k,a}(r\omega,s\eta;z)
:= (\widetilde{V}_k h_{k,a}(r,s;z;\cdot))(\omega,\eta),
\end{equation}
where $\widetilde{V}_k$ is introduced in \eqref{eqn:Vkhx} by using 
the Dunkl intertwining operator $V_k$ and 
$h_{k,a}(r,s;z;t)$ is defined in \eqref{eqn:hka}.

\subsubsection{Expansion formula}\hfill\\
For $a>0$, we will derive a series representation for the kernel $\Lambda_{k,a}$ in terms of $\Lambda_{k,a } ^{(m)}$ and the Poisson kernel of the space 
$\Hkm(\R^N).$

In light of the definitions of $\mathcal{I}(b,\nu;w;t)$ (see
\eqref{eqn:Ibw}) and $\Lambda_{k,a}^{(m)}(r,s;z)$ (see
\eqref{eqn:4.6}), 
we may rewrite \eqref{eqn:hka} as 
\begin{equation}\label{eqn:hLmd} 
h_{k,a}(r,s;z;t)=a^{({{2\langle k\rangle+N-2}\over a})}\Gamma\left( {2\langle k\rangle+N+a-2}\over a\right)\sum_{{m}\in \N} \Lambda_{k,a }^{(m)} (r,s;z) 
\Big( {{\langle k\rangle +m +{{N-2}\over 2}}\over {\langle k\rangle  +{{N-2}\over 2}}}\Big) 
C_{m } ^{\langle k\rangle+{{N-2}\over 2}}(t).\end{equation}
The above expansion formula \eqref{eqn:hLmd} is
the series expansion by Gegenbauer polynomials (Fact \ref{fact:Gegen})
corresponding to Lemma \ref{lem:Ghka}.

Now, applying the operator $\widetilde{V}_k$ to \eqref{eqn:hLmd},
we get
\begin{thm} \label{thm:4.14}
 For $a>0$ and  $z\in \C^+\setminus i\pi \Z ,$  we have
$$\Lambda_{k,a} (x,y;z)=a^{({{2\langle k\rangle+N-2}\over a})}\Gamma\left( {2\langle k\rangle+N+a-2}\over a\right)  \sum_{{m}\in \N} \Lambda_{k,a }^{(m)} (r,s;z) P_{k,   {m} }(\omega, \eta)  $$ where 
$x=r\omega$, $y=s\eta$, and
\begin{equation}\label{eqn:Poisson}
\index{Pa_{k, {m} }(\omega, \eta)@$P_{k, {m} }(\omega, \eta)$|main}%
P_{k, {m} }(\omega, \eta):=\Big( {{\langle k\rangle +m
+{{N-2}\over 2}}\over {\langle k\rangle  +{{N-2}\over 2}}}\Big)
\Bigl(\widetilde{V}_kC_m^{\langle k\rangle+\frac{N-2}{2}}\Bigr)(\omega,\eta).
\end{equation}
\end{thm}
\noindent
In Theorem \ref{thm:4.14}, the function $P_{k,m}(\omega,\eta)$ on
$S^{N-1} \times S^{N-1}$ 
is the \textit{Poisson kernel}, 
or the reproducing kernel, of the space  of spherical $k$-harmonic   polynomials of degree $m$, 
which is characterized by the
following proposition.
\begin{prop}\label{prop:Poisson}
$P_{k,m}(\omega,\eta)$ is the kernel function of the projection from
the Hilbert space
$L^2(S^{N-1},\vartheta_k(\eta)d\eta)$
to $\Hkm(\R^N)$, namely, 
for any $p\in\mathcal{H}_k^n(\R^N)$,
\[
d_k \int_{S^{N-1}} P_{k,m} (\omega,\eta) p(\eta) \vartheta_k(\eta)
d\sigma(\eta)
= \begin{cases}
      p(\omega)   & (n=m),
  \\
      0           & (n\ne m).
  \end{cases}
\]
\end{prop}

\begin{exa}\label{ex:Poisson}
For $N=1$,
$S^{N-1}$ consists of two points $(\pm 1)$,
and $\Hkm(\R^N) = 0$ if $m \ge 2$.
In this case, it is easy to see
\begin{align*}
& d_k = \frac{1}{2},
\\
&P_{k,m}(\omega,\eta) =
\begin{cases} 1 & (m=0), \\ \sgn(\omega\eta) & (m=1). \end{cases}
\end{align*}
\end{exa}

\begin{proof}[Proof of Proposition \ref{prop:Poisson}]
By the Funk--Hecke formula in the Dunkl setting (see Fact
\ref{fact:4.8}),
we have
\begin{align*}
d_k \int_{S^{N-1}}
\left( \widetilde{V}_k C_m^{\langle k\rangle+\frac{N-2}{2}} \right)
(\omega,\eta) p(\eta) \vartheta_k(\eta) d\sigma(\eta)
&= \mathcal{C}_{\langle k\rangle+\frac{N-2}{2},n}
   \left( C_m^{\langle k\rangle+\frac{N-2}{2}} \right) p(\omega)
\\
&= \begin{cases}
      \frac{\langle k\rangle+\frac{N-2}{2}}
           {m+\langle k\rangle+\frac{N-2}{2}}
      p(\omega) 
      & (n = m),
   \\
      0
      & (n \neq m).
   \end{cases}
\end{align*}
Here, we have used Fact \ref{fact:Gegen} and \eqref{eqn:bal}
for the last equality.
Hence, Proposition \ref{prop:Poisson} follows.
\end{proof}

\subsubsection{Integral representation of $\Omega_{k,a}(\gamma_z)$}\hfill\\
We are ready to prove the main result of this section.
Recall from Theorem \ref{thm:3.18} that $\Omega_{k,a} (\gamma_z) $
is a holomorphic semigroup consisting of Hilbert--Schmidt operators on 
$L^2(\R^N,\vartheta_{k,a}(x)dx)$ for $\Re z > 0$,
and is a one-parameter subgroup of unitary operators for $z\in i\,\R $.
Here is an integral representation of $\Omega_{k,a} (\gamma_z)$:
\begin{thm}\label{thm:4.9}
Suppose $a>0$ and $k$ is a non-negative 
multiplicity function on the root system
 $\mathcal{R}$ satisfying
\begin{equation}\label{eqn:kN12}
2\langle k\rangle+N > \max(1,2-a).
\end{equation}
\begin{enumerate}[\upshape 1)]
\item  
Suppose $\Re z > 0$.
Then, 
the Hilbert--Schmidt operator $\Omega_{k,a} (\gamma_z)$ on $L^2\big(\R^N,
 \vartheta_{k,a}(x)dx\big)$ is given by
 \begin{equation}\label{eqn:4.24}\Omega_{k,a} (\gamma_z)f(x)=c_{k,a}
 \int_{\R^N}f(y)  \Lambda_{k,a}(x,y;z)   
 \vartheta_{k,a}(y)dy,
\end{equation} 
where   $c_{k,a}$ is the constant defined in \eqref{eqn:cka} and 
the kernel function $\Lambda_{k,a}(x,y;z)$ is defined in \eqref{eqn:4.22}.
\item  
Suppose $z = i\mu$ $(\mu \in \R \setminus \pi \Z)$.
Then,  
the unitary operator  $\Omega_{k,a} (\gamma_{i\mu})$ on $L^2\big(\R^N, \vartheta_{k,a}(x)dx )$ is given by \begin{equation}\label{eqn:4.29}
 \Omega_{k,a} (\gamma_{i\mu})f(x)=c_{k,a} \int_{\R^N} f(y)\Lambda_{k,a}(x,y;i\mu)   \vartheta_{k,a}(y)dy.\end{equation}
\end{enumerate}
\end{thm} 

\begin{proof}
Thanks to Proposition \ref{pro:Lm},
it suffices  to show the following identity:
$$c_{k,a}\int_{S^{N-1}} \Lambda_{k,a} (r\omega,s\eta;z)
p(\eta)\vartheta_k(\eta)d\sigma(\eta) = \Lambda_{k,a } ^{(m)}(r,s;z)p(\omega),
$$
for all $p\in \Hkm(\R^N)$
and $m\in \N$. 
This follows from Theorem \ref{thm:4.14},
Proposition \ref{prop:Poisson},
and \eqref{eqn:cka-dk}.
Hence, Theorem \ref{thm:4.9} is proved.
\end{proof}

\subsection{The $a=1,2$ case}\label{sec:a12}
\hfill\break
As we have seen in Theorem \ref{thm:4.9},
 the kernel function $\Lambda_{k,a}(x,y;z)$ for the
holomorphic semigroup $\Omega_{k,a}(\gamma_z)$ is given as
\[
\Lambda_{k,a}(r\omega,s\eta;z)
= (\widetilde{V}_k h_{k,a}(r,s;z;\cdot)) (\omega,\eta).
\]
See \eqref{eqn:Vkhx} for the definition of $\widetilde{V}_k$.
In this section, 
we give a closed formula of $h_{k,a}(r,s;z;t)$ 
 for $a=1,2$,
and discuss the convergence of the integral \eqref{eqn:4.24} in Theorem
\ref{thm:4.9}.

\subsubsection{Explicit formula for $h_{k,a}(r,s;z;t)$ $(a=1,2)$}
\hfill\break
When $a=1,2,$  the series expansion in  \eqref{eqn:hLmd} can be expressed in terms of elementary  functions as follow. 
\begin{thm}\label{thm:a12}
Let $\langle k\rangle$ be defined in \eqref{eqn:2.6}, and
$\widetilde{I}_\nu$ the normalized $I$-Bessel function (see
\eqref{eqn:4.15}). 
Then, for $z\in\C^+\setminus i\pi\Z,$ we have: 
\begin{align}\label{eqn:4.21}
\index{hx_{k,a}(r,s;z;t)@$h_{k,a}(r,s;z;t)$|main}%
h_{k,a}(r,s;z;t)
={}&
 \frac{\exp(-\frac{1}{a}(r^a+s^a)\coth(z))}
                  {\sinh(z)^{\frac{2\langle k\rangle+N+a-2}{a}}}
\nonumber
\\
&\qquad \times
\begin{cases}
   \displaystyle
   \Gamma\Bigl(\langle k\rangle+\frac{N-1}{2}\Bigr)
   \widetilde{I}_{\langle k\rangle+\frac{N-3}{2}}
   \Bigl( \frac{\sqrt{2}(rs)^{\frac{1}{2}}}{\sinh z}
          (1+t)^{\frac{1}{2}} \Bigr)
   & (a=1),
 \\
   \displaystyle
   \exp \Bigl( \frac{rst}{\sinh z} \Bigr)
   & (a=2).
\end{cases}
\end{align}  
\end{thm}

\begin{proof}
In view of the definition \eqref{eqn:hka} of $h_{k,a}(r,s;z;t)$,
Theorem \ref{thm:a12} follows from formulas \eqref{eqn:I1t} and
\eqref{eqn:I2t} for $\mathcal{I}(a,\nu;w;t)$ in Example \ref{ex:Ibnwt}.
\end{proof}

\subsubsection{Absolute convergence of integral
representation}
\hfill\break
By using Theorem \ref{thm:a12},
we shall give an upper bound for $\Lambda_{k,a}(x,y;z)$.
We begin with the following:
\begin{lem}\label{lem:Iineq}
For $b=1,2$
\begin{equation}\label{eqn:Iineq}
|\mathcal{I}(b,\nu;w;t)|
\le e^{|\Re w|}
\end{equation}
for any $t\in[-1,1]$, $\nu>0$ and $w\in\C$.
\end{lem}

\begin{proof}
We have seen in \eqref{eqn:I1t} and \eqref{eqn:I2t} the explicit
formulas of $\mathcal{I}(b,\nu;w;t)$ for $b=1,2$.
Then \eqref{eqn:Iineq} is obvious for $b=1$,
and follows from the upper estimate \eqref{eqn:4.9} of the $I$-Bessel
function for $b=2$. 
\end{proof}

\begin{prop}\label{prop:Lineq}
Suppose $b$ is a positive number, 
for which the inequality \eqref{eqn:Iineq}
holds.
Let $a := \frac{2}{b}$. 
Then the function $\Lambda_{k,a}(x,y;z)$ (see \eqref{eqn:4.22})
satisfies the following inequalities:
\begin{enumerate}[\upshape 1)]
\item  
For $\Re z>0$, there exist positive constants $A, B$ depending on $z$
such that
\begin{equation}\label{eqn:Lz1}
\vert \Lambda_{k,a}(x,y;z)\vert
\le A \exp \left( -B(\Vert x\Vert^a + \Vert y\Vert^a) \right),
\quad\text{for any $x,y\in\R^N$}.
\end{equation}
\item  
For $z=i\mu+\epsilon$ $($$\mu\in\R \setminus \pi\Z$, $\epsilon\ge0$$)$,
\begin{equation}\label{eqn:Lz2}
\vert \Lambda_{k,a}(x,y;i\mu+\epsilon)\vert
\le \frac{1}{\vert \sin(\mu)\vert^{\frac{N+2\langle k\rangle+a-2}{a}}}.
\end{equation}
\end{enumerate}
\end{prop}

\begin{rem}
By Lemma \ref{lem:Iineq},
the assumption of Proposition \ref{prop:Lineq} is fulfilled 
for $b=1,2$.
We do not know if \eqref{eqn:Iineq} holds for $b$ other than $1$
and $2$.
\end{rem}

\begin{proof}
Suppose the inequality \eqref{eqn:Iineq} of
$\mathcal{I}(b,\nu;w;t)$ holds.
Then, by the definition of $h_{k,a}(r,s;z;t)$ in \eqref{eqn:hka},
the inequality \eqref{eqn:Iineq} brings us to the following estimate:
\[
|h_{k,a}(r,s;z;t)|
\le
\frac{1}{|\sinh(z)|^{\frac{2\langle k\rangle+N+a-2}{a}}}
\exp\Bigl(-\frac{1}{a}(r^a+s^a)\alpha(z)\Bigr)
\exp\Bigl(\frac{2}{a}(rs)^{\frac{a}{2}} \alpha(z)\beta(z)\Bigr).
\]
Here, we have used the functions $\alpha(z)$ and $\beta(z)$ defined in Lemma \ref{lem:abz}.
Since $r^a+s^a \ge 2(rs)^{\frac{a}{2}}$,
we have obtained:
\[
|h_{k,a}(r,s;z;t)|
\le
\frac{1}{|\sinh(z)|^{\frac{N+2\langle k\rangle+a-2}{a}}}
\exp\Bigl(-\frac{1}{a}(r^a+s^a) \alpha(z)(1-|\beta(z)|)\Bigr).
\]
We recall that $\Lambda_{k,a}(x,y;z)$ is defined by applying the
operator $\widetilde{V}_k$ to $h_{k,a}(r,s;z;\cdot)$
(see \eqref{eqn:4.22}).
Then, it follows from Proposition \ref{prop:VkB} that
\[
|\Lambda_{k,a}(x,y;z)|
\le
\Vert h_{k,a}(r,s;z;\cdot)\Vert_{L^\infty}
\le
\frac{1}{|\sinh(z)|^{\frac{N+2\langle k\rangle+a-2}{a}}}
\exp\Bigl(-\frac{1}{a} (\Vert x\Vert^a + \Vert y\Vert^a)
\alpha(z) (1-|\beta(z)|)\Bigr).
\]
Suppose now that $\Re z >0$.
Then, $\alpha(z)>0$ and $|\beta(z)|<1$ by Lemma
\ref{lem:abz}.
Hence, we have proved \eqref{eqn:Lz1}.

On the other hand,
suppose $z=i\mu+\epsilon$ $($$\mu\in\R \setminus \pi\Z$,
$\epsilon\ge0$$)$.
Then $\alpha(z)\ge0$, $\vert\beta(z)\vert<1$ by Lemma \ref{lem:abz},
and as we have seen in \eqref{eqn:shxy}
\[
\vert\sinh z\vert \ge \vert \sin \mu \vert.
\]
Hence, we have shown \eqref{eqn:Lz2}.
\end{proof}

Now we are ready to prove:
\begin{coro}\label{coro:4.10}
Suppose we are in the setting of Theorem \ref{thm:4.9}.
Let $a=1,2$.
\begin{enumerate}[\upshape 1)]
\item  
For $\Re z>0$,
 the right-hand side of 
\eqref{eqn:4.24} converges absolutely for any
$f\in L^2\big(\R^N, \vartheta_{k,a}(x)dx \big)$.
\item  
For $z=i\mu \in i(\R \setminus \pi\Z)$,
the right-hand side of \eqref{eqn:4.29}
converges absolutely for all
$f \in (L^1 \cap L^2)(\R^N,\vartheta_{k,a}(x)dx)$.
\end{enumerate}
\end{coro}

\begin{proof}
1) 
It follows from Proposition \ref{prop:Lineq} 1) that
\[
\Lambda_{k,a}(x,y;z)
\in L^2(\R^N \times \R^N, \vartheta_{k,a}(x) \vartheta_{k,a}(y)dx \, dy)
\]
for $\Re z >0$.
Therefore, Corollary is clear from the Cauchy--Schwartz inequality.

2) 
We shall substitute $z$ by $\epsilon+i\mu$ in \eqref{eqn:4.24} and let $\epsilon$ goes to $0$.

For the left-hand side of \eqref{eqn:4.24},
we use Theorem \ref{thm:3.18} 1), 
and get $$\lim_{\epsilon\rightarrow 0} \Omega_{k,a} (\gamma_{\epsilon+i\mu})=\Omega_{k,a}(\gamma_{i\mu}) .$$ 

For the right-hand side of \eqref{eqn:4.24}, 
thanks to Proposition \ref{prop:Lineq} 2),
we see
\[
\lim_{\epsilon\downarrow0}
\int_{\R^N}
\Lambda_{k,a}(x,y;i\mu+\epsilon)f(y)\vartheta_{k,a}(y)dy
= \int_{\R^N} \Lambda_{k,a}(x,y;i\mu)f(y)\vartheta_{k,a}(y)dy
\]
for $f\in L^1(\R^N,\vartheta_{k,a}(y)dy)$ by the Lebesgue dominated
convergence theorem.
Hence, we have shown that
\[
(\Omega_{k,a}(\gamma_{i\mu})f)(x)
= \int_{\R^N} \Lambda_{k,a}(x,y;i\mu)f(y)\vartheta_{k,a}(y)dy,
\]
and the right-hand side converges for any
$f\in (L^1 \cap L^2) (\R^N,\vartheta_{k,a}(y)dy)$.
Hence, Corollary \ref{coro:4.10} has been proved.
\end{proof}

\subsection{The rank one case}\label{subsec:4.3} \hfill\\For the one
  dimensional case, the only choice of the non-trivial reduced root
  system $\cal R$ is $\cal R=\{\pm 1\}$ in $\R$ up to scaling,
  corresponding to the Coxeter group $\CC=\{\rm id, \sigma\}\cong
  \Z/2\Z$ on $\R,$ where $\sigma(x)=-x.$ Here $\langle k\rangle=k.$
In this section we give a closed form of
$\Lambda_{k,a}(x,y;z)$ for $N=1$.

First of all, for $N=1$ and $a>0$, we note that we do not 
need Lemma \ref{lem:Ghka}, for which the assumption was
$2\langle k\rangle+N>1$.
Hence, instead of \eqref{eqn:kN12},
we simply need the following assumption:
\begin{equation}\label{eqn:k12}
a>0 \quad\text{and}\quad 2k > 1-a.
\end{equation}
The goal of this section is to
 find a closed formula of the kernel function $\Lambda_{k,a}(x,y;z)$
 (see \eqref{eqn:4.22}) for all
$a>0$ and for an arbitrary multiplicity function subject to \eqref{eqn:k12}.

\begin{prop}\label{fact:4.13}
Let $N=1$, $a>0$, $k\geq 0$ and $2k > 1-a$.
For $z\in \C^+\setminus i\pi\Z,$  the holomorphic semigroup $\Omega_{k,a}(\gamma_z)$ on $L^2(\R, \vert x\vert^{2k+a-2}dx)$ is given by 
$$\Omega_{k,a}(\gamma_z) f(x)= 2^{-1}a^{-({{2k-1}\over a})} \Gamma\Big({{2k+a-1}\over a}\Big)^{-1}
\int_\R f(y) \Lambda_{k,a}(x,y;z) \vert y\vert^{2k+a-2}dy,  $$ where 
 \begin{equation}\label{eqn:4.30}
\Lambda_{k,a}(x, y ;z) =\Gamma\Big( {{2k+a-1}\over a}\Big) 
{{e^{-{1\over a}(\vert x\vert^a+\vert y\vert^a)\coth(z)}}\over {\sh(z)^{{2k+a-1}\over a}}}
\Big(  \widetilde I_{ {2k-1}\over a}\Big( {2\over a}{{\vert xy\vert^{a\over 2}}\over {\sh(z)}}\Big)+ {1\over {a^{2\over a}}}{{xy}\over {   \sh(z) ^{2\over a}}} \widetilde  I_{ {2k+1}\over a}\Big( {2\over a}{{\vert xy\vert^{a\over 2}}\over {\sh(z)}}\Big)\Big).
\end{equation}  Here  $\widetilde I_\nu (w)$ denotes the normalized Bessel function \eqref{eqn:4.15}.
\end{prop}
\begin{proof}
By Theorem \ref{thm:4.14},
the kernel $\Lambda_{k,a}(x,y;z)$ can be recovered from a family of functions
$\{\Lambda_{k,a}^{(m)}(r,s;z): m\in\N\}$.
In the rank one case, 
$S^0$ consists of two points,
and correspondingly, 
Theorem \ref{thm:4.14} collapses to the following:
\begin{eqnarray*}
&&\Lambda_{k,a}(x, y ;z)\\
&=&{{c_{k,a}^{-1}}\over 2} \Big( \Lambda_{k,a}^{(0)}(\vert x\vert ,\vert y\vert;z)+ \Lambda_{k,a}^{(1)}(\vert x\vert ,\vert y\vert;z)\sgn(xy)\Big)\\
&=&{{c_{k,a}^{-1}}\over 2} 
{{e^{-{1\over a}(\vert x\vert^a+\vert y\vert^a)\coth(z)}}\over {\sh(z)}}
\vert xy\vert^{-k+(1/2)}\Big( I_{\lambda_{k,a,0}}\Big( {2\over a}{{\vert xy\vert^{a\over 2}}\over {\sh(z)}}\Big)+ I_{\lambda_{k,a,1}}\Big( {2\over a}{{\vert xy\vert^{a\over 2}}\over {\sh(z)}}\Big)\sgn(xy)\Big)\\
&=&{{c_{k,a}^{-1}}\over 2}
{{e^{-{1\over a}(\vert x\vert^a+\vert y\vert^a)\coth(z)}}\over {\sh(z)}}
\Big( {1\over {(a  \sh(z))^{\lambda_{k,a,0}}}} \widetilde I_{\lambda_{k,a,0}}\Big( {2\over a}{{\vert xy\vert^{a\over 2}}\over {\sh(z)}}\Big)+ {{xy}\over {(a  \sh(z))^{\lambda_{k,a,1}}}} \widetilde  I_{\lambda_{k,a,1}}\Big( {2\over a}{{\vert xy\vert^{a\over 2}}\over {\sh(z)}}\Big)\Big),
\end{eqnarray*}
where $$\lambda_{k,a,0}={{2k-1}\over a},\quad \lambda_{k,a,1}={{2k+1}\over a},\quad c_{k,a}^{-1}=2a^{({{2k-1}\over a})}\Gamma\Big( {{2k+a-1}\over a}\Big).$$
Here, we have used Example \ref{ex:Poisson} for the first equality,
and the formula  \eqref{eqn:4.8} of $\Lambda_{k,a}^{(m)}$ for the
second equality.
This finishes the proof of Proposition \ref{fact:4.13}.
\end{proof}

\section{The $(k,a)$-generalized Fourier transforms $\Fka$}
\label{sec:5}
The object of our study in this section is the $(k,a)$-generalized Fourier
transform given by
$$
\index{Fka@$\Fka$|main}%
\Fka 
 = e^{\frac{\pi i}{2}(\frac{2\langle k\rangle+N+a-2}{a})}
\exp\Bigl(\frac{\pi i}{2a}\bigl(\Vert x \Vert^{2-a}
                    \Delta_k - \Vert x \Vert^a\bigr) \Bigr).
$$
This is a unitary operator on the Hilbert space
$L^2(\R^N,\vartheta_{k,a}(x)dx)$.

As we mentioned in Introduction,
the unitary operator $\Fka$ includes some known transforms as special
cases:
\newline
\begin{tabular}{cll}
$\bullet$
&the Euclidean Fourier transform \cite{H2}
&($a=2$, $k\equiv0$),
\\
$\bullet$
&the Hankel transform \cite{KM}
&($a=1$, $k\equiv0$),
\\
$\bullet$
&the Dunkl transform $\mathcal{D}_k$
&($a=2$, $k>0$).
\end{tabular}
\newline
In this section, we develop the theory of the $(k,a)$-generalized
Fourier transform $\Fka$ for
general $a$ and $k$ by using the aforementioned idea of $\mathfrak{sl}_2$-triple.
The point of our approach is that we 
 interpret $\Fka$ not as an isolated operator but as a special value of the
unitary representation $\Omega_{k,a}$ of the simply connected, simple
Lie group $\widetilde{SL(2,\R)}$ at 
$\gamma_{\frac{\pi i}{2}}$
(see \eqref{eqn:gammaz}), 
or as the boundary value of the holomorphic semigroup.
Then, we see that many properties of the Euclidean Fourier
transforms can be extended to the $(k,a)$-generalized Fourier transform 
$\Fka$ by using
 the representation theory of $\widetilde{SL(2,\R)}$.
Our theorems for $\Fka$ include the inversion formula,
and a generalization of the 
Plancherel formula, the Hecke formula, the Bochner
formula, and the Heisenberg inequality for the uncertainty principle.

As in Diagram \ref{fig:Hs} of Introduction,
the Hilbert space $L^2(\R^N,\vartheta_{k,a}(x)dx)$
admits symmetries of
$\CC \times \widetilde{SL(2,\R)}$ for general $(k,a)$,
and even higher symmetries 
than $\CC \times \widetilde{SL(2,\R)}$ for particular values of $(k,a)$.
In fact, if $k \equiv 0$, then 
the Hilbert space $L^2(\R^N,\vartheta_{k,a}(x)dx)$ 
is a representation space of the Schr\"{o}dinger model of the
Weil representation (see \cite{xfol} and references therein)
 of the metaplectic group $Mp(N,\R)$
for $a=2$, and the $L^2$-model of the minimal representation
(see \cite{KO}) of the conformal group $O(N+1,2)$ for $a=1$.
The special value $a=2$ has a particular meaning
 also for general $k$ in the sense
that $\mathcal{F}_{k,2}$ is equal to
the Dunkl transform 
$\mathcal{D}_k$.
How about the $a=1$ case for general $k$?
The unitary operator
\begin{equation}\label{eqn:HFk}
\index{Hs_k@$\mathcal{H}_k$|main}%
\mathcal{H}_k 
:= \mathcal{F}_{k,1}
\end{equation}
may be regarded as the Dunkl
analogue of the Hankel-type transform $\mathcal{F}_{0,1}$ (see
Diagram \ref{fig:11} in Introduction).
As we have seen in Section \ref{sec:a12},
this unitary operator $\mathcal{H}_k$ can be written  by means of the Dunkl
intertwining operator $V_k$ and the classical Bessel functions
(see Section \ref{sec:5.3}).

\subsection{$\mathcal{F}_{k,a}$ as an inversion unitary element}
\label{subsec:5.1}
\hfill\break
The $(k,a)$-generalized
Fourier transform $\Fka$ 
on
$L^2\big(\R^N, \vartheta_{k,a}(x)dx\big)$
 is defined as
\begin{equation}\label{eqn:5.1}
\Fka :=e^{i{\pi\over 2} ({{2\langle k\rangle+N+a-2}\over a})}\Omega_{k,a}(\gamma_{i{\pi\over 2}}).
\end{equation}
Here, we recall from \eqref{eqn:gammaz} that 
$$
\gamma_{\frac{\pi i}{2}} 
= \Exp \Bigl( \frac{\pi}{2i} {\bf k} \Bigr)
= \Exp \Bigl( \frac{\pi}{2}
  \begin{pmatrix}0&-1\\1&0\end{pmatrix}\Bigr)
$$
is an element of the simply connected Lie group 
$\widetilde{SL(2,\R)}$,
and from Theorem \ref{thm:3.12} that $\Omega_{k,a}$ is a unitary
representation of $\widetilde{SL(2,\R)}$ on the Hilbert space 
$L^2(\R^N,\vartheta_{k,a}(x)dx)$.

In this subsection,
we discuss basic properties of $\Fka$ for general $k$ and $a$,
which are derived from the fact that $\gamma_{\frac{\pi i}{2}}$ is a
representative of 
the non-trivial (therefore, the longest) element of the Weyl group for
$\mathfrak{sl}_2$.
\begin{thm}\label{thm:Fkauni}
Let $a>0$ and $k$ be a non-negative 
multiplicity function on the root
system $\mathcal{R}$ satisfying the inequality
$a+2\langle k\rangle+N>2$ (see \eqref{eqn:a2k}).
\begin{enumerate}[\upshape (1)]
\item 
{\rm (Plancherel formula)} The $(k,a)$-generalized Fourier transform
$\Fka: L^2(\R^N,\vartheta_{k,a}(x)dx) \to
L^2(\R^N,\vartheta_{k,a}(x)dx)$
is a unitary operator.
That is,
$\Fka$ is a bijective linear operator satisfying
\[ 
\Vert \Fka(f)\Vert_k = \Vert f\Vert_k
\quad\text{for any $f\in L^2(\R^N,\vartheta_{k,a}(x)dx)$}.
\]
\item 
We recall from \eqref{3.16} that $\Phi_\ell^{(a)}(p,\cdot)$ is a
function on $\R^N$ defined
as
\[
\Phi_\ell^{(a)}(p,x)
= p(x) L_\ell^{(\lambda_{k,a,m})} \Bigl( \frac{2}{a}\Vert x\Vert^a\Bigr)
  \exp \Bigl( -\frac{1}{a}\Vert x\Vert^a\Bigr),
\quad x \in \R^N,
\]
for $\ell,m\in\N$ and $p\in\Hkm(\R^N)$.
Then, $\Phi_\ell^{(a)}(p,\cdot)$ is an eigenfunction of 
 $\Fka$:
\begin{equation}\label{5.6}
\Fka(\Phi_{\ell }^{(a)}(p,\cdot))=e^{-i{\pi }(\ell+{m\over
a})}  \Phi_{\ell }^{(a)}(p,\cdot).\end{equation}  
\end{enumerate}
\end{thm}

\begin{proof}
Since the phase factor in \eqref{eqn:5.1} is modulus one,
the first statement is an immediate consequence of the fact that
$\Omega_{k,a}$ is a unitary representation of 
$\widetilde{SL(2,\R)}$ (see Theorem \ref{thm:3.12}).

To see the second statement,
we recall from Proposition \ref{prop:3 7} 1) and Theorem \ref{thm:3.9} that 
$\Phi_\ell^{(a)}(p,\cdot)$
is an eigenfunction of 
$\omega_{k,a}({\bf k})$.
Then, the integration of \eqref{3.21 a} leads us to the identity \eqref{5.6}.
\end{proof}

\begin{coro}\label{cor:Fkafin}
The $(k,a)$-generalized Fourier transform
$\Fka$ is of finite order if and only if $a \in \Q$.
If $a$ is of the form $a=\frac{q}{q'}$,
where $q$ and $q'$ are positive integers, 
then
$$
(\Fka)^{2q} = \id.
$$
\end{coro}
\begin{proof}
We recall from Proposition \ref{prop:3 7} 3) that
$$
W_{k,a}(\R^N)  = \C\text{-span} \{\Phi_{\ell }^{(a)}(p,\cdot)\;|\; \ell\in \N, {m}\in \N, p\in
 \Hkm(\R^N)\}
$$
is a dense subspace in
$L^2(\R^N,\vartheta_{k,a}(x)dx)$. 
Hence, it follows from \eqref{5.6} that
the unitary operator $\Fka$ is of finite order if and only if $a \in
\Q$.
If $a=\frac{q}{q'}$,
then $(\Fka)^{2q}$ acts on $\Phi_{\ell}^{(a)}(p,\cdot)$ as a scalar
multiplication by
$$
\Bigl(e^{-i\pi(\ell+\frac{m}{a})}\Bigr)^{2q} = 1
$$
for any $m \in \N$ and any
$p \in \Hkm(\R^N)$.
Thus, we have proved $(\Fka)^{2q} = \id$.
\end{proof}

Corollary \ref{cor:Fkafin} implies particularly that
$\mathcal{H}_k := \mathcal{F}_{k,1}$ (see \eqref{eqn:HFk}) is of order
two, and the Dunkl transform $\mathcal{D}_k = \mathcal{F}_{k,2}$ is of
order four.
We pin down these particular cases as follows:
\begin{thm}[inversion formula]\label{thm:5.6}
Let $k$ be a non-negative multiplicity function
 on the root system $\mathcal{R}$.
\begin{enumerate}[\upshape 1)]
\item 
Let $r$ be any positive integer.
Suppose $2\langle k\rangle+N > 2-\frac{1}{r}$.
Then, $\mathcal{F}_{k,\frac{1}{r}}$ is an involutive unitary operator
on $L^2(\R^N, \vartheta_{k,\frac{1}{r}}(x)dx)$.
Namely, the inversion formula is given by
\begin{equation}\label{eqn:inv1}
(\mathcal{F}_{k,\frac{1}{r}})^{-1}
= \mathcal{F}_{k,\frac{1}{r}}.
\end{equation}
\item 
Let $r$ be any non-negative integer.
Suppose
\[
2\langle k\rangle+N > 2-\frac{2}{2r+1}.
\]
Then, $\mathcal{F}_{k,\frac{2}{2r+1}}$ is a unitary operator of order
four on
$L^2(\R^N, \vartheta_{k,\frac{2}{2r+1}}(x)dx)$.
The inversion formula is given as
\begin{equation}\label{eqn:inv2}
(\mathcal{F}_{k,\frac{2}{2r+1}}^{-1} f)(x)
= (\mathcal{F}_{k,\frac{2}{2r+1}} f)(-x). 
\end{equation}
\end{enumerate}
\end{thm}
\begin{proof}
The first statement has been already proved.
In the second statement, we remark that the inversion formula
\eqref{eqn:inv2} is stronger than the fact that
$(\Fka)^4 = \operatorname{id}$ for $a=\frac{2}{2r+1}$.
To see \eqref{eqn:inv2},
we use \eqref{5.6} to get
\begin{align*}
(\Fka)^2 \Phi_\ell^{(a)} (p,\cdot)
&= \exp(-m(2r+1)\pi i) \Phi_\ell^{(a)} (p,\cdot)
\\
&= (-1)^m \Phi_\ell^{(a)} (p,\cdot)
\end{align*}
if $a = \frac{2}{2r+1}$.
Since $p(-x) = (-1)^m p(x)$ for   $p\in \mathcal H_k^m(\R^N)$,
we have shown that \eqref{eqn:inv2} holds for any
$f \in W_{k,a}(\R^N)$.
Since $W_{k,a}(\R^N)$ is dense in
$L^2(\R^N,\vartheta_{k,a}(x)dx)$,
we have proved \eqref{eqn:inv2}.
\end{proof}

\begin{rem}\label{rem:5.7}
Theorem \ref{thm:5.6} 2) for $r=0$ (i.e.\ $\mathcal{F}_{k,2}=\mathcal{D}_k$,
the Dunkl transform) was proved  in
Dunkl \cite{D2},  and followed by de Jeu \cite{dJ} where  the
author proved the inversion formula 
 for $\C^+$-valued root multiplicity
functions $k$.
Our approach based on the $SL_2$ representation theory gives a new proof 
of the inversion formula  and the Plancherel formula  for $\Fka$ even
for $a=2$.
\end{rem}

\begin{rem}\label{rem:Fkauni}
We recall from Theorem \ref{thm:3.16}
 that $L^2(\R^N,\vartheta_{k,a}(x)dx)$ decomposes into a
discrete direct sum of irreducible unitary representations of 
$G=\widetilde{SL(2,\R)}$. 
Hence,
the square
$(\Fka)^2$ acts as a scalar multiple on each summand of
\eqref{eqn:Hmkpi} by Schur's lemma 
because $\gamma_{\pi i}$ is a central element of $G$ (see
\eqref{eqn:ZSL}) and
$
(\Fka)^2 = e^{i\pi(\frac{2\langle k\rangle+N+a-2}{a})}
 \Omega_{k,a}(\gamma_{\pi i})
$
by \eqref{eqn:5.1}. 
Since $\gamma_{\pi i}$ acts on the irreducible representation
$\pi(\lambda_{k,a,m})$ as a scalar
$e^{-\pi i(\lambda_{k,a,m}+1)}$ by Fact \ref{fact:SL2rep} (5),
$\Fka^2$ acts on it as the scalar
$$
e^{i\pi(\frac{2\langle k\rangle+N+a-2}{a})}
e^{-\pi i(\lambda_{k,a,m}+1)}
= e^{-\frac{2m\pi i}{a}}.
$$

This gives us an alternative proof of Corollary \ref{cor:Fkafin}.
\end{rem}

Next, we discuss intertwining properties of the $(k,a)$-generalized
Fourier transform $\Fka$ with differential
operators. 
Let $E = \sum_{j=1}^N x_j\partial_j$ be the Euler operator on
$\R^N$ as before. 
\begin{thm}\label{thm:Fka}
The unitary operator $\mathcal{F}_{k,a}$ 
satisfies the following intertwining relations on a dense subspace of
$L^2(\R^N,\vartheta_{k,a}(x)dx)$:
\begin{enumerate}[\upshape 1)]
\item 
$\mathcal{F}_{k,a} \circ E 
= -(E+N+2\langle k\rangle+a-2) \circ \mathcal{F}_{k,a}$.
\item 
$\mathcal{F}_{k,a} \circ \Vert x\Vert^a
= -\Vert x\Vert^{2-a} \Delta_k \circ \mathcal{F}_{k,a}$.
\item 
$\mathcal{F}_{k,a} \circ \Vert x\Vert^{2-a} \Delta_k
= -\Vert x\Vert^a \circ \mathcal{F}_{k,a}$.
\end{enumerate}
These identities hold in the usual sense, and also in the distribution sense in the space of distribution vectors of the unitary representation of $G$ on $L^2(\R^N, \vartheta_{k,a}(x)dx).$
\end{thm}

If we use $\xi$ (instead of $x$) for the coordinates of the target space of $\Fka$,
we may write Theorem \ref{thm:Fka} 2) and 3) as 
\begin{subequations}
  \renewcommand{\theequation}{\theparentequation\ \alph{equation}}%
\begin{equation}\label{5.4 b}
 \Fka( \Vert \cdot \Vert^a f)(\xi)=- \Vert \xi\Vert^{2-a} \Delta_k \Fka(f)(\xi),
\end{equation}
\begin{equation}\label{5.4 a}
\Fka(\Vert \cdot\Vert^{2-a} \Delta_k f)(\xi)=-\Vert \xi\Vert^a  \Fka(f)(\xi).
\end{equation}
\end{subequations}

\begin{proof}[Proof of Theorem \ref{thm:Fka}]
We observe that $\gamma_{\frac{\pi i}{2}}$ is a representative of the
longest Weyl group element,
and satisfies the following relations in $\mathfrak{sl}_2$:
$$
{\rm Ad}(\gamma_{\frac{\pi i}{2}}){\bf h} = -{\bf h}, \ 
{\rm Ad}(\gamma_{\frac{\pi i}{2}}){\bf e}^+ = -{\bf e}^-, \ 
{\rm Ad}(\gamma_{\frac{\pi i}{2}}){\bf e}^- = -{\bf e}^+,
$$
(see \eqref{eqn:3.1} for the definition of ${\bf e}^+$, ${\bf e}^-$, and 
${\bf h}$).
In turn,
we apply the identity
$$
\Omega_{k,a}(g) \omega_{k,a}(X) \Omega_{k,a}(g)^{-1}
= \omega_{k,a}(\operatorname{Ad}(g)X),
\quad (g\in G, \  X \in \g),
$$
to 
$\EE_{k,a}^+ = \omega_{k,a}({\bf e}^+)$,
$\EE_{k,a}^- = \omega_{k,a}({\bf e}^-)$, and 
$\HH_{k,a} = \omega_{k,a}({\bf h})$ (see \eqref{eqn:3.5}).
Then we have
\begin{align}
&\mathcal{F}_{k,a} \circ \HH_{k,a} = -\HH_{k,a} \circ \mathcal{F}_{k,a},
\label{eqn:HHka}
\\
&\mathcal{F}_{k,a} \circ \EE_{k,a}^+ = -\EE_{k,a}^- \circ \mathcal{F}_{k,a},
\nonumber
\\
&\mathcal{F}_{k,a} \circ \EE_{k,a}^- = -\EE_{k,a}^+ \circ \mathcal{F}_{k,a},
\nonumber
\end{align}
because $\mathcal{F}_{k,a}$ is a constant multiple of 
$\Omega_{k,a}\bigl(\gamma_{\frac{\pi i}{2}}\bigr)$
(see \eqref{eqn:5.1}).
Now, Theorem \ref{thm:Fka} is read from the explicit formulas of
$\EE_{k,a}^+$, $\EE_{k,a}^-$, and $\HH_{k,a}$ (see \eqref{eqn:3.3}).
\end{proof}

\subsection{Density of $(k,a)$-generalized Fourier transform $\Fka$}
\label{sec:Bka}
\hfill\break
By the Schwartz kernel theorem,
the unitary operator $\Fka$ can be expressed by means of a distribution
kernel.
By using the normalizing constant $c_{k,a}$ (see \eqref{eqn:cka}), 
we write the unitary
operator $\Fka$ on $L^2(\R^N,\vartheta_{k,a}(x)dx)$ as
an integral transform:
\begin{equation}\label{eqn:HB}
\Fka f(\xi) = c_{k,a} \int_{\R^N} 
\index{Ba_{k,a}(\xi,x)@$B_{k,a}(\xi,x)$|main}%
B_{k,a}(\xi,x) 
f(x)
\vartheta_{k,a}(x)dx.
\end{equation}
Comparing this with the integral expression of
$\Omega_{k,a}(\gamma_z)$ in Theorem \ref{thm:4.9},
we see that the distribution kernel $B_{k,a}(\xi,x)$ in \eqref{eqn:HB}
is given by
\begin{equation}\label{eqn:BLambda}
B_{k,a}(x,y)
= e^{i \pi (\frac{2\langle k\rangle+N+a-2}{2a})}
  \Lambda_{k,a}\Bigl(x,y; i\frac{\pi}{2}\Bigr)
\end{equation}
because $\Fka =  e^{i \pi (\frac{2\langle k\rangle+N+a-2}{2a})}\Omega_{k,a}(\gamma_{i\frac{\pi}{2}})$
 (see \eqref{eqn:5.1}).

Now, Theorem \ref{thm:Fka} is reformulated as the differential
equations that are satisfied by the distribution kernel $B_{k,a}(x,\xi)$ as follows:

\begin{thm}\label{thm:5.4}
The distribution $B_{k,a}(\cdot, \cdot)$ solves
 the following  differential-difference equation on $\R^N \times \R^N$ 
\begin{subequations}
  \renewcommand{\theequation}{\theparentequation\ \alph{equation}}%
 \begin{align}
\label{5.5 a}
&E^x B_{k,a}(\xi,x) = E^\xi B_{k,a}(\xi,x),
\\
\label{5.5 b}
&\Vert \xi\Vert^{2-a} \Delta_k^\xi B_{k,a} (\xi,x)
 = -\Vert x\Vert^a B_{k,a} (\xi,x),
\\
& \Vert x\Vert^{2-a} \Delta_k^x B_{k,a}(\xi,x)= -\Vert \xi\Vert^a B_{k,a}(\xi,x).
\label{5.5 c}
\end{align}
\end{subequations}
Here, the superscript in $E^x$, $\Delta_k^x$, etc indicates the relevant variable. 
\end{thm} 
\begin{rem}\label{rem:5.5}
 For $a=2,$ Theorem \ref{thm:5.4} was previously known as the
 differential equation of the Dunkl kernel (cf. \cite{D1}).
\end{rem}

\begin{proof}[Proof of Theorem \ref{thm:5.4}]
First we use the identity
\eqref{eqn:HHka} as operators on $\R^N$ for any $a>0$ and $k$.
It is convenient to write
\[
\HH_{k,a} = \frac{1}{a} (\ell+2E),
\]
where $E$ is the Euler operator and $\ell := N+2\langle k \rangle + a-2$.
Then, by \eqref{eqn:HB}, the identity \eqref{eqn:HHka} implies
\begin{equation}\label{eqn:HHf}
\int_{\R^N} \Bigl( (\ell+2E^x)f(x)\Bigr)
B_{k,a}(\xi,x) \vartheta_{k,a}(x)dx
= -\int_{\R^N} f(x) (\ell+2E^\xi) B_{k,a}(\xi,x) \vartheta_{k,a}(x)dx
\end{equation}
for any test function $f(x)$
(i.e.\ $f(x) \vartheta_{k,a}(x)^{\frac{1}{2}} \in \mathcal{S}(\R^N)$).

Now we recall that the density $\vartheta_{k,a}(x)$ (see
\eqref{eqn:01} for definition) is homogeneous of degree 
$a-2+2\langle k \rangle$ $(= \ell-N)$, we have
\begin{equation}\label{eqn:Etheta}
E^x \vartheta_{k,a}(x)
= (\ell-N)\vartheta_{k,a}(x).
\end{equation}
On the other hand,
it follows from
$\sum_{j=1}^N x_j \frac{\partial}{\partial x_j} 
 - \sum_{j=1}^N \frac{\partial}{\partial x_j} x_j
 = -N$
as operators,
we have
\begin{equation}\label{eqn:EN}
\int_{\R^N}(E^x f)(x) g(x)dx
= -\int_{\R^N} f(x)(N+E^x) g(x)dx.
\end{equation}
Combining \eqref{eqn:Etheta} and \eqref{eqn:EN},
we have
\[
\mbox{the left-hand side of \eqref{eqn:HHf}}
= -\int_{\R^N} f(x) (\ell B_{k,a}(\xi,x) + 2E^x B_{k,a}(\xi,x))
  \vartheta_{k,a}(x)dx.
\]
Hence, the identity \eqref{eqn:HHf} implies that the
distribution  kernel $B_{k,a}(\xi,x)$ satisfies the differential equation
\begin{equation}\label{eqn:EBka}
E^x B_{k,a}(\xi,x) = E^\xi B_{k,a}(\xi,x).
\end{equation}

Next, the identity \eqref{5.4 b} implies
$$
\int_{\R^N} B_{k,a}(\xi,x) \Vert x\Vert^a f(x) \vartheta_{k,a}(x)dx
= -\Vert \xi\Vert^{2-a} \Delta_k^\xi
   \int_{\R^N} B_{k,a}(\xi,x) f(x) \vartheta_{k,a}(x)dx
$$
for any test function $f$.
Hence the second differential equation \eqref{5.5 b} follows.

Finally, by the identity \eqref{5.4 a},
we have
$$
\int_{\R^N} B_{k,a}(\xi,x) 
  \Bigl(\Vert x\Vert^{2-a} \Delta_k^x f(x)\Bigr) \vartheta_{k,a}(x)dx
= -\Vert \xi\Vert^a
   \int_{\R^N} B_{k,a}(\xi,x) f(x) dx.
$$
Since $\Vert x\Vert^{2-a} \Delta_k^x$ is a symmetric operator 
on 
$L^2(\R^N,\vartheta_{k,a}(x)dx)$,
the left-hand side is equal to
$$
\int_{\R^N} \Vert x\Vert^{2-a}
 \Bigl( \Delta_k^x B_{k,a}(\xi,x) \Bigr)
 f(x) \vartheta_{k,a}(x) dx.
$$
Hence the third differential equation \eqref{5.5 c} is proved.
\end{proof}

We continue basic properties on the kernel $B_{k,a}(\xi,x)$ of the
$(k,a)$ generalized Fourier transform.
\begin{thm}\label{thm:Bka}
~
\begin{enumerate}[\upshape 1)]
\item 
$B_{k,a}(\lambda x,\xi) = B_{k,a}(x,\lambda\xi)$ \ for $\lambda>0$.
\item 
$B_{k,a}(hx,h\xi) = B_{k,a}(x,\xi)$ \ for $h \in \CC$.
\item 
  $B_{k,a}(\xi,x)=B_{k,a}(x,\xi).$
\item 
 $B_{k,a}(0,x)=1.$
\end{enumerate}
\end{thm}

\begin{proof}
1)
This statement follows from the differential equation \eqref{5.5 a} given in Theorem
\ref{thm:5.4}.

2)
Since $\Fka$ commutes with the action of the Coxeter group $\CC$,
the second statement follows from the fact that
$\vartheta_{k,a}(x)dx$ is a $\CC$-invariant measure.

3)
Putting $\mu = \frac{\pi}{2}$ in \eqref{eqn:hkaim},
we get
\begin{equation}\label{eqn:hkapi}
h_{k,a}\Bigl(r,s;\frac{\pi i}{2}; t\Bigr)
= e^{-\frac{2\langle k\rangle+N+a-2}{2a}\pi i}
  \mathcal{I}
  \left( \frac{2}{a}, \frac{2\langle k\rangle+N-2}{2};
         \frac{2(rs)^{\frac{a}{2}}}{ai} ; t
  \right).
\end{equation}
In particular,
we have
\[
h_{k,a}\Bigl(r,s;\frac{\pi i}{2};t\Bigr)
= h_{k,a}\Bigl(s,r;\frac{\pi i}{2};t\Bigr).
\]
In view of \eqref{eqn:4.22} and Proposition \ref{prop:VkB},
we conclude that
\[
\Lambda_{k,a}\Bigl(x,y;\frac{\pi i}{2}\Bigr)
= \Lambda_{k,a}\Bigl(y,x;\frac{\pi i}{2}\Bigr).
\]
Hence, the third statement has been proved.

4)
By Lemma \ref{lem:Ib0},
$h_{k,a}(r,0;\frac{\pi i}{2};t)
= e^{-\frac{2\langle k\rangle+N+a-2}{2a} \pi i}$.
Since the Dunkl intertwining operator $V_k$ satisfies
$V_k(\bid)=\bid$ ($\bid$ is the constant function on $\R^N$)
(see (I2) in Section \ref{sec:2}),
it follows from
\eqref{eqn:4.22} that
$\Lambda_{k,a}(x,y;\frac{\pi i}{2}) 
= e^{-\frac{2\langle k\rangle+N+a-2}{2a}\pi i}$.
Finally use
\eqref{eqn:BLambda}. 
\end{proof}

\subsection{Generalized Fourier transform $\Fka$ for special values at $a=1$ and $2$}\label{sec:5.3}
\hfill\break
In this subsection we discuss closed formulas of the kernel 
$B_{k,a}(x,y)$ of the $(k,a)$-generalized Fourier transform $\Fka$
(see \eqref{eqn:HB}) in the case $a=1,2$.
The $(k,a)$-generalized Fourier transform $\Fka$ reduces to the Dunkl
operator $\mathcal{D}_k$ if $a=2$,
and gives rise to a new unitary operator $\mathcal{H}_k$,
the Dunkl analogue of the Hankel transform if $a=1$.

We renormalize the Bessel function $J_\nu$ of the first kind as
\begin{equation}\label{eqn:Jtilde}
\index{Ja_\nu(w)@$\widetilde{J}_\nu(w)$|main}%
\widetilde{J}_\nu(w) 
:= \Bigl(\frac{w}{2}\Bigr)^{-\nu}  J_\nu(w)
=  \sum_{\ell=0}^\infty
   \frac{(-1)^\ell w^{2\ell}}{2^{2\ell} \ell! \Gamma(\nu+\ell+1)}.
\end{equation}
Then, from the definition \eqref{eqn:4.15} of $\widetilde{I}_\nu(z)$
we have
\[
\widetilde{J}_\nu(w) = \widetilde{I}_\nu(-iw) = \widetilde{I}_\nu(iw).
\]
By substituting $z=\frac{\pi i}{2}$ into \eqref{eqn:4.21},
we get the following formula:
\[
h_{k,a}(r,s;\frac{\pi i}{2};t)
= \begin{cases}
     \Gamma\Bigl(\langle k\rangle+\frac{N-1}{2}\Bigr)
        e^{-\frac{\pi i}{2}(2\langle k\rangle+N-1)}
        \widetilde{J}_{\langle k\rangle+\frac{N-3}{2}}
        (\sqrt{2}(rs)^{\frac{1}{2}} (1+t)^{\frac{1}{2}})
   &(a=1),
  \\
     e^{-\frac{\pi i}{2}(\langle k\rangle+\frac{N}{2})}
        e^{-irst}
   &(a=2). 
  \end{cases}
\]
Together with \eqref{eqn:BLambda} and \eqref{eqn:4.22},
we have:
\begin{prop}\label{prop:Bka12}
In the polar
coordinates  $x=r\omega$ and $ y=s\eta,$ 
the kernel $B_{k,a}(x,y)$ is given by
\begin{equation}\label{eqn:5.3}
B_{k,a}(r\omega,s\eta)
= \begin{cases}
  \Gamma\Big(\langle k\rangle+{{N-1}\over 2}\Big) \Bigl( \widetilde{V}_k   \Big(\widetilde J _{\langle k\rangle +{{N-3}\over 2}}\big(\sqrt{ 2rs (1+\cdot)}
 \big) \Big) \Bigr) (\omega,\eta)
   &(a=1), 
  \\
  \bigl(\widetilde{V}_k\big(e^{-i rs \cdot}\big)\bigr)(\omega,\eta)
   &(a=2).
  \end{cases}
\end{equation}
\end{prop}
As one can see form \eqref{eqn:5.3}, the kernel $B_{k,2}(x,y)$ coincides with the Dunkl kernel at $(x, -iy) $ (cf. \cite{D3}).
\begin{thm}\label{thm:5.2}
 Let $k$ be a non-negative root multiplicity function, 
$a=1$ or\/ $2$, and $x, y\in \R^N.$ Then
$\vert B_{k,a}(x,y)\vert \leq 1$.
\end{thm} 
\begin{proof}
Theorem \ref{thm:5.2} follows from the special case,
i.e.\ $\mu = \frac{\pi}{2}$,
of Proposition \ref{prop:Lineq} 2) because
$|B_{k,a}(x,y)|
= |\Lambda_{k,a}(x,y;i\frac{\pi}{2})|$
by \eqref{eqn:BLambda}.
\end{proof}
\begin{rem}\label{rem:5.3}
For $a=2$ it was shown that $|B_{k,2}(x,y)|$
 is uniformly bounded for $x, y \in \mathbb R^N$
  first by de Jeu\cite{dJ} and then by R\"osler
  \cite{R1} by $1$.
\end{rem}

We note that Theorem \ref{thm:5.2} implies the absolute convergence of
 the integral defining $\Fka,$ for $a=1,2,$  on $
(L^1\cap L^2) \big(\R^N, \vartheta_{k,a}(x)dx \big)$,
as we proved in Corollary  \ref{coro:4.10}.

\subsection{Generalized Fourier transform $\Fka$ in the rank-one case}
\label{sec:5.4}
\hfill\break
This section examines $\Fka$ and its kernel $B_{k,a}(x,y)$ in the
rank-one case.

Suppose $N=1$, $a>0 $, $k\geq 0$, and $2k>1-a$.
Then, by the explicit formula of the kernel $\Lambda_{k,a}$ 
(see Proposition \ref{fact:4.13}), 
followed by the formula \eqref{eqn:BLambda}, we get
\begin{align}
B_{k,a}(x, y) &=  e^{i{\pi\over 2}({{2k+a-1}\over a})} \Lambda_{k,a}(x,y;i{\pi\over 2})\nonumber\\
&  =  
\Gamma\Big( {{2k+a-1}\over a}\Big)\Big( \widetilde J_{\frac{2k-1}{a}}\Big({2\over a}\vert xy\vert^{a\over 2}\Big)+  {{xy}\over { (ai)^{ 2\over a} }}
\widetilde J_{\frac{2k+1}{a}}\Big({2\over a}\vert xy\vert^{a\over 2}\Big)\Big),\label{5.7} 
\end{align}
 where $\widetilde J_\nu(w) = \widetilde{I}_\nu(-iw)$  is the
normalized Bessel function given in \eqref{eqn:Jtilde}; here the branch of $i^{2\over a}$ is chosen so that $1^{2\over a}=1.$ Thus, for  $a>0,$ $k\in \R^+$ such that $2k>1-a,$ and $f\in L^2(\R, \vert x\vert^{2k+a-2}dx),$ the integral transform $\Fka$ takes the form
$$\Fka f(y)= 2^{-1}a^{-({{2k-1}\over a})} \int_\R f(x) \Big( \widetilde J_{\frac{2k-1}{a}}\Big({2\over a}\vert xy\vert^{a\over 2}\Big)+ {{xy}\over {(ai)^{2\over a}}} 
\widetilde J_{\frac{2k+1}{a}}\Big({2\over a}\vert xy\vert^{a\over 2}\Big)\Big) \vert x\vert^{2k+a-2}dx.$$    
\begin{rem}\label{rem:Fka1}
If we  set 
\begin{eqnarray*}
B_{k,a}^{\rm even}(x,y)&:=&{1\over 2}\big[  B_{k,a}(x,y)+B_{k,a}(x,-y)\big]\\
&=& \Gamma\Big( {{2k+a-1}\over a}\Big)  \widetilde J_{\frac{2k-1}{a}}\Big({2\over a}\vert xy\vert^{a\over 2}\Big). 
\end{eqnarray*}
Then,  the transform $\Fka(f)$ of an even function $f$ on the real line specializes to a Hankel type transform on $\R_+.$ 
\end{rem}

Let us find the formula \eqref{5.7} by an alternative approach.
First, for general dimension $\R^N$,
by composing \eqref{eqn:BLambda}, \eqref{eqn:4.22}, and
\eqref{eqn:hkapi},
we have
\[
B_{k,a}(x,y)
= \left( \widetilde{V}_k \mathcal{I} \Bigl( \frac{2}{a},
         \frac{2\langle k\rangle+N-2}{2};
         \frac{2(rs)^{\frac{a}{2}}}{ai}; \cdot \Bigr) 
  \right)
  (\langle \omega,\eta \rangle)
\]
in the polar coordinates $x=r\omega$, $y=s\eta$.
Furthermore, in the $N=1$ case,
a closed integral formula of the Dunkl intertwining operator $V_k$ is
known: 
\begin{equation}\label{eqn:VN}
(V_kf)(x)
= \frac{\Gamma(k+\frac{1}{2})}{\Gamma(k)\Gamma(\frac{1}{2})}
  \int_{-1}^1 f(tx)(1+t)(1-t^2)^{k-1} dt,
\end{equation}
see \cite[Theorem 5.1]{D2}.
Hence, we might expect that the formula \eqref{5.7} for the kernel
$B_{k,a}(x,y)$ could be recovered directly by using the integral formula
\eqref{eqn:VN} of $V_k$.
In fact this is the case.
To see this,
we shall carry out a compution of the following integral:
\begin{equation}\label{eqn:VB12}
B_{k,a}(x,y)
= \frac{\Gamma(k+\frac{1}{2})}{\Gamma(k)\Gamma(\frac{1}{2})}
  \int_{-1}^1
  \mathcal{I} \Bigl( \frac{2}{a}, \frac{2k-1}{2},
               \frac{2(rs)^{\frac{a}{2}}}{ai};
               t\langle \omega,\eta \rangle \Bigr)
  (1+t)(1-t^2)^{k-1} dt
\end{equation}
for $x=r\omega$, $y=s\eta$ ($r,s>0$, $\omega,\eta=\pm1$).

We notice that the summation \eqref{eqn:Ibw} for
$\mathcal{I}(b,\nu;w;t)$ is taken over $m=0$ and $1$ if $N=1$.
Hence we have
\begin{equation*}
\mathcal{I} \Bigl( \frac{2}{a}, \frac{2k-1}{2};
                   \frac{2(rs)^{\frac{a}{2}}}{ai};
                   t \langle \omega,\eta \rangle
            \Bigr)
= \Gamma \Bigl( \frac{2k+a-1}{a} \Bigr)
   \Biggl( \widetilde{J}_{\frac{2k-1}{a}}
           \Bigl( \frac{2\vert xy\vert^{\frac{a}{2}}}{a} \Bigr)
          + \frac{(2k+1)txy}{(ai)^{\frac{2}{a}}}
          \widetilde{J}_{\frac{2k+1}{a}}
           \Bigl( \frac{2\vert xy\vert^{\frac{a}{2}}}{a} \Bigr)
   \Biggr).
\end{equation*}

On the other hand,
by using the integral expression of the Beta function and the
duplication formula \eqref{eqn:dup} of the Gamma function,
we have
\[
\frac{\Gamma(k+\frac{1}{2})}{\Gamma(k)\Gamma(\frac{1}{2})}
\int_0^{1} t^m (1+t)(1-t^2)^{k-1} dt
=
\begin{cases}
     1                & (m=0), \\
     \dfrac{1}{2k+1}  & (m=1).
\end{cases}
\]
Substituting these formulas into
 the right-hand side of \eqref{eqn:VB12} we have completed an
 alternative proof of \eqref{5.7}.

\subsection{Master Formula and its applications}

\subsubsection{Master Formula}
\hfill\break
We state the following two reproducing properties of the kernel
$B_{k,a} $ of basic importance. 
\begin{thm}\label{thm:5.11}
  {\rm (Master Formula)} 
 Suppose $a>0$ and $k$ is a non-negative multiplicity function satisfying $2\langle k\rangle +N> \max(1,2-a)$.
\begin{enumerate}[\upshape 1)]
\item
 For $x,y\in \R^N,$ we have
\begin{align}\label{5.11}  
&c_{k,a} \int_{\R^N} \exp\Bigl(\frac{i}{a}\Vert u\Vert^a\Bigr) B_{k,a}(x,u)B_{k,a}(u,y)   
\vartheta_{k,a}(u) du
\nonumber
\\
&=  e^{i{\pi }({{2\langle k\rangle +N+a-2}\over a})}
\exp\Bigl(-{i\over a}(\Vert x\Vert^a+\Vert y\Vert^a)\Bigr)B_{k,a}(x,y).
\end{align}
\item 
 Let  $p$ be a homogeneous polynomial on $\R^N$ of degree $m.$ Then
we have
\begin{align}\label{eqn:reprB}
c_{k,a} \int_{\R^N} \exp \Bigl( -\frac{1}{a} \Vert u\Vert^a \Bigr)
\Bigl( \exp \Bigl( -\frac{1}{2a} \Vert \cdot \Vert^{2-a} \Delta_k
            \Bigr) p \Bigr) (u)
B_{k,a}(x,u) \vartheta_{k,a}(u)du
\nonumber
\\
= e^{-\frac{im\pi}{a}} 
  \exp \Bigl( -\frac{1}{a}\Vert x\Vert^a \Bigr)
  \Bigl( \exp \Bigl( -\frac{1}{2a} \Vert \cdot \Vert^{2-a} \Delta_k
              \Bigr) p \Bigr) (x).
\end{align}
\end{enumerate}
\end{thm}

\begin{rem}\label{rem:5.12}
 For $a=2,$ the reproducing properties \eqref{5.11} and \eqref{eqn:reprB}
were previously proved in Dunkl \cite[Theorem 3.2 and Proposition 2.1]{D3}.
In that case,  theses properties  played a crucial role in 
studying Dunkl analogues of Hermite polynomials (see \cite[Section 3]{RV}),
the properties of the heat kernel associated with the heat equation
 for the Dunkl operators (see \cite[Section 4]{RV}),
and in the construction of generalized Fock spaces (see \cite[Section 3]{BO}).
\end{rem}

\subsubsection{Proof of Theorem \ref{thm:5.11}}
\hfill\break
We begin with the proof of \eqref{5.11}.  From the semigroup law
$$\Omega_{k,a}(\gamma_{z_1})\Omega_{k,a}(\gamma_{z_2})=\Omega_{k,a}(\gamma_{z_1+z_2}),
\quad\text{for $\gamma_{z_1},\gamma_{z_2}\in\widetilde{\Gamma(  W)}$},$$ 
the integral representation of $\Omega_{k,a}(\gamma_z)$ (see Theorem
\ref{thm:4.9}) yields
\begin{equation}\label{5.13} 
c_{k,a}\int_{\R^N} \Lambda_{k,a}(x,u;i{\pi\over 4})\Lambda_{k,a}(u,y;i{\pi\over 4}) \vartheta_{k,a}(u) du=\Lambda_{k,a}(x,y;i{\pi\over 2}).\end{equation} 
We set 
\[
\mu := 2\langle k\rangle+N+a-2.
\]
In view of \eqref{eqn:hka},
a simple computation shows
\[
\frac{h_{k,a}(r,s;\frac{\pi
i}{4};t)}{h_{k,a}(2^{\frac{1}{a}}r,s;\frac{\pi i}{2};t)}
= 2^{\frac{\mu}{2a}} \exp\Bigl( \frac{i}{a} (r^a+s^a) \Bigr).
\]
Applying $\widetilde{V}_k$,
and using \eqref{eqn:4.22}, we get
\begin{align}\label{eqn:Lam4}
\Lambda_{k,a} \Bigl(x,u;\frac{\pi i}{4}\Bigr)
&= 2^{\frac{\mu}{2a}}
   \exp \Bigl( \frac{i}{a} \bigl( \Vert x\Vert^a + \Vert u\Vert^a
                           \bigr) \Bigr)
   \Lambda_{k,a} \Bigl( 2^{\frac{1}{a}} x,u;\frac{\pi i}{2} \Bigr)
\nonumber
\\
&= \bigl( 2e^{-\pi i}\bigr)^{\frac{\mu}{2a}}
   \exp \Bigl( \frac{i}{a} \bigl( \Vert x\Vert^a + \Vert u\Vert^a
                           \bigr) \Bigr)
   B_{k,a} \bigl( 2^{\frac{1}{a}}x,u \bigr).
\end{align}
In the second equality,
we have used \eqref{eqn:BLambda}.
By substituting \eqref{eqn:Lam4} and \eqref{eqn:BLambda} into
\eqref{5.13}, we get
$$c_{k,a}\int_{\R^N} \exp\Bigl({{2i}\over a}\Vert u\Vert^a\Bigr) B_{k,a}(2 ^{\frac1a} x,u)B_{k,a}(2 ^{\frac1a} u,y)  \vartheta_{k,a}(u) du=
2^{-{\mu\over a}} e^{i\pi{\mu\over a}}
 \exp\Bigl(-{i\over a}(\Vert x\Vert^a+\Vert y\Vert^a)\Bigr)B_{k,a}(x,y).$$
Since $B_{k,a}(2^{\frac1a} x,u) = B_{k,a}(x,2^{\frac1a} u)$
(see Theorem \ref{thm:Bka} 1)) and $\vartheta_{k,a}(u)du$ is
homogeneous degree $N+2\langle k\rangle+a-2=\mu$,
the left-hand side equals
\[
2^{-\frac{\mu}{a}} c_{k,a}
\int_{\R^N} \exp\Bigl(\frac{i}{a} \Vert u\Vert^a\Bigr)
B_{k,a}(x,u) B_{k,a}(u,y) \vartheta_{k,a}(u)du.
\]
Hence,
\eqref{5.11} is proved.

The remaining part of this subsection is devoted to the proof of the
second statement of Theorem \ref{thm:5.11}.

We recall from \eqref{eqn:3.3} and \eqref{eqn:3.5} that
\begin{align*}
&\EE_{k,a}^+ 
= \omega_{k,a} \begin{pmatrix} 0&1\\ 0&0\end{pmatrix}
= \frac{i}{a} \Vert x\Vert^a,
\\
&\EE_{k,a}^- 
= \omega_{k,a} \begin{pmatrix} 0&0\\ 1&0\end{pmatrix}
= \frac{i}{a} \Vert x\Vert^{2-a} \Delta_k
\end{align*}
are infinitesimal generators of the unitary representation
$\Omega_{k,a}$ of  $  \widetilde{SL(2,\R)}$ on
$L^2(\R^N,\vartheta_{k,a}(x)dx)$. 

We set
\begin{equation}\label{eqn:c0}
c_0 := \Exp i \begin{pmatrix}0&1\\0&0\end{pmatrix}
    \, \Exp \frac{i}{2} \begin{pmatrix}0&0\\1&0\end{pmatrix},
\end{equation}
and introduce the operator  
\begin{equation}\label{5.10}
\index{Bs_{k,a}@$\cal B_{k,a}$|main}%
\cal B_{k,a} 
:=\exp(i\,\EE_{k,a}^+) \exp\Bigl({i\over 2}\EE_{k,a}^-\Bigr)
=\exp\Bigl(- {1\over a} { \Vert x\Vert^a }\Bigr) \exp\Bigl(-{1\over {2a}}\Vert x\Vert^{2-a}\Delta_k\Bigr).
\end{equation}  
Then, the following identity in $\mathfrak{sl}_2$, 
\begin{equation}\label{eqn:chk}
\operatorname{Ad}(c_0) \, {\bf h} 
= \begin{pmatrix}\frac{1}{2}&i\\[.5ex] \frac{i}{2}&1\end{pmatrix}
  \begin{pmatrix}1&0\\[.5ex]0&-1\end{pmatrix}
  \begin{pmatrix}\frac{1}{2}&i\\[.5ex] \frac{i}{2}&1\end{pmatrix}
  ^{-1}
= -i \begin{pmatrix}0&1\\[.5ex]-1&0\end{pmatrix}
= {\bf k},
\end{equation}
leads us to the identity of operators: 
\begin{equation}\label{eqn:5.23}
\mathcal{B}_{k,a} \circ \omega_{k,a}({\bf h})
= \omega_{k,a}({\bf k}) \circ \mathcal{B}_{k,a}.
\end{equation}
Since $\HH_{k,a} = \omega_{k,a}({\bf h})$ acts on homogeneous
functions as scalar (see \eqref{eqn:3.3}),
we know a priori that
 homogeneous functions 
applied by $\mathcal{B}_{k,a}$ are eigenfunctions of $\omega_{k,a}({\bf k})$.
Here is an explicit formula:

\begin{prop}\label{lem:Bkabasis}
For $\ell,m\in\N$ and $p\in\Hkm(\R^N)$, 
\[
\mathcal{B}_{k,a}\bigl(p(x) \Vert x\Vert^{a\ell}\bigr)
= \Bigl( -\frac{a}{2} \Bigr)^\ell \ell! \Phi_\ell^{(a)} (p,x).
\]
\end{prop}
\begin{proof}
We recall from Lemma \ref{prop:slreduction} that the linear map
\begin{equation*}
T_a: C^\infty(\R^N) \otimes C^\infty(\R_+) \to C^\infty(\R^N
         \setminus \{0\}),
\ \ 
(p,\psi) \mapsto p(x) \psi(\Vert x\Vert^a)
\end{equation*}
satisfies the following identity on
$
\Hkm(\R^N) \otimes C^\infty(\R_+)$:
\begin{equation}
\exp \Bigl(\frac{i}{2}\EE_{k,a}^- \Bigr)\circ T_a
= T_a\circ \Bigl(\operatorname{id} \otimes 
                \exp \Bigl( -\frac{a}{2}
  \Bigl(r\frac{d^2}{dr^2}+(\lambda_{k,a,m}+1)\frac{d}{dr}\Bigr)\Bigr)\Bigr).
\label{eqn:Texp}
\end{equation}
Applying \eqref{eqn:Texp} to $p\otimes r^\ell$,
and using Theorem \ref{prop:homoLag}, we get
\begin{align*}
\exp \Bigl(\frac{i}{2}\EE_{k,a}^- \Bigr) \circ T_a(p\otimes r^\ell)(x)
={}& T_a(p \otimes \Bigl(-\frac{a}{2}\Bigr)^\ell \ell! 
L_\ell^{(\lambda_{k,a,m})} \Bigl( \frac{2}{a} r \Bigr) ) (x)
\\
={}& \Bigl( -\frac{a}{2}\Bigr)^\ell \ell! p(x)
L_\ell^{(\lambda_{k,a,m})} \Bigl( \frac{2}{a} \Vert x\Vert^a\Bigr).
\end{align*}
Hence,
\begin{align*}
\mathcal{B}_{k,a}(p(x) \Vert x\Vert^{a\ell})
={}& \exp (i\,\EE_{k,a}^+) \exp \Bigl(\frac{i}{2}\EE_{k,a}^-\Bigr)
   T_a (p \otimes r^\ell)
\\
={}& \Bigl( -\frac{a}{2}\Bigr)^\ell \ell! p(x)
   \exp \Bigl( -\frac{1}{a}\Vert x\Vert^a\Bigr)
   L_\ell^{(\lambda_{k,a,m})} \Bigr( \frac{2}{a}\Vert x\Vert^a\Bigr)
\\
={}& \Bigl( -\frac{a}{2}\Bigr)^\ell \ell! \Phi_\ell^{(a)}(p,x).
\end{align*}
Thus, Proposition \ref{lem:Bkabasis} has been proved.
\end{proof}

\begin{rem}\label{rem:Thm3.19} Let $p\in\Hkm(\R^N)$. 
By \eqref{eqn:3.3},
$\omega_{k,a}({\bf h})$ acts on 
$p(x)\Vert x\Vert^{a\ell}$ by the multiplication of the scalar $\lambda_{k,a,m}+1+2\ell$.
Hence,
$\omega_{k,a}({\bf k})$ acts on $\Phi_\ell^{(a)}(p,x)$ as the same scalar
$\lambda_{k,a,m}+1+2\ell$.
This gives an alternative proof of the formula \eqref{3.21 a} in
Theorem \ref{thm:3.9}.
\end{rem}

We now introduce the vector space
\begin{equation}\label{eqn:PaRN}
\mathcal{P}_a(\R^N)
:= \text{$\C$-span} \{p(x) \Vert x\Vert^{a\ell} :
   p \in \Hkm(\R^N)\text{ for some $m\in\N$, $\ell\in\N$} \}.
\end{equation}
For $a=2$, $\mathcal{P}_2(\R^N)$ coincides with the space $\mathcal{P}(\R^N)$ 
of polynomials on $\R^N$ owing to the following algebraic
direct sum decomposition (see \cite[Theorem 5.3]{BO}):
\[
\mathcal{P}(\R^N)
\simeq \bigoplus_{m=0}^\infty
\bigoplus_{\ell=0}^{[\frac{m}{2}]}
\Vert x\Vert^{2\ell} \mathcal{H}_k^{m-2\ell} (\R^N).
\]

We introduce an endomorphism of $\mathcal{P}_a(\R^N)$,
to be denoted by $(e^{-i\frac{\pi}{a}})^*$,
as
\begin{equation}\label{eqn:pbi}
\Bigl(e^{-i\frac{\pi}{a}}\Bigr)^* \Bigl(p(x) \Vert x\Vert^{a\ell}\Bigr)
:= e^{-i(\ell+\frac{m}{a})\pi} p(x) \Vert x\Vert^{a\ell},
\quad\text{for $p \in \Hkm (\R^N)$ and $\ell\in\N$}.
\end{equation}

\begin{rem}
The notation $(e^{-i\frac{\pi}{a}})^*$ stands for the `pull-back of
functions' on the complex vector space $\C^N$ given by
\[
\Bigl(e^{-i\frac{\pi}{a}}\Bigr)^* f(z)
= f\Bigl(e^{-i\frac{\pi}{a}} z\Bigr).
\]
However, taking branches of multi-valued functions into account,
we should
 note that
$\Bigl(e^{-i\frac{\pi}{a}}\Bigr)^* \ne \id$ for $a=\frac{1}{2}$.
\end{rem}

The next proposition is needed for later use. 
 \begin{prop}\label{claim:5.10} 
 For $a>0,$  the following diagram commutes
$$ \begin{CD} \mathcal P_a(\R^N) @) \cal B_{k,a} ))     L^2(\R^N, \vartheta_{k,a}(x)dx)\\
 @ V ( e^{-i{\pi\over a}})^*  VV @VV   \Fka V\\
  \mathcal P_a(\R^N) @) \cal B_{k,a} ))     L^2(\R^N, \vartheta_{k,a}(x)dx)
  \end{CD}$$
 \end{prop}
\begin{proof}  
The identity \eqref{eqn:5.23} in $\mathfrak{sl}_2$ lifts to the
identity
\[
\mathcal{B}_{k,a} \circ \Omega_{k,a}(\Exp t \, {\bf h})
= \Omega_{k,a}(\Exp t \, {\bf k}) \circ \mathcal{B}_{k,a},
\]
and in particular
\[
\mathcal{B}_{k,a} \circ \Omega_{k,a}(\Exp \frac{\pi}{2i}{\bf h})
= \Omega_{k,a}(\Exp \frac{\pi}{2i}{\bf k}) \circ \mathcal{B}_{k,a},
\]
on $\mathcal{P}_a(\R^N)$ where the both-hand sides make sense.
In terms of the $(k,a)$-generalized Fourier transform $\Fka$ (see
\eqref{eqn:5.1}),
we get
\[
\mathcal{B}_{k,a} \circ 
\exp\Bigl(\frac{\pi i}{2} \frac{2\langle k\rangle+N+a-2}{a}\Bigr)
\, \Omega_{k,a} \Bigl(\Exp \frac{\pi}{2i}{\bf h}\Bigr)
= \Fka \circ \mathcal{B}_{k,a}.
\]
On the other hand,
we recall from \eqref{eqn:3.5} that
$$
\omega_{k,a}({\bf h})
= \frac{2}{a} \sum_{j=1}^N x_j\partial_j +
  \frac{N+2\langle k\rangle+a-2}{a},
$$
and therefore its lift to the group representation is given by
\begin{equation}\label{eqn:exph}
(\Omega_{k,a}(\Exp t\, {\bf h})f)(x)
= \exp\Bigl(\frac{N+2\langle k\rangle+a-2}{a} t\Bigr) f(e^{\frac{2t}{a}} x).
\end{equation}
Substituting $t=\frac{\pi}{2i}$,
we get
$\mathcal{B}_{k,a} \circ \Bigl(e^{-i{\pi\over a}}\Bigr)^* = \Fka \circ \mathcal{B}_{k,a}$.
This completes the proof of Proposition \ref{claim:5.10}
\end{proof}

When $k\equiv 0$ and $a=2,$ $\cal B_{0,2}$ coincides with the inverse of the 
Segal--Bargmann transform restricted to $\cal P(\R^N)=\cal P_{0,2}(\R^N)$ (cf. \cite[p. 40]{xfol}).  
We may think of $\cal B_{k,a}$ as a $(k,a)$-generalized  Segal--Bargmann transform. 
We are ready to prove the second statement of Theorem \ref{thm:5.11}.

\begin{proof}[Proof of Theorem \ref{thm:5.11} 2)]   
In view of Proposition \ref{claim:5.10}, we have $\Fka\circ \cal B_{k,a} (p)=\cal B_{k,a} \circ (e^{-i{\pi\over a}})^*(p).$ 
Since $(e^{-i\frac{\pi}{a}})^* p(x) = e^{-\frac{im\pi}{a}} p(x)$ for a
homogeneous polynomial of degree $m$,
we get
\[
\Fka \circ \mathcal{B}_{k,a}(p)
= e^{-\frac{im\pi}{a}} \mathcal{B}_{k,a}(p).
\]
Hence, the reproducing property \eqref{eqn:reprB} is proved.
 \end{proof}
\subsubsection{Application of Master Formula}
\hfill\break
As an immediate consequence of Master Formula 
(see Theorem \ref{thm:5.11}), we have:
\begin{coro}\label{coro:5.14}
 {\rm (Hecke type identity)} If in addition to the
assumption in Theorem \ref{thm:5.11} 2),  the polynomial   $p$ is $k$-harmonic of degree $m$, 
then \eqref{eqn:reprB} reads 
\begin{equation}\label{5.14} 
\Fka(e^{-\frac{\Vert \,\cdot\,\Vert^a}{a}} p) (\xi)= e^{-i{\pi\over a} m}  e^{-{1\over a} \Vert \xi\Vert^a} p(\xi).
\end{equation}
\end{coro}
Corollary \ref{coro:5.14} may be regarded as a Hecke type identity for
the $(k,a)$-generalized Fourier transform $\Fka.$ 
An alternative way to prove this identity would be to substitute $0$ for  $\ell$  in \eqref{5.6}.  

The    identity  \eqref{5.14} is a particular case of Theorem \ref{thm:5.15} below. 
For this, we will denote by 
$\text {\itbf H}_{a,\nu}$ the
 classical Hankel transform of one variable 
defined by
\begin{equation}\label{5.15}
\text {\itbf H}_{a,\nu}(\psi  )(s):=\int_0^\infty \psi  (r) \widetilde J_\nu \Big({2\over a}( {rs})^{a\over 2}\Big) r^{a(\nu+1)-1} dr,\end{equation} 
for a function $\psi$ defined on $\R_+$.
Here, $\widetilde{J}_\nu$ is the normalized Bessel function 
$\widetilde J_\nu(w)=\big({w\over 2}\big)^{-\nu} J_\nu(w) $  (see \eqref{eqn:4.15}).
Then the $(k,a)$-generalized Fourier transform $\Fka$ satisfies the
following identity:
\begin{thm}\label{thm:5.15}
 {\rm (Bochner type identity)} If $f \in (L^1 \cap
 L^2)(\R^N,\vartheta_{k,a}(x)dx)$ is of the form 
$f(x)= p(x) \psi  (\Vert x\Vert )$  for some $p\in \Hkm(\R^N)$ and a
 one-variable function $\psi$ on $\R_+$,
 then 
$$\Fka(f)(\xi) = a^{-({{2m+2\langle k\rangle+N-2}\over a})}e^{-i{\pi\over a} m}  p(\xi) \text { \itbf H}_{a,\frac{2m +2\langle k\rangle+N-2}{a}}(\psi  ) (\Vert \xi\Vert). $$ 
In particular, if $f $ is radial,  then $\Fka(f)$ is also radial.
\end{thm}
\begin{rem}\label{rem:5.16}
The original Bochner identity for the Euclidean Fourier transform corresponds to the case $a=2$ and $k
\equiv 0$.
For $a=2$ and $k>0$,
Theorem \ref{thm:5.15} corresponds to the Bochner identity for the
Dunkl transform which was proved in \cite{B}.
For $a=1$ and $k \equiv 0$
it is the Bochner identity for the Hankel-type transform on $\R^N$ (see \cite{KM}).
\end{rem}
\begin{proof}[Proof of Theorem \ref{thm:5.15}]
It follows from \eqref{eqn:4.8} that
\begin{align}
\Lambda_{k,a}^{(m)} \Bigl(r,s;\frac{\pi i}{2}\Bigr)
&= \exp \Bigl( -\frac{\pi i}{2} (\lambda_{k,a,m}+1)\Bigr)
   (rs)^{-\langle k\rangle-\frac{N}{2}+1} 
   J_{\lambda_{k,a,m}} \Bigl( \frac{2}{a}(rs)^{\frac{a}{2}}\Bigr)
\nonumber
\\
&= a^{-\lambda_{k,a,m}} (rs)^m 
   \exp \Bigl( -\frac{\pi i}{2} (\lambda_{k,a,m}+1)\Bigr)
   \widetilde{J}_{\lambda_{k,a,m}} \Bigl(\frac{2}{a} (rs)^{\frac{a}{2}}\Bigr).
\label{eqn:LawJ}
\end{align}
We set $\psi_m(r):=r^m\psi(r)$.
Since $p$ is homogeneous of degree $m$,
we have
\[
p \left( \frac{x}{\Vert x\Vert} \right) \psi_m (\Vert x\Vert)
= p(x) \psi(\Vert x\Vert).
\]
From the definition of the unitary operator
$\Omega_{k,a}^{(m)} (\gamma_z)$
(see \eqref{eqn:4.2}),
we get
\begin{align*}
\Omega_{k,a}(\gamma_z)f(x)
&= p(x) \Vert x\Vert^{-m} \Omega_{k,a}^{(m)} (\gamma_z) \psi_m (\Vert x\Vert)
\\
&= p(x) \Vert x\Vert^{-m} \int_0^\infty \Lambda_{k,a}^{(m)}
   (\Vert x\Vert,s;z) \psi_m(s) s^{2\langle k\rangle+N+a-3} ds.
\end{align*}
Substituting \eqref{eqn:LawJ} into the above formula with
$z = \frac{\pi i}{2}$,
we get
\begin{align*}
\Omega_{k,a}(\gamma_{\frac{\pi i}{2}}) f(x)
&= a^{-\lambda_{k,a,m}}
   \exp \Bigl( -\frac{\pi i}{2} (\lambda_{k,a,m}+1) \Bigr) p(x)
   \int_0^\infty \widetilde{J}_{\lambda_{k,a,m}}
   \left( \frac{2}{a} (\Vert x\Vert s)^{\frac{a}{2}} \right)
   \psi(s) s^{2m+2\langle k\rangle+N+a-3} ds
\\
&= a^{-\lambda_{k,a,m}}
   \exp \Bigl( -\frac{\pi i}{2} (\lambda_{k,a,m}+1) \Bigr) p(x)
   \text{\itbf H}_{a,\lambda_{k,a,m}} (\psi) (\Vert x\Vert).
\end{align*}
Now, Theorem \ref{thm:5.15} follows from \eqref{eqn:5.1}.
\end{proof}

\subsection{DAHA and $SL_2$-action}\label{subsec:5.6}
\hfill\break
In this subsection we discuss some link between the representation
$\Omega_{k,a}$ of $\widetilde{SL(2,\R)}$ in
the $a=2$ case and the (degenerate) rational DAHA (double affine Hecke algebra).
To be more precise,
we shall see that our representation $\Omega_{k,2}$ of
$\widetilde{SL(2,\R)}$  induces the representation of
$SL(2,\C)$  on the algebra generated by Dunkl's operators,
multiplication operators,
and the Coxeter group (see \eqref{eqn:HHgen} below). This induced action on the operators coincides essentially with a special case of the $SL(2,\Z)$-action discovered by
Cherednik \cite{Ch} and that of the $SL(2, \C)$-action by Etingof and Ginzburg \cite{EG}. Note that our approach (but not the result) is new in that we use our action on functions to derive the action on the operators in the Hecke algebra.  The authors are grateful to E. Opdam for bringing their attention to
this link.

We begin with an observation
that if $\Omega$ is a representation of a group $G$
on a vector space $W$ 
then we can define an automorphism of the associative algebra
$\End(W)$ by
\begin{equation}\label{eqn:conjG}
A \mapsto \Omega(g) A \Omega(g)^{-1},
\quad
g \in G.
\end{equation}

We shall consider this induced action for
$G = \widetilde{SL(2,\R)}$, $\Omega = \Omega_{k,2}$
(see Theorem \ref{thm:3.12}),
$W = $ the vector space consisting of appropriate functions on $\R^N$.

\begin{rem}\label{rem:EndW}
We do not specify the class of functions here.
Instead, we shall use the formula \eqref{eqn:conjG} to define
algebraically the $G$-action on a certain subspace of\/ $\End(W)$.
The point here is that the $G$-action on such a subspace will be
well-defined even when the group $G$ may not preserve $W$.
\end{rem}

We begin with a basic fact on Dunkl operators on $\R^N$.
For $\xi\in \R^N,$ we define the multiplication operator $M_\xi$ by 
$$
M_\xi f(x):=\langle \xi, x\rangle f(x).
$$

Choose an orthonormal basis $\xi_1,\ldots, \xi_N$ in $\R^N.$ 
As in Section \ref{subsec:DLap},
we will use the abbreviation $T_i(k)$ for Dunkl operators $T_{\xi_i}(k)$,
and $M_j$ for $M_{\xi_j}$.
Then we have 
the following commutation relations:
\begin{equation}\label{eqn:TM}
[T_{i}(k), M_{j}]=\delta_{ij} +2\sum_{\alpha\in \cal R^+} k_\alpha {{\langle \alpha, \xi_i\rangle 
\langle \alpha, \xi_j\rangle }\over {\Vert \alpha\Vert^2}}
r_\alpha,\quad \text{for any $1\leq i,j\leq N$}.
\end{equation}
Since the formula \eqref{eqn:TM} is symmetric with respect to $i$ and
$j$,
we have:
\begin{equation}\label{eqn:TM12}
[T_i(k), M_j]
= [T_j(k), M_i]
\quad\text{for any $1 \le i, j \le N$}.
\end{equation}
Furthermore, we have the following formulas:
\begin{lem}\label{lem:DMT}
Let $\xi \in \R^N$ and $s \in \C$.
\begin{enumerate}[\upshape 1)]
\item  
$[\Delta_k,M_\xi] = 2T_\xi(k)$.
\item  
$e^{s\Delta_k} M_\xi e^{-s\Delta_k}
 = M_\xi + 2s T_\xi(k)$.
\end{enumerate}
The first statement is due to Dunkl \cite[Proposition 2.2]{D1}, but we give its proof below for the reader's convenience. 
\end{lem}

\begin{proof}
1)
It is sufficient to prove the formula for
$\xi = \xi_j$ $(j=1,\dots,N)$.

By using \eqref{eqn:TM},
we have
\begin{align*}
[T_i^2(k), M_{\xi_j}]
&= T_i(k) [T_i(k), M_j] + [T_i(k), M_j] T_i(k)
\\
&= 2\delta_{ij} T_i(k) + 2\sum_{\alpha\in\mathcal{R}^+}
   k_\alpha \frac{\langle \alpha,\xi_i\rangle \langle \alpha,\xi_j \rangle}
                 {\Vert \alpha\Vert^2}
   (T_i(k) r_\alpha + r_\alpha T_i(k)).
\end{align*}
Summing them up over $i$, and using the following relations:
\begin{align*}
& \sum_{i=1}^N \langle \alpha,\xi_i \rangle T_i(k)
  = T_\alpha(k),
\\
& T_\alpha(k) r_\alpha + r_\alpha T_\alpha(k)
  = 0,
\quad\text{(see (D1) in Section \ref{subsec:2.1})},
\end{align*}
we get
\[
[\Delta_k,M_{\xi_j}] = 2T_j(k).
\]

2)
The second statement is straightforward from the first statement.
\end{proof}

Let us consider the induced action of $G$ on $\End(W)$ (see Remark
\ref{rem:EndW}). 
\begin{prop}\label{prop:SLEnd}
We fix a non-zero $\xi \in \R^N$.
\begin{enumerate}[\upshape 1)]
\item 
The induced action of\/ $\Omega_{k,2}$ by \eqref{eqn:conjG}
 preserves the two dimensional subspace
\[
\C_\xi^2 := \C M_\xi + \C T_\xi (k).
\]
\item 
The resulting representation of $\widetilde{SL(2,\R)}$ on $\C_\xi^2$
descends to $SL(2,\R)$,
and extends holomorphically to $SL(2,\C)$.
\item 
Via the basis $\{M_\xi,T_\xi (k) \}$,
the representation of $SL(2,\C)$ on $\C_\xi^2$ is given by
\begin{equation}\label{eqn:SLMT}
\varphi: SL(2,\C) \to GL_\C(\C_\xi^2),
\quad
\begin{pmatrix} A & B \\ C & D \end{pmatrix}
\mapsto
\begin{pmatrix} A & -iB \\ iC & D \end{pmatrix}.
\end{equation}
\end{enumerate}
\end{prop}

\begin{proof}
1)
Since $\widetilde{SL(2,\R)}$ is generated by
$\Exp(t\kern.03em{\bf e}^+)$ and $\Exp(t\kern.03em{\bf e}^-)$ $(t \in \R)$,
it is sufficient to prove that the subspace
$\C M_\xi + \C T_\xi (k)$ is stable by the induced action of these
generators.
In light of the formula \eqref{eqn:3.3},
we have
\[
\Omega_{k,2}(\Exp t\kern.03em{\bf e}^+)
=
\Exp(t\kern.03em\EE_{k,2}^+)
= \exp \left( \frac{it}{2} \Vert x\Vert^2 \right).
\]
Obviously, this action  commutes with the multiplication operator
$M_\xi$.
On the other hand,
applying \eqref{3.9} with $a=2$ and $\lambda=-\frac{it}{2}$,
we get
\[
\exp \left( \frac{it}{2} \Vert x\Vert^2 \right)
\circ T_\xi(k) \circ \exp \left( -\frac{it}{2} \Vert x\Vert^2 \right)
= T_\xi(k) - it M_\xi.
\]

Hence, $\Exp(t\kern.03em\EE_{k,2}^+)$ preserves the two-dimensional subspace
$\C M_\xi + \C T_\xi (k)$,
and its action is given by the following matrix form
\begin{equation}\label{eqn:Endt}
\Exp(t\kern.03em{\bf e}^+) \mapsto
\begin{pmatrix} 1 & -it \\ 0 & 1 \end{pmatrix}
\end{equation}
with respect to the basis $\{ M_\xi,T_\xi (k) \}$.

Next, we consider the action
\[
\Omega_{k,2}(\Exp(t\kern.03em{\bf e}^-)) 
=
\Exp(t\kern.03em\EE_{k,2}^-)
= \exp \left( \frac{it}{2} \Delta_k\right).
\]
Obviously, it commutes with Dunkl's operator $T_\xi (k)$.
On the other hand,
applying Lemma \ref{lem:DMT} 2) with $s = \frac{it}{2}$,
we get
\[
\exp \left( \frac{it}{2} \Delta_k \right)
\circ M_\xi \circ \exp \left( -\frac{it}{2} \Delta_k \right)
= M_\xi + it T_\xi (k).
\]
Hence, $\Exp (t\kern.03em\EE_{k,2}^-)$ also preserves the subspace
$\C M_\xi + \C T_\xi (k)$,
and its action is given as
\begin{equation}\label{eqn:End-}
\Exp (t\kern.03em{\bf e}^-)
\mapsto \begin{pmatrix} 1 & 0 \\ it & 1 \end{pmatrix}.
\end{equation}
Thus, we have proved the first statement.

2)
The center of $\widetilde{SL(2,\R)}$ consists of the elements
$\Exp(in\pi{\bf k})$ $(n \in \Z)$ (see \eqref{eqn:ZSL}).
Let us compute the action of $\Exp(t\kern.03em{\bf k})$ on
$\C_\xi^2 = \C M_\xi + \C T_\xi (k)$.
For this, we recall from \eqref{eqn:chk} that
\[
\Exp(t\kern.03em{\bf k}) = c_0 \Exp(t\kern.03em{\bf h}) c_0^{-1}
\quad\text{for $t \in \C$}.
\]
In view of the formulas \eqref{eqn:Endt} and \eqref{eqn:End-},
the element 
$c_0 
= \Exp i\begin{pmatrix} 0 & 1 \\ 0 & 0 \end{pmatrix}
  \Exp \frac{i}{2}\begin{pmatrix} 0 & 0 \\ 1 & 0 \end{pmatrix}
$
(see \eqref{eqn:c0}) acts on $\C_\xi^2$ as
\[
c_0 \mapsto \begin{pmatrix} 1 & i \\[.5ex] 0 & 1 \end{pmatrix}
 \begin{pmatrix} 1 & 0 \\[.5ex] \frac{i}{2} & 1 \end{pmatrix}
= \begin{pmatrix} \frac{1}{2} & 1 \\[.5ex] \frac{i}{2} & 1 \end{pmatrix}.
\]
It follows readily from the formula
\[
\left(\Omega_{k,2}(\Exp(t\kern.03em{\bf h}))f\right)(x)
= \exp \left( \frac{N+2\langle k\rangle}{2} t \right)
  f(e^t x)
\]
(see \eqref{eqn:exph}) that the action on
$\Exp(t\kern.03em{\bf h})$ on $\C^2_\xi$ is given by
\[
\Exp(t\kern.03em{\bf h}) \mapsto \begin{pmatrix} e^t & 0 \\ 0 & e^{-t} \end{pmatrix}.
\]
Therefore, $\Exp(t\kern.03em{\bf k})$ acts on $\C_\xi^2$ by the formula:
\[
\Exp(t\kern.03em{\bf k})
\mapsto \begin{pmatrix} 
\frac{1}{2} & i \\[.5ex] \frac{i}{2} & 1 \end{pmatrix}
        \begin{pmatrix} e^t & 0 \\[.5ex] 0 & e^{-t} \end{pmatrix}
        \begin{pmatrix} \frac{1}{2} & i \\[.5ex] \frac{i}{2} & 1 \end{pmatrix}^{-1}
= \begin{pmatrix} \cosh(t) & -i \sinh(t) \\[.5ex] 
                  i \sinh(t) & \phantom{-}\cosh(t) \end{pmatrix}.
\]
In particular,
if $t = in\pi$,
then
$\Exp(in\pi{\bf k})$ acts as $(-1)^n \operatorname{id}$ on $\C_\xi^2$.
Thus, the action of $\widetilde{SL(2,\R)}$ descends to 
$\widetilde{SL(2,\R)}/2\Z \simeq SL(2,\R)$.
Then, clearly, this two-dimensional representation
 extends holomorphically to $SL(2,\C)$.
Hence, the second statement is proved.

3) 
Since a representation of $SL(2,\C)$ is uniquely determined by the
generators $\Exp(t\kern.03em{\bf e^+})$ and $\Exp(t\kern.03em{\bf e}^-)$
$(t \in \C)$,
the third statement follows from \eqref{eqn:Endt} and \eqref{eqn:End-}.
\end{proof}

Let $H\!\!H$ be the algebra
generated by
\begin{equation}\label{eqn:HHgen}
\{ M_\xi, T_\xi : \xi \in \R^N \} \cup \CC,
\end{equation}
where $\CC$ is the Coxeter group.
Its defining relations are given by 
the commutativity of the Dunkl operators $T_\xi (k)$ (see (D2) in Section \ref{sec:2}), 
the commutativity of the multiplication operators $M_\xi$,
the commutation  relations \eqref{eqn:TM}, 
and the following $\CC$-equivariance: 
$$h \circ T_\xi (k)\circ h^{-1}= T_{h \xi}(k),\quad h\circ  M_\xi
\circ h^{-1}= M_{h \xi} ,\quad \text{for any $h\in \CC,\; \xi\in \R^N$},$$
see \cite{CM, EG}.

We recall from Proposition \ref{prop:SLEnd} 3) that
the matrix representation of the $SL(2,\C)$-action on 
$\C_\xi^2 = \C M_\xi + \C T_\xi (k)$ does not depend on
$\xi \in \R^N \setminus \{0\}$. 
Then, a simple computation relied on \eqref{eqn:TM12} yields
\begin{align*}
& [g \cdot T_i(k), g \cdot T_j(k)] = 0 = g \cdot [T_i(k), T_j(k)],
\\
& [g \cdot M_i, g \cdot M_j] = 0 = g \cdot [M_i, M_j],
\end{align*}
for any $1 \le i, j \le N$ and for any $g \in SL(2,\C)$.
Likewise, we get from \eqref{eqn:TM}
\[
[g \cdot T_i(k), g \cdot M_j] =    [T_i(k), M_j].
\]

Furthermore, the representation $\Omega_{k,a}$ of $\widetilde{SL(2,\R)}$
commutes with the action of the Coxeter group $\CC$.
Therefore, the action of $SL(2,\C)$ on $\C_\xi^2$ 
$(\xi \in \R^N \setminus \{0\})$ and the trivial action on the Coxeter
group $\CC$ extends to an automorphism of 
$H\!\!H$ because all the defining relations of $H\!\!H$ are preserved by $SL(2,\C)$.

Hence, we have proved:
\begin{thm}\label{thm:SLDAHA}
The representation $\Omega_{k,2}$ of $\widetilde{SL(2,\R)}$ induces
the above action of
$SL(2,\C)$ on the algebra $H\!\!H$ as automorphisms.
\end{thm}

\begin{rem}
The $SL(2,\C)$-action on the algebra $H\!\!H$ is essentially the same with the one given in \cite[Corollary 5.3]{EG}.
\end{rem}

\begin{rem}\label{rem:PSL}
As we saw in the proof of Proposition \ref{prop:SLEnd},
the center $\begin{pmatrix} -1 & 0 \\ 0 & -1 \end{pmatrix}$ of
$SL(2,\C)$ acts on $\C_\xi^2$ as $-\operatorname{id}$.
Therefore, $PSL(2,\C)$ acts on $H\!\!H$ as projective
automorphisms. 
\end{rem}

In order to compare the $SL(2,\Z)$-action on $H\!\!H$ defined by
Cherednik \cite{Ch} we consider the following automorphism of
$SL(2,\C)$:
\begin{equation}\label{eqn:iota}
\iota: SL(2,\C) \to SL(2,\C),
\quad
\begin{pmatrix} A & B \\ C & D \end{pmatrix}
\mapsto
\begin{pmatrix} A & -iB \\ iC & D \end{pmatrix},
\end{equation}
and twist the $SL(2,\C)$-action on $H\!\!H$ (see Theorem
\ref{thm:SLDAHA}) by $\iota$.
This means that the new action takes the form
\begin{equation}\label{eqn:twSL}
\varphi \circ \iota :
SL(2,\C) \to GL_\C(\C_\xi^2),
\quad
\begin{pmatrix} A & B \\ C & D \end{pmatrix}
\mapsto
\begin{pmatrix} A & -B \\ -C & D \end{pmatrix},
\end{equation}
on the generators
$\C_\xi^2 = \C M_\xi + \C T_\xi(k)$
(see \eqref{eqn:SLMT}).

We write $\tau_1$ and $\tau_2$ for the automorphisms of $H\!\!H$
corresponding to  the generators
$\begin{pmatrix} 1 & 1 \\ 0 & 1 \end{pmatrix}$
and
$\begin{pmatrix} 1 & 0 \\ 1 & 1 \end{pmatrix}$
of $SL(2,\Z)$.
Then, by \eqref{eqn:twSL},
$\tau_1$ and $\tau_2$ are given by
$$\tau_1:  M_\xi \mapsto M_\xi,\quad T_\xi\mapsto
T_\xi-M_\xi,\quad h\mapsto h,$$
$$\tau_2:  T_\xi \mapsto T_\xi,\quad M_\xi\mapsto M_\xi-T_\xi,
\quad h\mapsto h,$$ 
which coincide with the one given in \cite{Ch}.

Of particular importance in Cherednik \cite{Ch}
 is the
automorphism 
$$
\sigma:= \tau_1\tau_2^{-1}\tau_1=\tau_2^{-1}
\tau_1\tau_2^{-1},
$$ 
which corresponds to the action of
$\begin{pmatrix}0&1\\ -1&0\end{pmatrix}$.
The automorphism $\sigma$ is characterized by
$$
\sigma(T_\xi)= -M_\xi, 
\quad
\sigma( M_\xi)= T_\xi
\quad\text{and}\quad
\sigma (h)=h
\quad\text{for all $h\in \CC$}.
$$

{}From our view point,
these automorphisms on $H\!\!H$ can be obtained as the
conjugations of the action on the function space (see \eqref{eqn:conjG}).
In view of the formulas (see \eqref{eqn:iota}):
\begin{alignat*}{2}
& \iota \begin{pmatrix} 1 & 1 \\ 0 & 1 \end{pmatrix}
&& = \begin{pmatrix} 1 & -i \\ 0 & 1 \end{pmatrix},
\\
& \iota \begin{pmatrix} 1 & 0 \\ 1 & 1 \end{pmatrix}
&& = \begin{pmatrix} 1 & 0 \\ i & 1 \end{pmatrix},
\\
& \iota \begin{pmatrix} 0 & 1 \\ -1 & 0 \end{pmatrix}
&& = \begin{pmatrix} 0 & -i \\ -i & 0 \end{pmatrix}
= \Exp \left( \frac{\pi i}{2} {\bf h} \right) \gamma_{\frac{\pi i}{2}},
\end{alignat*}
we may interpret that $\tau_1$, $\tau_2$, and $\sigma$ are given by
the conjugations of
\begin{align*}
&\Omega_{k,2} (\Exp(-i {\bf e}^+))
= \exp \Bigl( \frac{1}{2} \Vert x\Vert^2 \Bigr),
\\
&\Omega_{k,2} (\Exp(i {\bf e}^-))
= \exp \Bigl( -\frac{1}{2} \Delta_k \Bigr),
\end{align*}
\begin{align}\label{eqn:Fre}
\cal F_{k,2} ^{re}
&:=e^{\frac12\Vert x\Vert^2} \circ \exp\Bigl(\frac{1}{2}\Delta_k\Bigr)\circ e^{\frac12\Vert x\Vert^2}= \exp\Bigl(\frac{1}{2}\Delta_k\Bigr)\circ e^{\frac12\Vert x\Vert^2} \circ \exp\Bigl(\frac{1}{2}\Delta_k\Bigr)
\\
&= \Omega_{k,2} \Bigl(\Exp \frac{\pi i}{2} {\bf h}\Bigr) 
   \Omega_{k,2}(\gamma_{\frac{\pi i}{2}}),
\nonumber
\end{align}
respectively.
Recalling the formulas:
\begin{align*}
&(\Omega_{k,2}\Bigl(\Exp \frac{\pi i}{2} {\bf h}\Bigr) f)(x)
= \exp \left( \frac{i\pi(N+2\langle k \rangle)}{4} \right) f(ix),
\\
&(\Omega_{k,2} ( \gamma_{\frac{\pi i}{2}} ) f)(x)
= \exp \left( -\frac{i\pi(N+2\langle k \rangle)}{4} \right) 
  (\mathcal{F}_{k,2}f)(x),
\end{align*}
we have
$$
\cal F_{k,2} ^{re} f(x)=\cal F_{k,2} f(ix).
$$ 
Hence, $\sigma$ may be interpreted as an algebraic version of the Dunkl
transform. 
We notice that the formula \eqref{eqn:Fre} fits well into 
 Master Formula \eqref{eqn:reprB} for $a=2,$ which we may rewrite as 
$$ c_{k,2}\int_{\R^N}e^{-{1\over 2}\Vert u\Vert^2} \big(\exp\Bigl(-{1\over {2 }} \Delta_k\Bigr) p\big)(u) B_{k,2}( ix,u) \prod_{\alpha\in \cal R}\vert \langle \alpha, u\rangle\vert^{k_\alpha} du=  e^{{1\over 2}\Vert x\Vert^2}p(x).$$

\subsection{The uncertainty inequality for the  transform $\Fka$} \hfill\\
The Heisenberg uncertainty principle may be formulated by means of 
the so-called Heisenberg  inequality for
the Euclidean Fourier transform on $\R$.
Loosely, the more a function is concentrated, the more its Fourier
transform is spread. We refer the reader to an excellent survey
\cite{Folland} for various mathematical aspects of the Heisenberg uncertainty principle. 
In this section we extend the Heisenberg inequality to the 
$(k,a)$-generalized Fourier transform $\Fka$ on $\R^N$.

Let  $\Vert \cdot \Vert_k$ be the $L^2$-norm with
respect to the measure $\vartheta_{k,a}(x)dx$ on $\R^N$ (see
\eqref{eqn:01}).
The  goal of this subsection is to prove the following
multiplicative inequality:
  \begin{thm}\label{thm:5.20}
 {\rm{ (Heisenberg type inequality)}} 
For all   $f\in L^2 (\R^N, \vartheta_{k,a}(x)dx)$ the $(k,a)$-generalized Fourier 
 transform $\Fka$  satisfies  
\begin{equation}\label{eqn:Heisen}
\big\Vert \Vert \cdot\Vert ^{a\over 2} f\big\Vert_k\;\big\Vert
\Vert \cdot\Vert ^{a\over 2} \Fka(f)\big\Vert_k\geq  \Big({{ 2\langle
k\rangle +N+a-2}\over 2} \Big) \Vert f\Vert_k^2.
\end{equation}
The equality holds if
and only if the function $f$ is of the form $f(x)=\lambda \exp(-c\Vert x\Vert^a)$ for some $\lambda \in \C$ and $c\in\R_+$.
\end{thm}
\begin{rem}\label{rem:5.21}
The inequality \eqref{eqn:Heisen} for $k\equiv 0$ and $a=2$ is the original
Heisenberg inequality for the Euclidean Fourier transform.
The inequality for $k>0$ and $a=2$ is the Heisenberg type inequality
for the Dunkl transform,
which was proved first by R\"{o}sler \cite{R2} and then by Shimeno \cite{Sh}. In physics terms we can think of the function $f(x)=\lambda \exp(-c\Vert x\Vert^a)$ where the equality holds in the above theorem as a ground state; indeed when $a=c=1,$ $N=3,$ and $k\equiv 0,$ it is exactly the wave function for the Hydrogen atom with the lowest energy. 
\end{rem}

In order to prove Theorem \ref{thm:5.20} we begin with the following
additive inequality:
\begin{lem}\label{claim:5.18}
\begin{enumerate}[\upshape (1)]
\item  
For all  $f\in L^2\big(\R^N, \vartheta_{k,a}(x)dx\big)$
\begin{equation}\label{5.16}
\big\Vert  \Vert \cdot\Vert^{a\over 2} f\big\Vert_k^2+ \big\Vert  \Vert \cdot\Vert^{a\over 2} \Fka(f)\big\Vert_k^2\geq  ({ 2\langle k\rangle +N+a-2} ) \Vert f\Vert_k^2. \end{equation}
\item  
The equality holds in \eqref{5.16} if and only if $f(x)$ is a scalar
multiple of\/ $\exp(-\frac{1}{a}\Vert x\Vert^a)$.
\end{enumerate}
\end{lem}
\begin{proof}
By Theorem \ref{thm:Fka}(3) and Theorem \ref{thm:Fkauni}(1),
we get
\begin{align*}
\left\Vert \Vert\cdot\Vert^{\frac{a}{2}} \Fka f\right\Vert_k^2
& = \langle\!\langle \Vert x\Vert^a \Fka f, \Fka f \rangle\!\rangle_k
\\
& = -\langle\!\langle \Fka(\Vert x\Vert^{2-a}\Delta_k f), \Fka f\rangle\!\rangle_k
\\
& = -\langle\!\langle\Vert x\Vert^{2-a}\Delta_k f, f\rangle\!\rangle_k.
\end{align*}
Hence, the left-hand side of \eqref{5.16} equals
\begin{equation}\label{eqn:afFf}
\langle\!\langle(\Vert x\Vert^a - \Vert x\Vert^{2-a}\Delta_k)f,f\rangle\!\rangle_k
= \langle\!\langle-\Delta_{k,a}f,f\rangle\!\rangle_k.
\end{equation}
It then follows from Corollary \ref{cor:kaLap} that the self-adjoint
operator $-\Delta_{k,a}$ has only discrete spectra,
of which the minimum is $2\langle k\rangle+N-2+a$.
Therefore, we  have proved
\[
\mbox{\eqref{eqn:afFf}}
\ge (2\langle k\rangle+N-2+a) \Vert f\Vert_k^2.
\]
Thus, the inequality \eqref{5.16} has been proved.
Further, the equality holds if and only if $f$ is an eigenfunction of
$-\Delta_{k,a}$ corresponding to the minimum eigenvalue
$2\langle k\rangle+N-2+a$, namely,
$f$ is a scalar multiple of
$\exp(-\frac{1}{a}\Vert x\Vert^a)$
(i.e.\ by putting $\ell=m=0$ in the formula \eqref{3.16} of
$\Phi_\ell^{(a)}(p,x)$).
Hence, Lemma \ref{claim:5.18} has been proved.
\end{proof}

\begin{proof}[Proof of Theorem \ref{thm:5.20}]
Now, for $c>0,$ we set $f_c(x):=f(cx).$ Using the fact that the density
$\vartheta_{k,a}$ is homogeneous of degree $2\langle k\rangle+a-2,$ we get
\[
\big\Vert  \Vert \cdot\Vert^{a\over 2} f_c\big\Vert_{k}^2=c^{-2\langle k\rangle-N-2a+2}  \big\Vert  \Vert \cdot\Vert^{a\over 2} f\big\Vert_{k}^2, 
\]
and 
 \[
\Vert f_c\Vert_k^2=c^{-2\langle k\rangle-N-a+2} \Vert f\Vert_k^2. 
\]
Furthermore, we lift the formula in Theorem \ref{thm:Bka}(1) to the
formula
\[
(\Fka f_c)(x) = c^{-(N+2\langle k\rangle+a-2)}
(\Fka f)\Bigl(\frac{x}{c}\Bigr),
\]
from which we get
\[
\big\Vert  \Vert \cdot\Vert^{a\over 2} \Fka(f_c)\big\Vert_{k}^2=c^{-2\langle k\rangle-N+2}  \big \Vert  \Vert \cdot\Vert^{a\over 2} \Fka(f)\big\Vert_{k}^2.  \]
Thus, if we substitute $f_c$ for $f$ in Lemma \ref{claim:5.18}, we obtain  $$c^{-a}  \big\Vert   \Vert \cdot\Vert^{a\over 2} f\big\Vert_{k}^2+c^a  \big\Vert  \Vert \cdot\Vert^{a\over 2} \Fka(f)\big\Vert_{k}^2\geq 
  ({ 2\langle k\rangle +N+a-2} ) \Vert f\Vert_k^2.$$ Obviously  the minimum  value of the left-hand side (as a function of $c\in \R_+$) is 
 $$ 2\big\Vert  \Vert \cdot\Vert ^{a\over 2}f\big\Vert_k\;\big\Vert  \Vert \cdot\Vert ^{a\over 2}\Fka(f)\big\Vert_k.$$
Hence, Theorem \ref{thm:5.20} has been proved.
\end{proof}

\medbreak

$\langle$Acknowledgement$\rangle$

It is a great pleasure to thank C. Dunkl and
E. Opdam for their valuable comments on a very preliminary version of
this article.
The second author is partially supported by Grant-in-Aid for Scientific 
Research (B) (18340037, 22340026),
Japan Society for the Promotion
of Science, Harvard University, Max Planck Institut f\"ur Mathematik,
and the Alexander Humboldt Foundation of Germany.


\renewcommand\indexname{List of Symbols}
\begin{theindex}

  \item $\Gamma(W)$, \main{32}
  \item $\Delta_k$, \main{14}
  \item $\Delta_{k,a}$, \main{5}
  \item $\Lambda_{k,a}^{(m)}(r,s;z)$, \main{35}
  \item $\Phi _{ \ell }^{(a)}(p,x)$, \main{23}
  \item $\widetilde {\Phi }_{\ell }^{(a)}(p,x)$, \main{26}
  \item $\Omega_{k,a}$, \main{30}
  \item $\gamma_z$, \main{27}, \main{33}
  \item $\vartheta_k(x)$, \main{12}
  \item $\vartheta_{k,a}(x)$, \main{5}
  \item $\lambda_{k,a,m}$, \main{19}
  \item $\mu_x^k$, \main{13}
  \item $\pi(\lambda)$, \main{28}
  \item $\pi_K(\lambda)$, \main{28}
  \item $\omega_{k,a}$, \main{18}

  \indexspace

  \item \par\indexspace$\langle\kern-.2em\langle \ ,\ \rangle\kern-.2em\rangle _k$, 
		\main{29}

  \indexspace

  \item $B_{k,a}(\xi,x)$, \main{54}
  \item $\cal B_{k,a}$, \main{60}

  \indexspace

  \item $C(G)$, \main{28}
  \item $C_m^\nu(t)$, \main{40}
  \item $\C^+$, \main{33}
  \item $\C^{++}$, \main{33}
  \item $\CC$, \hyperpage{5}, \main{11}
  \item $\mathcal{C}_{\nu,m}$, \main{42}
  \item $c_{k,a}$, \main{45}

  \indexspace

  \item $d_k$, \main{15}
  \item $\mathcal {D}_k$, \main{7}

  \indexspace

  \item $E$, \main{15}
  \item $\EE^-_{k,a}$, \main{17}
  \item $\EE_{k,a}^+$, \main{17}
  \item $\widetilde {\mathbb {E}}_{k,a}^+$, \main{19}
  \item $\widetilde {\mathbb {E}}_{k,a}^-$, \main{19}
  \item ${\bf e}^+$, \main{17}
  \item ${\bf e}^-$, \main{17}

  \indexspace

  \item $\Fka$, \hyperpage{7}, \main{51}
  \item $f_{\ell,m}^{(a)}(r)$, \main{24}

  \indexspace

  \item $\HH_{k,a}$, \main{17}
  \item $\widetilde {\mathbb {H}}_{k,a}$, \main{19}
  \item $\mathcal{H}_k$, \hyperpage{7}, \main{52}
  \item $\Hkm(\R^N)$, \main{15}
  \item $h_{k,a}(r,s;z;t)$, \main{48}
  \item ${\bf h}$, \main{17}

  \indexspace

  \item $I_\lambda(w)$, \main{36}
  \item $\widetilde I_\lambda (w)$, \main{36}
  \item $\mathcal{I}_{k,a}(z)$, \main{5}

  \indexspace

  \item $\widetilde{J}_\nu(w)$, \main{57}

  \indexspace

  \item $ \langle k\rangle$, \main{12}
  \item ${\bf  k}$, \main{18}

  \indexspace

  \item $L_\ell^{(\lambda)}(t)$, \main{21}

  \indexspace

  \item ${\bf  n}^+$, \main{18}
  \item ${\bf  n}^-$, \main{18}

  \indexspace

  \item $P_{k, {m} }(\omega, \eta)$, \main{47}

  \indexspace

  \item $\cal R$, \hyperpage{5}, \main{11}

  \indexspace

  \item $T_\xi(k)$, \main{11}

  \indexspace

  \item $V_k$, \main{12}

  \indexspace

  \item $W_{k,a}(\R^N)$, \main{23}

\end{theindex}


\begin{thebibliography}{99}
\bibitem{AAR}
G. Andrews, R. Askey and R. Roy, \textit{Special Functions}, Cambridge, 1999.
\bibitem{B} S. Ben Said, {\it On the integrability of a representation of $\s\l(2,\R)$},  J. Funct. Anal. {\bf 250} (2007), 249--264.
\bibitem{BKO-CRAS}
S. Ben Said, T. Kobayashi and B. {\O}rsted,
\href{http://dx.doi.org/10.1016/j.crma.2009.07.015}{\textit{Generalized Fourier Transforms $\mathcal{F}_{k,a}$}},
C. R. Acad. Sci. Paris, Ser. I \textbf{347} (2009), 1119-1124.
\bibitem{BO} S. Ben Said and B. {\O}rsted, {\it Segal--Bargmann transforms associated with finite Coxeter groups},   Math. Ann.  {\bf 334}  (2006),   281--323.
\bibitem{Ch}  I. Cherednik, {\it Macdonald's evaluation conjectures and difference Fourier transform},   Invent. Math.  {\bf 122}  (1995),   119--145.
\bibitem{CM} I. Cherednik and Y.  Markov,  {\it Hankel transform via double Hecke algebra},  Iwahori--Hecke algebras and their Representation Theory (Martina-Franca, 1999), M. W. Baldoni and D. Barbasch (eds.), Lecture Notes in Math. {\bf 1804}, Springer-Verlag, Berlin, 2002, 1--25.
\bibitem{DB}
L. Debnath and D. Bhatta,
\textit{Integral transforms and their applications},
 Second Edition, 
Chapman \& Hall/CRC, Boca Raton, FL, 2007. 
\bibitem{D88}
C. F. Dunkl,
\textit{Reflection groups and orthogonal polynomials on the sphere,}
Math. Z. \textbf{197} (1988), 33--60.
\bibitem{D1} C. F. Dunkl, {\it Differential-difference operators
associated to reflection groups,} Trans. Amer. Math. Soc. {\bf
311} (1989), 167--183.
\bibitem{D2} C. F. Dunkl,  {\it Integral kernels with reflection group invariance,} 
Canad. J. Math. {\bf 43} (1991), 1213--1227.
\bibitem{D3} C. F. Dunkl, {\it Hankel transforms associated to finite reflection groups,}   Proc. of the special session on hypergeometric functions on domains of positivity, Jack polynomials and applications. Proceedings, Tampa 1991, Contemp. Math. {\bf 138} (1992), 123--138.
\bibitem{D03} C. F. Dunkl, {\it A Laguerre polynomial orthogonality and the hydrogen atom},  Anal. Appl. (Singap.)  {\bf 1}  (2003),    177--188.
\bibitem{D08} C. F. Dunkl, 
\href{http://dx.doi.org/10.1007/s11537-008-0819-3}{\textit{Reflection groups in analysis and applications}},
Japan. J. Math. (3rd ser.) {\bf 3} (2008), 215--246.
\bibitem{DX} C. F. Dunkl and Y. Xu, {\it Orthogonal Polynomials of
Several Variables}, Cambridge University Press, Cambridge, 2001.
\bibitem{DJO} C. F. Dunkl, M. de Jeu and E. Opdam, {\it Singular polynomials for finite reflection groups,} Trans. Amer. 
Math. Soc. {\bf 346} (1994), 237--256.
\bibitem{E}
P. Etingof,
\href{http://arxiv.org/abs/0903.5084}{\textit{A uniform proof of the Macdonald--Mehta--Opdam identity for
finite Coxeter groups}}, 
arXiv 0903.5084 
\bibitem{EG} P. Etingof and V. Ginzburg, {\it  Symplectic reflection algebras, Calogero--Moser space, and deformed Harish-Chandra homomorphism},  Invent. Math.  {\bf 147}  (2002), 243--348.
\bibitem{xfol} G. B. Folland, Harmonic Analysis in Phase Space, Ann. of Math. Stud. {\bf 122}, Princeton University Press, Princeton, NJ, 1989.
\bibitem{Folland} G. B. Folland and A. Sitaram,  {\it The uncertainty principle: a mathematical survey},   J. Fourier Anal. Appl.  {\bf 3}  (1997),   207--238.
\bibitem{GY} L. Gallardo and M.  Yor,  
\href{http://dx.doi.org/10.1214/009117906000000133}{\textit{A chaotic representation property of the multidimensional Dunkl processes}}  Ann. Probab.  {\bf 34}  (2006),   1530--1549.
\bibitem{GG} I. M. Gelfand and S. G. Gindikin, {\it Complex manifolds whose spanning trees are real semisimple Lie groups, and analytic discrete series of representations},  Funct. Anal. Appl.  {\bf 11} (1977),  19--27.
\bibitem{GR} I. S. Gradshteyn and I. M. Ryzhik, {\it  Table of Integrals, Series, and Products},  Corrected and enlarged edition edited by Alan Jeffrey. Incorporating the fourth edition edited by Yu. V. Geronimus  and M. Yu. Tseytlin. Translated from the Russian. Academic Press [Harcourt Brace Jovanovich, Publishers], New York-London-Toronto, Ont., 1980.
\bibitem{H} G. J. Heckman, {\it A remark on the Dunkl
differential-difference operators,} Harmonic Analysis on Reductive
Groups, (eds. W. Barker and P. Sally), Progr. Math. {\bf
101}, Birkh\"auser, 1991, 181--191.
\bibitem{Helg} S. Helgason, {\it Groups and Geometric Analysis,} Academic Press, New York, 1984.
\bibitem{Hy} J. E. Humphreys, {\it Reflection Groups and Coxeter Groups}, 
Cambridge Stud. Adv. Math. {\bf 29}, Cambridge 
University Press, Cambridge, 1990.
\bibitem{HN} J. Hilgert and K-H. Neeb, {\it Lie Semigroups and their Applications},  Lecture Notes in Math. {\bf 1552}, Springer-Verlag, Berlin, 1993.
\bibitem{HN2} J. Hilgert and K-H. Neeb,  {\it Positive definite spherical functions on Olshanski domains},  Trans. Amer. Math. Soc.  {\bf 352}  (2000),   1345--1380.
\bibitem{H1} R. Howe, {\it On the role of the Heisenberg group in harmonic analysis},   Bull. Amer. Math. Soc. (N.S.)  {\bf 3}  (1980),   821--843.
\bibitem{H2} R. Howe, {\it The oscillator semigroup},  The Mathematical Heritage of Hermann Weyl (Durham, NC, 1987), R. O. Wells, Jr. (ed.), Proc. Sympos. Pure Math. {\bf 48}, Amer. Math. Soc., Providence, RI, 1988, 61--132.
\bibitem{HT} R. Howe and E.-C. Tan, {\it Nonabelian Harmonic Analysis. Applications of ${\rm SL}(2,R)$},  Universitext, Springer-Verlag, New York, 1992. 
\bibitem{dJ} M. F. E. de Jeu,  {\it The Dunkl transform},  Invent. Math.  {\bf 113}  (1993),    147--162.
\bibitem{KN}
A. W. Knapp and D. A. Vogan, Jr.,
\textit{Cohomological Induction and Unitary Representations},
Princeton University Press, 1995. 
\bibitem{Kdisc}T. Kobayashi,
\textit{Discrete decomposability of the restriction of $A_\mathfrak q(\lambda)$
 with respect to reductive subgroups and its applications}, 
{\href{http://dx.doi.org/10.1007/BF01232239}{%
Part I}}. Invent. Math. {\bf 117} (1994), 181--205;
{\href{http://dx.doi.org/10.2307/120963}{
Part II}}. Ann. of Math. \textbf{147} (1998), 709--729;
{\href{http://dx.doi.org/10.1007/s002220050203}{Part III}}.
Invent. Math. {\bf 131} (1998), 229--256.
\bibitem{Kaspm97}
T. Kobayashi,
\textit{Discretely decomposable restriction of unitary representations of
reductive Lie groups, ---examples and conjectures},
Advanced Studies in Pure Mathematics \textbf{26} (2000), 99--127.
\bibitem{KM05} T. Kobayashi and G. Mano, 
\href{http://www.springerlink.com/content/pt65717043175035/fulltext.pdf}{%
\textit{Integral formulas for the minimal representations for $O(p, 2)$}}, 
Acta Appl. Math. {\bf 86} (2005), 103--113.
\bibitem{KM} T. Kobayashi and G. Mano, {\it The inversion formula and
holomorphic extension of the minimal representation of the conformal
group},   Harmonic Analysis, Group Representations, Automorphic Forms
and Invariant Theory: In honor of Roger Howe, (eds. J. S. Li,
E. C. Tan, N. Wallach and C. B. Zhu), World Scientific,
2007, 159--223
(cf. \href{http://uk.arxiv.org/abs/math.DG/0607007}{math.RT/0607007)}).
\bibitem{KM-intopq} T. Kobayashi and G. Mano, 
{%
\textit{Integral formula of the unitary inversion operator for the minimal representation of $O(p,q)$}},
\href{http://projecteuclid.org/euclid.pja/1176126886}{ 
Proc. Japan Acad. Ser. A} {\bf 83} (2007), 27--31.
\bibitem{KM-intopq-ams} T. Kobayashi and G. Mano, 
\textit{The Schr\"odinger model for the minimal representation of the
indefinite orthogonal group $O(p,q)$}, vi + 132pp.  the
Mem. Amer. Math. Soc. vol. 212, 
{\href{http://dx.doi.org/10.1090/S0065-9266-2011-00592-7}{no. 1000}} (2011), cf.  
{\href{http://uk.arxiv.org/abs/0712.1769}{arXiv:0712.1769}}. 
\bibitem{KO} T. Kobayashi and B. {\O}rsted, {\it 
\href{http://dx.doi.org/10.1016/S0001-8708(03)00013-6}{%
Analysis on the minimal representation of $ O(p,q)$. II. Branching laws}},   
Adv. Math. {\bf 180} (2003),   513--550.
\bibitem{Kos} B. Kostant, {\it On Laguerre polynomials, Bessel functions, and Hankel transform and a series in the unitary dual of the simply-connected covering group of $SL(2, \mathbb R)$,} Represent. Theory {\bf 4} (2000), 181--224.
\bibitem{L} S. Lang,  $SL_2(\mathbb R)$,  Reprint of the 1975 edition. Graduate Texts in Mathematics, 105. Springer-Verlag, New York, 1985.
\bibitem{bigMac} I. G. Macdonald, {\it Some conjectures for root systems,} SIAM J. Math. Anal. 
{\bf 13} (1982), 988--1007.
\bibitem{Mano} G. Mano, {\it A continuous family of unitary representations with two hidden symmetries---an example,} RIMS Kokyuroku Bessatsu, Representation Theory and Analysis on Homogeneous Spaces (Conference Proceedings at RIMS, 2006 August), H. Sekiguchi (ed.), 2008, 137--144 (in Japanese).
\bibitem{Nels} E. Nelson, {\it Analytic vectors,} Ann. of Math. {\bf 70} (1959), 572--615.
\bibitem{Op} E. Opdam, {\it Dunkl operators, Bessel functions and discriminant of a finite Coxeter group,} Comp. Math. {\bf 85} (1993), 333--373.
\bibitem{Ols} G. Olshanski, {\it Invariant cones in Lie algebras, Lie semigroups and the holomorphic discrete series}, 
Funct. Anal. Appl. {\bf 15} (1981),   275--285.
\bibitem{RR} R. Ranga Rao,  {\it Unitary representations defined by boundary conditions---the case of ${\s}{\l}(2, \R)$},   Acta Math.  {\bf 139}  (1977),  185--216.
\bibitem{OSS} B. {\O}rsted,  P. Somberg and V. Soucek,  {\it The Howe duality for the Dunkl version of the Dirac operator},  Adv. Appl. Clifford Algebr. {\bf 19} (2009),   403--415. 
\bibitem{R1} M. R\"{o}sler, {\it Positivity of Dunkl's intertwining operator}, Duke Math. J. {\bf 98} (1999), 445--463.
\bibitem{R2} M. R\"{o}sler, \textit{An uncertainty principle for the Dunkl transform,}   Bull. Austral. Math. Soc.  {\bf 59}  (1999),    353--360.
\bibitem{RV} M. R\" {o}sler and  M. Voit, {\it Markov processes related with Dunkl operators},  Adv. in Appl. Math.  \textbf{21}  (1998),    575--643.
\bibitem{Sato}
M. Sato,  \textit{Theory of hyperfunctions I},
  J. Fac. Sci. Univ. Tokyo, \textbf{8}  (1959), 
139-193.
\bibitem{Sh} N. Shimeno, {\it A note on the uncertainty principle for the Dunkl transform}, 
J. Math. Sci. Univ. Tokyo {\bf 8} (2001),  33--42.
\bibitem{Sta} R. J. Stanton, {\it Analytic extension of the holomorphic discrete series},   Amer. J. Math.  {\bf 108}  (1986),  1411--1424.
\bibitem{TX05}
\bibitem{T} K. Trim\`eche, {\it Paley--Wiener theorems for the
Dunkl transform and Dunkl translation operators,} Integral Transforms Spec. Funct. {\bf 13} (2002), 17--38.
\bibitem{V} M. Vergne, {\it On Rossmann's character formula for discrete series},  Invent. Math.  {\bf 54}  (1979),   11--14.
\bibitem{Velin} 
N. Ja. Vilenkin, {\it Special Functions and the Theory of Group Representations},
Translated from the Russian by V. N. Singh, 
Transl. Math. Monogr. {\bf 22}, Amer. Math. Soc., Providence, RI, 1968.
\bibitem{Weil} A. Weil, {\it Sur certains groupes d'op\`erateurs unitaires,}  Acta Math. {\bf 111} (1964), 143--211.
\bibitem{W} Z. X. Wang,  and D. R. Guo,   {\it Special Functions}, Translated from the Chinese by Guo and X. J. Xia, World Scientific Publishing Co., Inc., Teaneck, NJ, 1989.
\bibitem{X1} Y. Xu,  {\it Funk--Hecke formula for orthogonal polynomials on spheres and on balls}, Bull. London Math. Soc. {\bf 32} (2000),  447--457.


\end{thebibliography}
\end{document}